\documentclass[a4paper,11pt,leqno]{article}

\usepackage[utf8]{inputenc}
\usepackage[T1]{fontenc} 
\usepackage[english]{babel}
\usepackage{verbatim}
\usepackage{calligra}
\usepackage{enumerate}
\usepackage{dsfont}
\usepackage{amsfonts}
\usepackage{amsmath}
\usepackage{amsthm}
\usepackage{amssymb}
\usepackage{mathtools}
\usepackage{mathrsfs}
\usepackage{color}
\usepackage{graphicx}
\usepackage{bm}
  \usepackage[all]{xy}
\usepackage{hyperref}
\usepackage{comment}
\hypersetup{
    colorlinks=true,
    linkcolor=black,
    filecolor=black,      
    urlcolor=black,
}

\usepackage{tikz}

\usepackage{fullpage}

\usepackage{graphicx}		
\usepackage[hypcap=false]{caption}

\theoremstyle{plain}
\newtheorem{defi}{Definition}
\newtheorem{thm}[defi]{Theorem}
\newtheorem{prop}[defi]{Proposition}

\theoremstyle{definition}

\newtheorem{rmk}[defi]{Remark}

\newcommand{\mc}{\mathcal}

\newcommand{\what}{\widehat}

\newcommand{\D}{{\rm div}\,}

\newcommand{\R}{\mathbb{R}}

\newcommand{\N}{\mathbb{N}}
\newcommand{\Z}{\mathbb{Z}}

\newcommand{\curl}{{\rm curl}\,}

\newcommand{\dx}{ \, {\rm d} x}
\newcommand{\dt}{ \, {\rm d} t}
\newcommand{\ds}{ \, {\rm d} s}
\newcommand{\dy}{ \, {\rm d} y}

\allowdisplaybreaks

\title{\textsc{Critical Well-Posedness of a Euler-Chemotaxis System IN 2D}}
\date{\today}
\author{Dimitri \textsc{Cobb}\\
\textsl{Université Paris-Est Créteil -- Laboratoire d'Analyse et de Mathématiques Appliquées} \\
\url{dimitri-alexandre.cobb@u-pec.fr}}

\begin{document}

\maketitle

\begin{abstract}
    This article is about well-posedness results for a 2D Euler-Chemotaxis-Consumption system with initial data in scale critical spaces, namely existence, uniqueness and critical initial data to solution continuity. This problem is related to the well-known family of Keller-Segel equations, although it bears major differences which will be crucial for our study.

    In most recent results, solutions to such Chemotaxis-Consumption systems have used a free-energy balance law to obtain well-posedness. This has the disadvantage of requiring initial data to be positive almost everywhere. Our main goal, in addition to sharp well-posedness, is to show that this positivity condition is not needed for existence and uniqueness. In fact, we will never use the concept of free-energy in our proof.

    Our argument revolves around an understanding of the scaling of the different parts of the equation, and especially the subcritical interaction terms. This will allow us to upgrade basic low regularity estimates (\textit{e.g.} mass conservation) into critical bounds. In the process, we find out that closing $L^2$ estimates leads to the Besov space $\dot B^0_{\infty, 2}$ as the most natural choice of regularity for the initial chemo-attractant concentration. Extensive use of real interpolation also plays a key role for deriving Chemin-Lerner bounds for the solutions.
\end{abstract}

{\footnotesize \textbf{Key words:} Chemotaxis, PDE coupling, Critical Spaces, Besov Spaces, Real Interpolation

\textbf{2020 Mathematics Subject Classification:} 35M31 (main), 76Z10, 35K58, 35Q31, 76B03}



\section{Introduction}

The purpose of this article is proving the existence and uniqueness of global in time solutions of the following Euler-chemotaxis system of equations. The corresponding notation are introduced immediately below.
\begin{equation}\label{ieq:EC}
    \begin{cases}
        \partial_t u + (u \cdot \nabla) u + \nabla \pi = ng\\
        \partial_t n - \Delta n + u \cdot \nabla n + n^3 - n + \D(n \nabla c) = 0\\
        \partial_t c - \Delta c + u \cdot \nabla c \pm nc = 0\\
        \D(u) = 0
    \end{cases}
\end{equation}
The equations are set on the whole 2D Euclidean space $x \in \R^2$. This system models the motion of mobile micro-organisms (micro-swimmers) in an ideal fluid, in response a chemical stimulus. This process, called chemotaxis, helps the micro-swimmers to find nutrients which they then consume, or to aggregate in large structures as they follow a chemical signal designed to attract others of their kind.

\medskip

Let us define the notation used above, and give summary explanations on the structure of the system (see below in the introduction for more details). In the equations above, $n(t,x) \geq 0$ is a concentration of a micro-organism species, which reproduces according to a ``logistic'' law. The organisms are suspended in a fluid, whose velocity field is $u(t,x) \in \R^2$ and pressure field is $\pi(t,x)$. The fluid is in turn affected by the presence of the micro-organism by gravity: above $g = (0, -1)$ is the constant gravity vector field. In addition to being advected (and subject to diffusion) by the fluid, the particles undergo a process of chemotaxis: motion due to the presence of a chemical stimulus, whose concentration is denoted by $c(t,x) \geq 0$. Through the particle flux $n \nabla c$, the micro-organisms migrate towards regions of high chemical concentration. Finally, the chemical is advected by the fluid, subject to diffusion, and is consumed (or produced) by the micro-organisms, through the term $\pm nc$.

\medskip

In the remainder of this introduction, we will give a quick overview of chemotaxis systems to provide context for our result. We will then discuss scaling matters, before stating and explaining our main result, Theorem \ref{it:main}. This will be followed by a summary of the arguments used in the proof.

\subsubsection*{AI Statement}

No AI was used, either in the process of finding ideas and proofs, performing any computations, or writing the \LaTeX \; file.

\subsection{Context on Chemotaxis}

In this Subsection, we give some context on chemotaxis systems: where they come from, what are there interesting qualitative behaviors, and give a brief review of the literature. 

\subsubsection{Keller-Segel Equations}

To introduce chemotaxis, we start by Keller-Segel equations. This is a family of PDEs whose simplest representative is the following parabolic/elliptic system
\begin{equation}\label{ieq:KellerSegel}
    \begin{cases}
        \partial_t n - \Delta n + \D (n \nabla c) \\
        - \Delta c = n.
    \end{cases}
\end{equation}
In the above, $n(t,x) \geq 0$ is a density of mobile micro-organisms, (slime-mold, micro-swimmers, cancer cells, \textit{etc.}), whose motion is due to natural diffusion (Fick's law) and the response to an attractive chemical signal of density $c(t,x) \geq 0$. The first equation in \eqref{ieq:KellerSegel} was derived by Keller and Segel \cite{KS1970, KS1971} to study slime mold aggregation. The second equation above introduces non-linearity in the equations: the chemo-attractant $c(t,x)$ is produced by the micro-organisms, so that they are attracted to each other, and instantly diffuses across the whole space $\R^d$. Therefore, system \eqref{ieq:KellerSegel} should be seen as a competition between two phenomena, a dissipative diffusion process, and an attractive chemotaxis motion.

\medskip

The 2D case, which is the focus of this article, is especially interesting. In two dimensions, both dissipation and chemotaxis act on the same scale as mass conservation. Observe that the total micro-organism mass is a constant of the evolution:
\begin{equation}\label{ieq:MassCons}
    \int n(t,x) \dx = \int n_0 (x) \dx := M.
\end{equation}
Then, examine the evolution of the second momentum of $n$. A quick computation\footnote{Replace c = $ (- \Delta)^{-1}n$ by its expression using a fundamental solution, and note that $\frac{x \cdot (x - y)}{|x - y|^2} = 1 - \frac{y \cdot (y - x)}{|y - x|^2}$.} gives
\begin{equation*}
    \frac{\rm d}{\dt} \int |x|^2 n(t,x) \dx = 4M \left( 1 - \frac{M}{8\pi} \right).
\end{equation*}
This equations shows that there are two broad types of behavior. If $M > 8 \pi$, then the momentum vanishes at some time $T > 0$, so that the solution $n(t,x)$ becomes singular in finite time. If $M \leq 8 \pi$, this does not happen, and solutions are global. We owe this remark to \cite{DLM}.

\medskip

There is massive literature on Keller-Segel equations, of which \eqref{ieq:KellerSegel} is but one example, so we cannot hope to do justice to the work done on these systems. We will be content to give a few references linked to our interests.

\medskip

Well-posedness (for a slightly different parabolic/elliptic problem) has been studied in \cite{Biler1995, Biler1998, KS2008, KSW2012, Iwabuchi2011, Iwabuchi2014}. Among these, the article of Kozono and Sugiyama \cite{KS2008} is especially important, as it establishes well-posedness results for initial data $n_0 \in L^1(\R^2)$, which is a critical space in 2D with regards to the scaling of \eqref{ieq:KellerSegel}. Indeed, if $(n, c)$ is a solution, then, for all $\lambda > 0$,
\begin{equation}\label{ieq:KellerSegelScaling}
    n_\lambda(t, x) = \lambda^2 n(\lambda^2 t, \lambda x) \qquad \text{and} \qquad c_\lambda(tx) = c(\lambda^2 t, \lambda x)
\end{equation}
also is a solution. The results of Iwabuchi \cite{Iwabuchi2011, Iwabuchi2014} are also worth singling out, as they provide well-posedness in critical Besov spaces $n_0 \in \dot B^{-1-d/p}_{p, \infty}$ for $p < + \infty$, which are equally critical and satisfy, again in 2D,
\begin{equation*}
    L^1(\R^2) \hookrightarrow \dot B^0_{1, \infty} (\R^2) \subset \dot B^{-2 + \frac{2}{p}}_{p, \infty}(\R^2).
\end{equation*}

\medskip

Study of the blowup was done in \cite{JL1992, Nagai1995, BKLN2006, BCM2008}. In an important article \cite{HV1996}, Herrero and Vel\'azquez provide a description of the singularity: the solution $n$ converges to the sum of a Dirac mass with an absolutely continuous (with respect to $\dx$) background. This process is called chemotactic collapse, and represents aggregation of the micro-organisms to form a macroscopic structure. This has prompted the study of measure-valued (weak) solutions by Senba and Suzuki \cite{SS2002}.

\medskip

When the $c$ equation above is replaced by a parabolic equation, we have a parabolic/parabolic Keller-Segel system (not unlike our PDE \eqref{ieq:EC}), which may feature very different properties: global well-posedness when $M < 8 \pi$, and possible but not automatic blowup when $M > 8 \pi$. Concerning well-posedness and blow-up, we refer to \cite{CC2008, Mizoguchi2013, BGK2015, HV1996b, SS2001, Mizoguchi2020, LemRieu}. Characterization of the blowup in terms of chemotactic collapse is done in \cite{HV1997, NSS2000}.

\medskip

Finally, since blowup (and therefore chemotactic collapse) in the parabolic/parabolic system may or may not occur beyond the critical mass $M > 8 \pi$, it becomes interesting to find other mechanisms which might prevent blowup. One that is especially relevant to our equation is the introduction of a reproduction term in the micro-organism density equation, which becomes
\begin{equation*}
    \partial_t n - \Delta n = \kappa n - \lambda n^2 - \D(n \nabla c),
\end{equation*}
where $\lambda, \kappa > 0$ is considered to be a parameter of the model. While the extra term allows for an increase in population, its quadratic part $-n^2$ should be understood as a damping term. Accordingly, Winkler \cite{Winkler2010, Winkler2018} shows that this prevents blowup of solutions, provided $\lambda$ is sufficiently large when compared to the initial datum, while it does not prevent blowup in the case of the elliptic-parabolic system.

\subsubsection{Chemotaxis-Consumption Systems}

Our system bears a major difference with the Keller-Segel equations discussed above, in that the chemical signal can be produced \textit{or} absorbed by the micro-organisms, depending on the sign in the third equation in \eqref{ieq:EC}. In this paragraph, we focus on the case where the chemical is absorbed (for instance, when $c(t,x)$ represents a nutrient density to which the organisms are attracted), so that the Keller-Segel equations should be replaced by the Chemotaxis-Consumption equations,
\begin{equation}\label{ieq:ChemCon}
    \begin{cases}
        \partial_t n - \Delta n + \D(n \nabla c) = 0 \\
        \partial_t c - \Delta c + nc = 0.
    \end{cases}
\end{equation}
Since the chemical mass diminishes and diffuses, we should expect it to be harder for blowup to occur (also observe that the $c$ equation has a maximum principle). This is indeed the case. We refer to the helpful survey \cite{LW2023} for clear explanations an a complete bibliography, while we shall remain content with a few selected works. 

\medskip

In almost all recent papers (with the exception of the very recent article \cite{LZ2026}, see below), proving well-posedness for the chemotaxis-consumption system used a free-energy decay law, which was first shown by Duan, Lorz and Markowich \cite{DLM}. The decay law reads
\begin{equation}\label{ieq:FE}
    \frac{\rm d}{\dt} \int \left( n \log n + \frac{|\nabla c|^2}{2c} \right) \dx + \int \left( \frac{|\nabla n|^2}{n} + c \big|\nabla^2 \log c \big| - \frac{n |\nabla c|^2}{c} \right) \dx = 0,
\end{equation}
where, although not obvious at first glance, the decay rate is in fact non-negative. Using this structure, the authors of \cite{DLM} were able to show small data well-posedness in 3D, while the 2D case was later achieved by Y. Tao, Winkler \cite{TW2012}, Winkler \cite{Winkler2012} and Jiang, Wu and Zheng \cite{JWZ2015}. In all dimensions, well posedness in high regularity Besov spaces is given by Zhang in \cite{Zhang2013}. For the rest of this paragraph, we restrict our attention to the 2D case.

\medskip

Using estimate \eqref{ieq:FE} is possible under the condition that the free-energy is initially finite. In particular, the initial data $n_0$ and $c_0$ have to be almost everywhere positive, and must satisfy $n_0 \in L\log L$ and $\nabla \sqrt{c_0} \in L^2$, where $L \log L$ is the Orlicz space defined as being the linear hull of the set of functions $f$ with $|f| \log^+ |f| \in L^1$. Accordingly, Zhang and Zheng \cite{ZZ2014} prove existence and uniqueness of solutions with data satisfying very close assumptions: $n_0, c_0 > 0$ a.e. and
\begin{equation*}
    n_0 \in L^1 \cap L^2 \qquad \nabla \sqrt{c_0} \in L^2 \qquad \text{and} \qquad c_0 \in L^1 \cap L^\infty.
\end{equation*}
These spaces are consistent with the scaling of the equations (the assumption $n_0 \in L^2$ seems to come from the fact that \cite{ZZ2014} examines a fluid-chemotaxis coupling, see next paragraph). The Chemotaxis-Consumption system is indeed also invariant under transformations \eqref{ieq:KellerSegelScaling}, so that the spaces $n_0 \in L^1$ and $c_0 \in \dot H^1 \cap L^\infty$ are scale critical in 2D. 

\medskip

At this point, it should be noted that the positivity assumptions $n_0, c_0 > 0$, while natural with respect to free-energy, may seem a bit artificial: the PDE should be solvable (perhaps globally) even if there are zones without micro-organisms or chemicals. Secondly, we remark that, while the space $c_0 \in \dot H^1$ is critical, it is not endpoint. This leaves open the question (which we will examine here) whether $\dot H^1$ can be replaced by an endpoint space, such as $\dot B^0_{\infty, 2}$.

This remark about endpoint spaces is reinforced by a very interesting paper of Zhang \cite{Zhang2016A}, a continuation criterion expressed in terms of $L^\infty$ norms only: solutions may be continued beyond a time $T < + \infty$ if and only if
\begin{equation*}
    \int_0^T \Big( \| \nabla n \|_{L^\infty} + \| \nabla c \|_{L^\infty}^2 \Big) \dt < + \infty \,.
\end{equation*}
It should be noted that the spaces implied $(n, c) \in L^1 (\dot W^{1, \infty}) \times L^2_T(\dot W^{1, \infty})$ are have consistent scaling with the parabolic nature of \eqref{ieq:ChemCon} and the assumptions $(n_0, c_0) \in L^2 \times L^\infty$ of \cite{ZZ2014}.

\medskip

We mention the very recent prerpint \cite{Hou2026}, which proves global well-posedness under a smallness condition on $\| c_0 \|_{L^\infty}$, for initial data in \eqref{ieq:Hou} with $m = 2$. Finally, we should point out that the 3D problem has received a lot of attention. In an earlier work, Y. Tao \cite{JTao2011} found small data global well-posedness, later followed by Baghaei and Khelghati \cite{BK2017}. Very rough (\textit{i.e.} renormalized) weak solutions were constructed by Li and Wang in \cite{LW2019}, while a putative blowup is examined by Jiang, Wu and Zheng in \cite{JWZ2018}.

\subsubsection{Fluid Couplings}

In our last bibliographical paragraph, we discuss chemotaxis-fluid couplings, as our main system \eqref{ieq:EC}.\footnote{In this case, the micro-organisms should be seen as micro-swimmers, such as bacteria with flagella.} The issue at hand is still whether chemotactic collapse may happen or not. Obviously, the presence of an additional equation would not prevent singularity formation, unless special regimes are studied (taking $u = 0$ in the system reduces it to previous problems), however it is interesting to know if there are some regimes for which chemotactic collapse will \textit{not} occur, or in other words whether the solutions are global under some condition. Although we are primarily focused on the fluid-chemotaxis-consumption model, we should mention the paper of Cao and Kang \cite{CK2021}, find global solutions of a Keller-Segel system coupled with the Euler equations, under a smallness assumption on the concentration.

\medskip

Many of the works we cited in the section above examine in fact chemotaxis-consumption-fluid models, namely \cite{DLM, Winkler2012, JWZ2015, Zhang2013, ZZ2014, LZ2026}. There, the system was paired with the Navier-Stokes equations by using a gravity term identical to the one in \eqref{ieq:EC},
\begin{equation}\label{ieq:FluidCoupl}
    \partial_t u - \nu \Delta u + (u \cdot \nabla) u + \nabla \pi = n g.
\end{equation}
Here, $\nu > 0$ is the fluid viscosity. Therefore, analogous results hold: global 2D solutions were found in \cite{Winkler2012, JWZ2015} under some appropriate assumption on the velocity field, often related the domain of the Stokes operator. In terms of regularity, the best result so far still is in the paper of \cite{ZZ2014}, where well-posedness is found for initial data satisfying $n_0, c_0 > 0$ almost everywhere and
\begin{equation}\label{ieq:ZData}
    n_0 \in L^1 \cap L^2, \qquad \nabla \sqrt{c_0} \in L^2, \qquad c_0 \in L^1 \cap L^\infty, \qquad \text{and} \qquad u_0 \in L^2.
\end{equation}
The assumption $u_0 \in L^2$ on the initial velocity is consistent with the scaling transformations \eqref{ieq:KellerSegelScaling} and the coupled Navier-Stokes equation \eqref{ieq:FluidCoupl}. We also point out that the existence of global 3D weak solutions has been proved by Lie and Lorz \cite{LL2011} and Winkler \cite{Winkler}, this last result holding under scale critical assumptions.

\medskip

Closer to our interests lies the coupling with Euler equations. In that case, the equations become exactly \eqref{ieq:EC}, with the exception of the reproduction term. We mention three important results related to well-posedness. The first one 3D local well-posedness in high order Sobolev spaces, proved by Zhang in \cite{Zhang2016B}. The second is from Na \cite{Na2024}, and is 2D local well-posedness for Sobolev data
\begin{equation}\label{ieq:Hou}
    n_0 \in H^m, \qquad c_0 \in H^{m+1}, \qquad \text{and} \qquad u_0 \in H^{m+1},
\end{equation}
where $m > 3$. Lastly and thirdly, the article \cite{LZ2026} of Liu and Zhang is important in two respects, as it proves global well-posedness for a minor variant of our system \eqref{ieq:EC}, \textit{without requiring the positivity assumption} $n_0, c_0 > 0$, albeit with very regular data
\begin{equation}\label{ieq:LZdata}
    n_0 \in L^1 \cap H^1, \qquad c_0 \in L^1 \cap H^2, \qquad \text{and} \qquad u_0 \in H^1 \cap B^1_{\infty, 1}.
\end{equation}
All these results involve very high regularity assumptions, which are largely above the scaling of the problem.

\subsection{Our Problem: The Euler-Chemotaxis-Consumption System}

The system we study in this article is a coupling between the Chemotaxis Consumption equations \eqref{ieq:ChemCon} and the Euler equations, with the addition of a reproduction term $n(1-n^2)$ in the micro-swimmer equation. We start by explaining some basic differences with stand-alone Chemotaxis-Consumption systems, and then see how this affects scaling considerations.

\subsubsection{Conservation Laws and the Reproduction Term}

Most global well-posedness results rely on global estimates, which sometimes can be derived from physical principles, such as conservation laws. In the case of the Euler or the Navier-Stokes equations, vorticity or energy conservation laws give 2D global well-posedness.

\medskip

Writing an energy balance in \eqref{ieq:FluidCoupl} requires the gravity term $ng$ to lie in the space $L^1_{\rm loc}(L^2)$. However, mass conservation \eqref{ieq:MassCons} is insufficient by itself to give such a bound, even when using maximal regularity in the micro-swimmer equation, and relying on other estimates for $c$.

To circumvent this problem, we introduce a reproduction term $R_3 := n(1 - n^2)$ in the first equation in \eqref{ieq:EC}, which induces some damping in the equations. In biological terms, this means that the micro-swimmers reproduce under two constraints: at small populations, the growth speed at some point $x$ and time $t$ is proportional to the density $n(t,x)$, and there is some kind of saturation phenomenon, which causes the organisms to die when the density is above the threshold $n=1$. There are many choices for such a law, the most common being the Verhulst law $R_2 := n(1-n)$, which corresponds to the so-called ``logistic'' evolution. 

The choice of a cubic law $R_3 = n(1-n^2)$ way seem surprising, as the cube $n^3$ is not consistent with the scaling of the Chemotaxis-Consumption equations \eqref{ieq:KellerSegelScaling}, which we have used so far (contrary to the Verhulst law $R_2$). However, it is motivated by important considerations related to the Euler equations.

\medskip

Firstly, the Euler equations (\eqref{ieq:FluidCoupl} with $\nu = 0$) are known to be well-posed only under the condition that the velocity field satisfies an Osgood condition (\textit{e.g.} if $u$ is $\log$-Lipschitz). This regularity follows readily from $L^\infty$ bounds on the vorticity $\omega = \partial_1 u_2 - \partial_2 u_1$, which solves the equation
\begin{equation}\label{ieq:VortEQ}
    \partial_t \omega + u \cdot \nabla \omega = g \cdot \nabla^\perp n\, .
\end{equation}
To achieve this, we will therefore need a bound of the type $\nabla n \in L^1(L^\infty)$. Now, the norm $\| \nabla n \|_{L^1(L^\infty)}$ is not critical for the scale transformation \eqref{ieq:KellerSegelScaling}. This indicates that a different set of scaling transformations might be more relevant to study \eqref{ieq:EC}, as we will see in the next paragraph. 

The space $L^1(\dot W^{1, \infty})$ shares the scaling of the ``energy'' space $L^\infty(L^2) \cap L^2(\dot H^1)$, which leads us to expect this space to play an important role in our analysis. This space is critical with respect to the 2D cubic non-linear heat equation
\begin{equation*}
    \partial_t n - \Delta n = \pm n^3,
\end{equation*}
which makes the cubic reproduction term $R_3$ especially relevant.

\subsection{Scaling Considerations and Regularity Upgrading}\label{ss:Scaling}

This paragraph is aimed at better understanding system \eqref{ieq:EC}, the relative importance of the different terms, and the type of result we can expect, and the strategy of the proof. By using the vorticity equation \eqref{ieq:VortEQ}, equations \eqref{ieq:EC} can be recast as
\begin{equation}\label{ieq:ECV}
    \begin{cases}
        \partial_t \omega + u \cdot \nabla \omega = g \cdot \nabla^\perp n \\
        \partial_t n - \Delta n + u \cdot \nabla n + n^3 - n + \D(n \nabla c) = 0\\
        \partial_t c - \Delta c + u \cdot \nabla c + nc = 0\\
        u = \nabla^\perp \Delta^{-1} \omega.
    \end{cases}
\end{equation}
This form, in particular the vorticity equation, will be of constant use in the article.

\medskip

As we have said in the paragraph above, the scaling transformations \eqref{ieq:KellerSegelScaling} are not appropriate for solving the Euler coupling. Instead, consider the following one parameter family of scaling transformations: for any $\lambda > 0$, set
\begin{equation}\label{ieq:Scaling}
    n_\lambda (t, x) = \lambda n(\lambda^2 t, \lambda x), \qquad c_\lambda (t, x) = c(\lambda^2 t, \lambda x) \qquad \text{and} \qquad u_\lambda(t, x) = \lambda u(\lambda^2 t, \lambda x).
\end{equation}
The scaling of the parameters $t$ and $x$ is related to the parabolic nature of the last two equations of \eqref{ieq:EC}. Now, while it is \textit{not} true that equations \eqref{ieq:EC} are invariant under the scaling transformation (\textit{i.e.} that $(n_\lambda, c_\lambda, u_\lambda)$ is a solution as soon as $(n, c, u)$ is), it is mostly true. More precisely, remark that the system would be scaling invariant if we removed three terms, which are typed in bold font in the equations immediately below:
\begin{equation*}
    \begin{cases}
        \partial_t u + (u \cdot \nabla) u + \nabla \pi = \textbf{ng}\\
        \partial_t n - \Delta n + u \cdot \nabla n + n^3 - \textbf{n} + \D(n \nabla c) = 0\\
        \partial_t c - \Delta c + u \cdot \nabla c \pm \textbf{nc} = 0\\
        \D(u) = 0.
    \end{cases}
\end{equation*}
They are \textit{sub-critical} with respect to \eqref{ieq:Scaling}, which means that they will involve bounds of lower order, or can be treated as perturbation terms of less importance. This fundamental remark will be at the very core of our proof strategy.

\medskip

Again, as explained in the previous paragraph, the scaling transformation is closely related to the type of well-posedness result that might be expected for \eqref{ieq:EC}, as, optimal well-posedness spaces generally have norms which are invariant under the scaling transformations. Examples of such scale-invariant spaces are
\begin{equation}\label{ieq:CriticalSpaces}
    n \in L^\infty(L^2) \qquad c \in L^\infty(L^\infty) \qquad \text{and} \qquad u \in L^\infty(L^2),
\end{equation}
or
\begin{equation}\label{ieq:TargetScale}
    n \in L^\infty(L^2) \cap L^2(\dot H^1), \qquad c \in L^\infty(L^\infty) \cap L^2(\dot W^{1, \infty}) \qquad \text{and} \qquad  \omega \in L^1(L^\infty),
\end{equation}
Before moving onward, we should make a comment on the hyperbolic nature of the vorticity equation, the first one in \eqref{ieq:ECV}. For this type of equation, uniqueness estimates require velocity fields which are (log)-Lipschitz. This is related to the fact that the Euler equations, unlike most semi-linear equations, have a two-parameter family of scaling invariance $u(t,x) \longmapsto \lambda^{-1} \mu\,  u(\mu t, \lambda x)$. This further constrains the optimal spaces, which must share the scaling of $L^1(\dot W^{1, \infty})$ in both the time and space parameters. Therefore, we should keep in mind that the optimal space for well-posedness will involve Lipschitz (or possibly log-Lipschitz) regularity of the velocity field, as is the case in \eqref{ieq:TargetScale}.

\subsection{Main Result}

In this Section, we state our main theorem, as well as give some minor comments. The remarks of the previous paragraphs already show that there is probably is some kind of optimality in terms of scaling. We will expand on this in the remainder of the introduction.

\begin{thm}\label{it:main}
        Consider a set of initial data $(n_0, c_0, u_0)$ which satisfy $n_0, c_0 \geq 0$ and
    \begin{equation}\label{ieq:OurData}
        n_0 \in L^1 \cap L^2 \cap \dot B^{-1}_{\infty, 1}, \qquad c_0 \in L^2 \cap \dot B^0_{\infty, 2}, \qquad \text{and} \qquad u_0 \in L^2 \text{ with } \omega_0 \in L^{6/5} \cap L^\infty.
    \end{equation}
    Here $\omega_0 = \curl(u_0)$. We refer to the appendix for a definition of the spaces $\dot B^{-1}_{\infty, 1}$ and $\dot B^0_{\infty, 2}$. Under the assumptions above, the following points hold:
    \begin{enumerate}[(i)]
        \item The system \eqref{ieq:EC} has a unique global solution $(n, c, u)$ that satisfies, for any time $T > 0$,
        \begin{equation}\label{ieq:Bounds1}
            \begin{split}
                &n \in L^\infty([0, T[ ; L^2) \cap L^2([0, T[ ; \dot H^1) \cap L^1([0, T[ ; \dot W^{1, \infty}) \\ 
                &c \in L^\infty([0, T[ ; L^2 \cap \dot B^0_{\infty, 2}) \cap L^2([0, T[ ; \dot H^1 \cap \dot W^{1, \infty}) \\
                &u \in L^\infty([0, T[ ; L^2) \text{ with } \curl (u) \in L^\infty([0, T[ ; L^{6/5} \cap L^\infty).
            \end{split}
        \end{equation}

        \item Moreover, the solution enjoys the following regularity properties: for any $T > 0$,
        \begin{equation}\label{ieq:Bounds2}
            \begin{split}
                & n \in L^\infty([0, T[; L^1) \cap L^3([0, T[  ; L^3) \cap L^1 \big([0, T[; \dot W^{1, \frac{6}{5}} \cap \dot B^1_{\infty, 1} \big) \\
                & c \in L^2 ([0, T[ ; \dot B^1_{\infty, 1}) \cap \widetilde{L^\infty_T}(\dot B^0_{\infty, 2}).
            \end{split}
        \end{equation}
        
        \item Consider $R > 0$, and $B_X(0, R)$ a ball of radius $R$ in the initial data space
        \begin{equation*}
            (n_0, c_0, u_0) \in X := L^1 \cap L^2 \cap \dot B^{-1}_{\infty, 1} \times L^2 \cap \dot B^0_{\infty, 2} \times L^2 \cap \{ \omega_0 \in L^{6/5} \cap L^\infty \}.
        \end{equation*}
        The the initial datum to solution map $(n_0, c_0, u_0) \longmapsto (n, c, u)$ restricted to the ball $B_X(0, R)$ is continuous in the following topology: for any $0 < s < 1$,
        \begin{equation}\label{ieq:StabSubCr}
            \begin{split}
                &\dot H^{-s} \times L^2 \cap \dot B^0_{\frac{2}{s}, 2} \times L^2  \longrightarrow \\
                & \qquad  L^\infty_{\rm loc}(\dot H^{-s}) \cap L^2_{\rm loc}(\dot H^{1-s}) \times L^\infty_{\rm loc}(L^2) \cap L^2_{\rm loc}(\dot H^1 \cap \dot B^1_{\frac{2}{s}, 2}) \times L^\infty_{\rm loc}(L^2),
            \end{split}
        \end{equation}
        with explicit continuity bounds (see Theorem \ref{t:Stability}). Remark that $\dot B^0_{\frac{s}{2}, 2} \subset L^2 \cap \dot B^0_{\infty, 2}$ (Proposition \ref{p:ComplexBesovInterpolate}) and that $\dot H^{-s} \subset L^1 \cap L^2$ (Proposition \ref{p:BesovInterpolate}).

        \item Moreover, with the previous assumptions, we also have continuity on $B_X(0, R)$ of the initial data to solution map in the critical topology: for any $0 < s < 1$,
        \begin{equation}\label{ieq:StabCrit}
            \begin{split}
                &\dot H^{-s} \cap L^2 \times L^2 \cap \dot B^0_{\infty, 2} \times L^2  \longrightarrow \\
                & \qquad  L^\infty_{\rm loc}(\dot H^{-s} \cap L^2) \cap L^2_{\rm loc}(\dot H^{1-s} \cap \dot H^1) \times L^\infty_{\rm loc}(L^2) \cap L^2_{\rm loc}(\dot H^1 \cap \dot B^1_{\infty, 1}) \times L^\infty_{\rm loc}(L^2),
            \end{split}
        \end{equation}
        again with explicit bounds (see Theorem \ref{t:t:StabilityCritical}).
    \end{enumerate}
\end{thm}

\begin{rmk}
    On the other hand, the Besov assumptions $n_0 \in \dot B^{-1}_{\infty, 1}$ and $c_0 \in \dot B^0_{\infty, 2}$ may seem surprising or artificial. In the rest of this introduction, we will explain why these are in fact natural, and probably optimal.
\end{rmk}

\begin{rmk}
    The initial data regularity \eqref{ieq:OurData} should be compared with previous results in fluid-Chemotaxis-Consumption couplings, in particular \eqref{ieq:ZData} from \cite{ZZ2014}. On the one hand, our assumption on the initial chemo-attractant concentration $c_0 \in \dot B^0_{\infty, 2}$ is more general than $c_0 \in \dot H^1$, while sharing the same scaling. On the other hand, we have an extra assumption $n_0 \in \dot B^{-1}_{\infty, 1}$ when compared to \eqref{ieq:ZData}. As we will see a bit below, this is related to the fact that we work with the Euler equations, and almost certainly cannot be dispensed with.
\end{rmk}

\begin{rmk}
    Finally, we relax the positivity assumption $n_0, c_0 > 0$ present in \cite{ZZ2014}, which enables us to consider data which cancel on open areas, and hence have infinite initial free energy. While this is the case for the initial data in \cite{LZ2026}, our method allows for much rougher initial data than \eqref{ieq:LZdata}.
\end{rmk}

\subsection{Global Estimates other than the Free Energy}\label{ss:GlobalQuestion}

In order to achieve global well-posedness, we will be looking for global estimates on the solution which are at or above critical regularity (\textit{e.g.} \eqref{ieq:CriticalSpaces} or \eqref{ieq:TargetScale}). However, as we want to allow for infinite free-energy data, estimate \eqref{ieq:FE} is no longer available. We therefore look for other types of bounds. We start by listing a few basic conservations laws, and comment on their usefulness for well-posedness.

\medskip

\textbf{Mass conservation.} Just as in \eqref{ieq:MassCons}, we have a mass conservation law for the micro-swimmer density. By integrating the $n$ equation in \eqref{ieq:EC}, we see that (recall that $\D(u) = 0$)
\begin{equation}\label{ieq:MassCNew}
    \frac{\rm d}{\dt} \int n \dx + \int n^3 \dx = \int n,
\end{equation}
which shows that $n\in L^\infty_{\rm loc}(L^1) \cap L^3_{\rm loc}(L^3)$. While this will be useful, these spaces are below critical regularity. We will need higher regularity bounds.

\medskip

\textbf{Chemo-attractant $L^2$ bounds.} By multiplying the third equation of \eqref{ieq:EC} by $c$ and integrating, we find, after some integration by parts,
\begin{equation}\label{ieq:ChemCNew}
    \frac{1}{2} \frac{\rm d}{\dt} \int c^2 \dx + \int |\nabla c|^2 \dx \pm \int nc^2 \dx = 0.
\end{equation}
In the case where the sign $\pm$ is positive, the non-negativity of $nc^2$ immediately gives $c \in L^\infty(L^2) \cap L^2(\dot H^1)$. We will see this also holds in the case of a negative sign. However, once again, this space is below critical regularity \eqref{ieq:TargetScale}.

\medskip

\textbf{Chemo-attractant maximum principle.} Again, considering the case where the sign $\pm$ in \eqref{ieq:EC} is positive, the $c$ equation satisfies a maximum principle, which gives the uniform bound $\| c \|_{L^\infty} \leq \| c_0 \|_{L^\infty}$. This is an encouraging sign, as the bound $c \in L^\infty(L^\infty)$ is critical. However, as it turns out, we will never need such a bound (which behaves badly with respect to maximal regularity for the heat equation anyway).

\medskip

\textbf{Micro-Swimmer $L^2$ estimate.} From \eqref{ieq:TargetScale}, we see that the ``energy'' space $n \in L^\infty(L^2) \cap L^2(\dot H^1)$ is critical. This prompts us to multiply the second equation in \eqref{ieq:EC} by $n$ and integrate. We have
\begin{equation}\label{ieq:nEstimateL2}
    \begin{split}
        \frac{1}{2} \frac{\rm d}{\dt} \int n^2 \dx + \int |\nabla n|^2 \dx + \int n^4 \dx & = \int n \dx + \int n \nabla n \cdot \nabla c \dx \\
        & \leq \| n \|_{L^2}^2 + \| n \|_{L^2} \| \nabla n \|_{L^2} \| \nabla c \|_{L^\infty}.
    \end{split}
\end{equation}
Observe that closing the estimates would require to already possess a bound of the type $\nabla c \in L^2_{\rm loc}(L^\infty)$, which is a critical space \eqref{ieq:TargetScale}. Obtaining such an estimate will be possible, although quite technical.

\medskip

Overall, we see that none of the ``basic'' bounds described above are satisfactory.

\subsection{Choice of Functional Setting and Sketch of Proof}

A crucial fact will come to our rescue here. As we have seen in Subsection \ref{ss:Scaling}, the coupling between the fluid and the suspended micro-organisms is made through the subcritical gravity term $ng$. To illustrate how this works, we write the kinetic energy balance in the fluid:
\begin{equation*}
    \frac{\rm d}{\dt} \int \frac{1}{2} |u|^2 \dx + \int \D \Big( \frac{1}{2} |u|^2 u \dx + \pi u \Big) = \int n g \cdot u \dx.
\end{equation*}
The integral involving the divergence disappears, as it is a total derivative. On the other hand, the mass conservation equation \eqref{ieq:MassCNew} shows that $n \in L^\infty_{\rm loc}(L^1) \cap L^3_{\rm loc}(L^3)$ (this is a subcritical bound), and so
\begin{equation*}
    \frac{1}{2} \frac{\rm d}{\dt} \| u \|_{L^2} \leq \| n \|_{L^1}^{1/4} \| n \|_{L^3}^{3/4},
\end{equation*}
and Grönwall's lemma implies that $u \in L^\infty_{\rm loc}(L^2)$, which is a critical bound \eqref{ieq:CriticalSpaces}. This fact is highly significant: by using the subcriticality of the gravity term, we have obtained a critical bound for the velocity from a subcritical bound on the micro-organism density.

\medskip

This indicates a possible strategy of proof: use the subcritical gravity coupling $ng$ to upgrade the subcritical bounds \eqref{ieq:MassCNew} and \eqref{ieq:ChemCNew} to a family of higher regularity estimates, until we reach critical regularity. There will be three such cycles. In the end, we will need bounds of the type \eqref{ieq:TargetScale}. These are summarized in the diagram \eqref{ieq:Diagram} below, where each continuous arrow represents the use of the subcritical gravity term to upgrade the level of the estimates. Observe that the final bounds for $\omega$ are in fact supercritical (because of the $L^\infty_{\rm loc}$ regularity with respect to the time variable in the hyperbolic vorticity equation), whereas the bounds for $n$ and $c$ are critical. In the remainder of this paragraph, we give a more detailed discussion on how the different bounds are proved, and the analytical tools we use to prove them.
\begin{equation}\label{ieq:Diagram}
    \xymatrix{
        \substack{c \in L^\infty(L^2) \cap L^2(\dot H^1) \\ n \in L^\infty_{\rm loc}(L^1) \cap L^3_{\rm loc}(L^3)} \ar@/^/[rrrd]^{(i)} \\
        & & & \substack{u \in L^\infty_{\rm loc}(L^2) \\ \nabla n \in L^1_{\rm loc}(L^{6/5})} \ar@/_/[llld]^{(ii)}\\
        \substack{\omega \in L^\infty_{\rm loc}(L^{6/5}) \\ \nabla c \, \in L^2_{\rm loc}(L^\infty) \\ n \in L^\infty_{\rm loc}(L^2) \cap L^2_{\rm loc}(\dot H^1)} \ar@{.>}[d]^{(iii)} & & & \\
        \substack{\nabla n \in L^1_{\rm loc}(L^\infty)} \ar@/^/[rrrd]^{(iv)} & & & \\
        & & & \substack{\omega \in L^\infty_{\rm loc}(L^\infty)}
    }
\end{equation}

We have already explained, just above, how step ($i$) can be proved. Finding $u \in L^\infty_{\rm loc}(L^2)$ regularity, it is then possible to use the basic bounds \eqref{ieq:MassCNew} and \eqref{ieq:ChemCNew} to prove that the micro-organism density satisfies a (still subcritical) bound $\nabla n \in L^1_{\rm loc}(L^{6/5})$. This is done with the help of (Lebesgue) maximal regularity estimates for the heat equation, (see Proposition \ref{p:MaxRegularity} in the appendix).

\medskip

Then, using the vorticity equation (and hence the subcritical gravity term) directly provides that $\omega \in L^{\infty}_{\rm loc}(L^{6/5})$. This corresponds to arrow ($ii$). Our next goal is to close the critical bounds $n \in L^\infty_{\rm loc}(L^2) \cap L^2_{\rm loc}(\dot H^1)$ in \eqref{ieq:nEstimateL2}, which require a preliminary estimate of the form $\nabla c \in L^2_{\rm loc}(L^\infty)$. Bounds of this type are usually technical to prove. The main reason for the difficulty is that Lebesgue $L^q(L^p)$ maximal regularity estimates fail in the case of endpoint exponents $p, q \in \{ 1, \infty\}$. As an example, consider the heat semi-group applied to a bounded initial datum $c_0 \in L^\infty(\R^2)$, which reads
\begin{equation*}
    \nabla e^{t\Delta} c_0(x) = \frac{1}{\sqrt{t}} \int c_0(x - y) \frac{y}{2\sqrt{t}} \exp \left( - \frac{|y|^2}{4t} \right) \frac{\dy}{4 \pi t}.
\end{equation*}
Observe that straightforwardly bounding $c_0$ by the supremum norm provides the inequality $\| \nabla e^{t \Delta} c_0 \|_{L^\infty} \leq \frac{1}{2 \sqrt{t}} \| c_0 \|_{L^\infty}$. This shows that the solution has Lorentz regularity with respect to the time variable $\nabla e^{t \Delta} c_0 \in L^{2, \infty}( L^\infty)$, but no better, and these estimates are optimal as is demonstrated by the 1D example of the function $c_0(x_1, x_2) = {\rm sign}(x_1)$, for which we have $\partial_1 e^{t \Delta} c_0 (x) = C/\sqrt{t}$. It is therefore necessary to use a more adapted functional framework. Function spaces for which the heat equation is better behaved are usually Besov spaces ot ${\rm BMO}$ based spaces \cite{OS2021, DWY, Chemin1999}. For us, the most natural choice, if we aim at obtaining the bound $\nabla e^{t \Delta} c_0 \in L^2(L^\infty)$, would be the space $X$ of initial concentrations such that
\begin{equation*}
    \| c_0 \|_{X}^2 :=  \int_0^{+ \infty} \| \nabla e^{t \Delta} c_0 \|_{L^\infty}^2 \dt < + \infty .
\end{equation*}
This space is well known to being the homogeneous Besov space $X = \dot B^0_{\infty, 2}$, and the relation above is the characterization of $\dot B^0_{\infty, 2}$ in terms of the heat semi-group (see Proposition \ref{p:BesovChar}). As we will see, the corresponding chemical concentration estimates are $c \in L^2_{\rm loc}(\dot B^1_{\infty, 1})$ (this remark comes from \cite{AH2025}).

\medskip

Even so, deriving $c \in L^2_{\rm loc}(\dot B^1_{\infty, 1})$ estimates is not completely straightforward. This is largely because the heat equation is better behaved in spaces that differ a bit from $L^2_{\rm loc}(\dot B^1_{\infty, 1})$, namely the class of Chemin-Lerner spaces $\widetilde{L^2_{\rm loc}}(\dot B^1_{\infty, 1}) \subsetneq L^2_{\rm loc}(\dot B^1_{\infty, 1})$, which were introduced in \cite{CL1995}. Unfortunately, these spaces are notoriously bad-behaved in Grönwall-type arguments. Dealing with Chemin-Lerner spaces often involve careful analysis (see \textit{e.g.} \cite{FX2015}). For example, in \cite{AG2020} Arsénio and Gallagher prove maximal regularity estimates in Lebesgue-Besov spaces such as our space of interest $L^2_{\rm loc}(\dot B^1_{\infty, 1})$. However, here we explore a different method, and resort to interpolation. Chemin-Lerner spaces can indeed be seen as real interpolation spaces\footnote{The author is grateful to Anatole Gaudin for this remark, as well as help on interpolation of time-space function spaces.}
\begin{equation}\label{ieq:CLabstract}
    \widetilde{L^q}(\dot B^{s}_{p, 1}) = \left[ \widetilde{L^q}(\dot B^{s_1}_{p, \infty}), \widetilde{L^q}(\dot B^{s_2}_{p, \infty}) \right]_{\theta, 1} \, ,
\end{equation}
where $1 \leq p, q \leq \infty$ and $s_1 < s = \theta s_1 + (1 - \theta)s_2$ for any $\theta \in ]0, 1[$. We prove the corresponding interpolation inequality in Proposition \ref{p:Interpolation} below.

\medskip

With the bound $\nabla c \in L^2_{\rm loc}(L^\infty)$ proved, we also have access to the critical energy estimate $n \in L^\infty_{\rm loc}(L^2) \cap L^2_{\rm loc}(\dot H^1)$. To conclude our estimates and prove, with step ($iv$), that $\omega \in L^\infty_{\rm loc}(L^\infty)$, it therefore remains to prove that $\nabla n \in L^1_{\rm loc}(L^\infty)$. This is step ($iii$), and it likewise deserves some explanation. In order to obtain the desired micro-organism density bound $\nabla n \in L^1_{\rm loc}(L^\infty)$, we look for initial data which lie in the space of functions such that
\begin{equation*}
    \| n_0 \|_{\dot B^{-1}_{\infty, 1}} := \int_0^{+\infty} \| \nabla e^{t \Delta} n_0 \|_{L^\infty} \dt < + \infty.
\end{equation*}
Once again, this is a Besov space $\dot B^{-1}_{\infty, 1}$, which has the same scaling as $L^2$, although $L^2 \not \subset \dot B^{-1}_{\infty, 1}$. The bound $\nabla n \in L^1_{\rm loc}(L^\infty)$ will also be proved with real interpolation techniques.

\begin{rmk}
    The discussion above shows that the Besov assumptions $c_0 \in \dot B^0_{\infty, 2}$ and $n_0 \in \dot B^{-1}_{\infty, 1}$ are not simply technical or artificial, but naturally emerge from the need to close basic estimates, and that they are most likely optimal.
\end{rmk}

\begin{rmk}
    Concerning the optimality of $n_0 \in \dot B^{-1}_{\infty, 1}$, it may be noted that we do not actually need $\omega \in L^1(L^\infty)$ for uniqueness in the Euler equations, but that it is enough for $\omega$ to be log-Lipschitz (by the transport estimates in \cite{Danchin2005}). For instance, a bound of the type $\omega \in L^1( \dot B^0_{\infty, \infty})$ would be enough, raising the question of knowing whether we could be content with proving that $\nabla n \in L^1(\dot B^0_{\infty, \infty})$.
    
    However, obtaining a $\dot B^0_{\infty, \infty}$ bound for $\omega$ requires knowing that the velocity is Lipschitz $\nabla u \in L^1(L^\infty)$ \cite{Vishik1998, HK2008}, which will prevent us from closing the estimates. This issue is not unlike the notoriously hard global well-posedness question for the Riesz transform problem
    \begin{equation*}
        \partial_t \omega + u \cdot \nabla \omega = \partial_1 (- \Delta)^{-1/2} \omega, \qquad \text{with } u = \nabla^\perp \Delta^{-1} \omega,
    \end{equation*}
    which is poorly behaved in $L^\infty$ \cite{EK2025}.
\end{rmk}

\textbf{Concerning Uniqueness.} To conclude the introduction, we point out that uniqueness of solutions also is delicate for solutions of critical regularity, especially in our case where our regularity assumptions \eqref{ieq:TargetScale} imply that the velocity field might not be Lipschitz, but only log-Lipschitz. 

We resort to the implementation of the celebrated Yudovich method for solutions of 2D Euler with bounded vorticity \cite{Yudovich, YudovichEN}. This track of proof relies on $L^\infty_{\rm loc}(L^2)$ norms of the velocity difference function $\delta u = u_2 - u_1$, which are scaling critical. The technical difficulties involved have prompted us to first prove stability estimates for $n$ and $c$ in subcritical spaces \eqref{ieq:StabSubCr}, before upgrading to critical spaces \eqref{ieq:StabCrit} later on.

A second feature of our critical setting is that we will have to use, and prove, critical Chemin-Lerner bounds of regularity just above that of $\delta u$. To do this, we will resort to interpolation. In order to obtain the best modulus of continuity possible, we will prove a logarithmic interpolation version of \eqref{ieq:CLabstract},
\begin{equation*}
    \| u \cdot \nabla c \|_{\widetilde{L^2}(\dot B^{-1}_{\infty, 1})} \lesssim \| u \|_{L^\infty(L^2)} \| \nabla c \|_{L^2(L^\infty)} \log \left( 1 + \frac{ \| \omega \|_{L^\infty(L^p)} \| \nabla c \|_{L^2(L^2)} }{ \| u \|_{L^\infty(L^2)} \| \nabla c \|_{L^2(L^\infty)}} \right),
\end{equation*}
where $p > 1$. This approach is not unlike that of \cite{Danchin2005B}, although our target space is Chemin-Lerner.

\subsubsection*{Summary of the Article}

Our article is divided into five main sections and an appendix. Section \ref{s:Global} contains the proof of the global \textit{a priori} estimates, which we have sketched above. Section \ref{s:Unique} focuses on proving uniqueness of solutions by means of stability estimates. Section \ref{s:Existence} the proceeds to rigorously prove the existence of solutions with the help of an iterative scheme. Section \ref{s:CritUn} proves stability estimates in critical spaces. Finally, the manuscript ends with an appendix where are gathered some harmonic analysis results from the literature.

\subsubsection*{Notation and Conventions}

We gather some more or less standard notation we will use in the sequel.

\begin{itemize}
    \item The notation $C$ will refer to a generic constant, which depends only on irrelevant parameters. Whenever we wish to give the exact dependencies of the constant, we will write $C (\, . \,)$.

    \item If $A$ and $B$ are two non-negative quantities, the proposition $A \lesssim B$ means that $A \leq CB$.

    \item unless otherwise mentioned, all functions spaces are set on $\R^2$. In order to simplify notations, we will omit the reference to $\R^2$ when writing such a space. For example, we note $L^p(\R^2) = L^p$.

    \item If $X \subset \mc S'$ is a Banach space, we denote $\| \, . \, \|_X$ the associated norm and $\langle \, . \, , \, . \, \rangle_{X' \times X}$ the associated duality bracket.

    \item For any time $T \in \R_+ \cup \{+\infty\}$ and $p \in [0, T]$, the symbol $L^p_T(X)$ will refer to the Banach space $L^p([0, T[ ; X)$, and likewise when replacing $L^p$ by $C^0$, \textit{etc}. When $T = + \infty$, we simply abridge $L^p(X)$.

\end{itemize}

\subsubsection*{Acknowledgments}

The author wishes to thank Anatole Gaudin for helpful exchanges about real interpolation and Chemin-Lerner spaces, as well as Raphaël Danchin for instructive discussions about parabolic problems.

This work was partially supported by the Agence Nationale de la Recherche (ANR) Project SUSPENSIONS, grant ANR-24-CE92-0028.

\section{Global Estimates}\label{s:Global}

In this Section, we prove global \textit{a priori} estimates: assuming the existence of smooth and sufficiently decaying solutions, we attempt to derive bounds that hold on every time interval $[0, T]$. Since different scales are involved in our problem, we will not be able to find bounds in spaces of the type $L^p(X)$. Accordingly, we will dismiss any constant which depends on $T$ only as being irrelevant. We aim at proving the following result.

\begin{thm}
    Consider a smooth and decaying solution $(n, c, u)$ of \eqref{ieq:EC},with initial data satisfying \eqref{ieq:OurData}. Then, for any $T > 0$, we have \eqref{ieq:Bounds1} and \eqref{ieq:Bounds2} with explicit bounds depending only on the initial data and on $T$.
\end{thm}

\subsection{Basic Estimates}

In this paragraph, we will have a precise discussion concerning the basic estimates mentioned in Subsection \ref{ss:GlobalQuestion}.

\begin{prop}\label{p:Basic}
    Consider a smooth and decaying solution of \eqref{ieq:EC} with initial data $n_0 \in L^1$ and $c_0 \in L^2$, which are assumed to be non-negative $n_0, c_0 \geq 0$. Then, the solutions are also non-negative $n, c \geq 0$ and, for any $T > 0$, they satisfy the bounds
    \begin{equation*}
        \begin{split}
            & \| n \|_{L^\infty_T(L^1)} + \| n \|_{L^3_T(L^3)}^3 \lesssim \| n_0 \|_{L^1}\\
            & \| c \|_{L^\infty_T(L^2)} + \| c \|_{L^2_T(\dot H^1)} \lesssim \| c_0 \|_{L^2} \exp \big( Ce^T \| n_0 \|_{l^1} \big).
        \end{split}
    \end{equation*}
    Note that this bound holds regardless of the sign $\pm$ in \eqref{ieq:EC}.
\end{prop}

\begin{proof}
    Our first remark is that the micro-organism density and the concentration are both non-negative functions at all times, provided they are initially, and that the solutions are smooth. This can be seen by expanding the second and third equation in \eqref{ieq:EC}
    \begin{equation*}
        \begin{split}
            \partial_t n + (u + \nabla c) \cdot \nabla n - \Delta n = n\Big(1 - n^2 + \Delta c \Big).
        \end{split}
    \end{equation*}
    Observe that the operator $\partial_t + (u + \nabla c) - \Delta$ satisfies a maximum principle, while the right-hand side in the above is a multiple of $n$. This forbids $n$ from becoming negative provided the initial datum satisfies $n_0(x) \geq 0$. Likewise, a minimum principle applied to the chemical concentration equation
    \begin{equation*}
        \big( \partial_t + u \cdot - \Delta \big)c = \mp nc
    \end{equation*}
    shows that $c(t, x) \geq 0$ provided the initial datum satisfies $c_0(x) \geq 0$. We will always assume that $n, c \geq 0$.

    \medskip

    \textbf{Population Estimate:} By integrating the $n$ equation in \eqref{ieq:EC}, we find a differential equation solved by the total micro-swimmer population, namely
    \begin{equation*}
        \frac{\rm d}{\dt} \int n + \int n^3 = \int n.
    \end{equation*}
    Therefore, by resorting to Grönwall's inequality, we see that the micro-swimmer population $\int n = \| n(t) \|_{L^1}$ always remains finite, provided that it is finite initially with $\| n \|_{L^1} \leq e^t \| n_0 \|_{L^1}$. In addition, integrating the differential inequality above gives:
    \begin{equation}\label{ieq:MassConserv}
        \| n \|_{L^\infty_T(L^1)} + \| n \|_{L^3_T(L^3)} \leq e^T \| n_0 \|_{L^1}.
    \end{equation}
    This is the place where our choice of cubic reproduction law $n(1 - n^2)$ plays a role in the analysis of the equations: solutions naturally have $n \in L^3_{\rm loc}(L^3)$ regularity.

    \medskip

    \textbf{Chemical ``Energy'' Bounds:} integrating the chemical concentration equation would also show that the total chemical mass is also finite at every time, provided it is at initial time $c_0 \in L^1$. However, much better bounds are available. Multiplying the third equation in \eqref{ieq:EC} by $c$ and integrating yields, after integration by parts,
    \begin{equation*}
        \frac{\rm d}{\dt} \int c^2 + \int |\nabla c|^2 = \mp \int nc^2.
    \end{equation*}
    In the case where the chemical is a nutrient that is consumed by the micro-organisms, the sign $\pm$ in \eqref{ieq:EC} is positive, and the right-hand side above becomes negative. This immediately gives a global bound
    \begin{equation*}
        c \in L^\infty(\R_+ ; L^2) \cap L^2(\R_+ ; \dot H^1).
    \end{equation*}
    In the case where the chemical is a signal produced by the micro-swimmers, it is also possible to show a comparable bound. The Gagliardo-Nirenberg inequality (Proposition \ref{p:GN} below) provides
    \begin{equation*}
        \begin{split}
            \int nc^2 \leq \| n \|_{L^3} \| c \|_{L^6}^2 & \leq C \| n \|_{L^3} \| c \|_{L^2}^{2/3} \| \nabla c \|_{L^2}^{4/3} \\
            & \leq \frac{1}{2} \| \nabla c \|_{L^2}^2 + C \| n \|_{L^3}^3 \| n \|_{L^2}^2.
        \end{split}
    \end{equation*}
    It is then possible to use Grönwall's inequality, recalling that $n \in L^3_{\rm loc}(L^3)$, to show that
    \begin{equation*}
        \sup_{[0, t]} \| c \|_{L^2}^2 + \int_0^t \| \nabla c \|_{L^2}^2 \leq \| c_0 \|_{L^2} \exp \big( C e^t \| n_0 \|_{L^1} \big).
    \end{equation*}
    In both cases, we have the following bound
    \begin{equation}\label{ieq:ConcEN}
        c \in L^\infty_{\rm loc}(\R_+ ; L^2) \cap L^2_{\rm loc}(\R_+ ; \dot H^1).
    \end{equation}
\end{proof}

\subsection{Critical Velocity Energy Estimate}

From the ``basic'' estimates of Proposition \ref{p:Basic}, we use the subcriticality of the gravity coupling term to obtain critical velocity bounds. This is step $(i)$ in our proof scheme \eqref{ieq:Diagram}.

\begin{prop}\label{p:CriticalVelEN}
    Assume that $(n, c, u)$ is a smooth solution of the Euler-chemotaxis system \eqref{ieq:EC} satisfying the assumptions of Proposition \ref{p:Basic}. In particular $n \in L^\infty_T(L^1) \cap L^3_T(L^3)$ for every $T > 0$. 
    
    If the initial velocity lies in $u_0 \in L^2$, then $u \in L^\infty_T(L^2)$ and the following bound holds:
    \begin{equation*}
        \| u \|_{L^\infty_T(L^2)} \leq \| u_0 \|_{L^2} + 2 T^{4/3} \| n \|_{L^\infty_T(L^1)}^{1/4} \| n \|_{L^3_T(L^3)}^{3/4}.
    \end{equation*}
\end{prop}

\begin{proof}
    The proof is based on an energy bound for the fluid part of the system:
    \begin{equation*}
        \partial_t u + (u \cdot \nabla) u + \nabla \pi = ng.
    \end{equation*}
    By taking the dot product with the velocity $u$ and integrating in space we get a conservation of energy relation,
    \begin{equation*}
        \frac{1}{2} \frac{\rm d}{\dt} \int |u|^2 = \int n g \cdot u \leq \| n \|_{L^2} \| u \|_{L^2}.,
    \end{equation*}
    which is equivalent to (by factoring a norm $\| u \|_{L^2}$)
    \begin{equation}\label{eq:velENEQ1}
        \frac{1}{2} \frac{\rm d}{\dt} \| u \|_{L^2} \leq \| n \|_{L^2}.
    \end{equation}
    This means that we will be able to close the estimate if we can show that the micro-swimmer density lies in $L^1_T(L^2)$. Recall that conservation of mass implies that $n \in L^\infty_T(L^1) \cap L^3(L^3)$. We interpolate between these two spaces by using Hölder's inequality, which states that
    \begin{equation*}
        \| n \|_{L^2} \leq \| n \|_{L^1}^{1/4} \| n \|_{L^3}^{3/4}.
    \end{equation*}
    This inequality implies that
    \begin{equation*}
        \| n \|_{L^1_T(L^2)} \leq \| n \|_{L^\infty_T(L^1)}^{1/4} \int_0^T \| n \|_{L^3}^{3/4} \dt \leq T^{4/3} \| n \|_{L^\infty_T(L^1)}^{1/4} \| n \|_{L^3_T(L^3)}^{3/4}.
    \end{equation*}
    Plugging this inequality in \eqref{eq:velENEQ1} gives the velocity bound we were looking for,
    \begin{equation*}
        \| u \|_{L^\infty_T(L^2)} \leq \| u_0 \|_{L^2} + 2 T^{4/3} \| n \|_{L^\infty_T(L^1)}^{1/4} \| n \|_{L^3_T(L^3)}^{3/4}.
    \end{equation*}
    
\end{proof}

\subsection{Slightly Supercritical Vorticity Estimates}

Now that we have critical bounds for the velocity from Proposition \ref{p:CriticalVelEN}, we can use them to improve the bounds on $n$. Using one more time the subcriticality of the gravity coupling term, this gives estimates for $u$ that are just above critical regularity.

\begin{prop}\label{p:VortSuper}
    Let $T > 0$. Assume that $(n, c, u)$ is a smooth solution of \eqref{ieq:EC}, associated to initial data $(n_0, c_0, u_0) \in L^1 \times L^2 \times L^2$ such that $\curl(u_0) = \omega_0 \in L^{6/5}$. In particular, the bounds in Propositions \ref{p:Basic} and \ref{p:CriticalVelEN} are assumed to hold: 
    \begin{equation*}
        n \in L^\infty_T(L^1) \cap L^3_T(L^3), \qquad c \in L^\infty(L^2) \cap L^2(\dot H^1) \qquad \text{and} \qquad u \in L^\infty_T(L^2).
    \end{equation*}
    Then $\nabla n \in L^1_T(L^{6/5})$ and $\omega \in L^\infty_T(L^{6/5})$. More precisely, the following two estimates hold:
    \begin{equation*}
        \begin{split}
            \| \nabla n \|_{L^1_T(L^{6/5})} & \lesssim \| n_0 \|_{L^1} + \| n \|_{L^3_T(L^3)} \| u \|_{L^\infty_T(L^2)} + \| n \|_{L^3_T(L^3)} \| c \|_{L^2(\dot H^1)} \\
            & \qquad \qquad \qquad \qquad  \qquad \qquad  + \| n \|_{L^\infty_T(L^1)}^{3/4} \| n \|_{L^3_T(L^3)}^{1/4} + \| n \|_{L^3_T(L^3)}^3
        \end{split}
    \end{equation*}
    and
    \begin{equation*}
        \begin{split}
            \| \omega \|_{L^\infty_T(L^{6/5})} & \lesssim \| \omega_0 \|_{L^{6/5}} + \| n_0 \|_{L^1} \\
            & \qquad + \| n \|_{L^3_T(L^3)} \| u \|_{L^\infty_T(L^2)} + \| n \|_{L^3_T(L^3)} \| c \|_{L^2(\dot H^1)} \\
            &  \qquad + \| n \|_{L^\infty_T(L^1)}^{3/4} \| n \|_{L^3_T(L^3)}^{1/4} + \| n \|_{L^3_T(L^3)}^3.
        \end{split}
    \end{equation*}
\end{prop}

\begin{proof}
    Recall that the vorticity equation reads
    \begin{equation*}
        \partial_t \omega + u \cdot \nabla \omega = \nabla^\perp n \cdot g,
    \end{equation*}
    which is a transport equation with measure preserving flow. In order to obtain $L^p$ estimates for the vorticity, we must first find a $L^1(L^p)$ estimate for the derivative $\nabla n$. The proof therefore centers on the micro-swimmer density equation, which can be recasted with Duhamel's formula as
    \begin{equation*}
        \nabla n = \nabla e^{t \Delta}n_0 + \int_0^t \nabla e^{(t - s)\Delta} \Big( n - n^3 - \D (n \nabla c) - u \cdot \nabla n \Big) \, {\rm d}s.
    \end{equation*}
    Our main tool will be the maximal regularity estimates for the heat equation enclosed in Propositions \ref{p:MaxRegularity} and Proposition \ref{p:heat}. Consequently, we bound the right-hand side terms in the equation above. We start by the convective term $u \cdot \nabla n = \D (nu)$. Proposition \ref{p:MaxRegularity} gives
    \begin{equation}\label{eq:VortEQ1}
        \left\| \int_0^t e^{(t - s)\Delta} \nabla \D (nu) \, {\rm d}s \right\|_{L^{3}_T(L^{6/5})} \lesssim \| nu \|_{L^{3}_T(L^{6/5})} \lesssim \| n \|_{L^3_T(L^3)} \| u \|_{L^\infty_T(L^2)}.
    \end{equation}
    For the chemotactic term, we have, likewise
    \begin{equation}\label{eq:VortEQ2}
        \left\| \int_0^t e^{(t - s)\Delta} \nabla \D (n \nabla c) \, {\rm d}s \right\|_{L^{6/5}_T(L^{6/5})} \lesssim \| n \nabla c \|_{L^{6/5}_T(L^{6/5})} \lesssim \| n \|_{L^3_T(L^3)} \| \nabla c \|_{L^2(L^2)}.
    \end{equation}
    Next comes the reproduction terms $n - n^3$. We would like the final bound to be in a space similar to the ones above (\textit{i.e.} some $L^q_T(L^{6/5})$) in order to simplify the vorticity estimates. First look at the linear part $n$. Note that, for all functions $f$, of time (with values in some Banach space), we have
    \begin{equation*}
        \| f \|_{L^{4/3}([0, T])} \leq |f(0)| + \| f \|_{\dot{W}^{1, \frac{4}{3}}([0, T])}.
    \end{equation*}
    This implies that we may bound a Duhamel term by it's time derivative, which is controlled by maximal regularity (Proposition \ref{p:MaxRegularity}). We deduce
    \begin{equation}\label{eq:VortEQ3}
        \begin{split}
            \left\| \int_0^t e^{(t - s)\Delta} n \, {\rm d}s \right\|_{L^{4/3}_T(L^{6/5})} & \lesssim \left\| \partial_t \int_0^t e^{(t - s)\Delta} n \, {\rm d}s \right\|_{L^{4/3}_T(L^{6/5})} \\
            & \lesssim \| n \|_{L^{4/3}_T(L^{6/5})} \\
            & \lesssim \| n \|_{L^\infty_T(L^1)}^{3/4} \| n \|_{L^3_T(L^3)}^{1/4}.
        \end{split}
    \end{equation}
    We now examine the cubic term $n^3$. Given our regularity assumptions, the best bound we can find for it is the endpoint space $n^3 \in L^1_T(L^1)$, for which Lebesgue maximal regularity fails. We must therefore resort to maximal regularity estimates in Chemin-Lerner spaces. Note that $L^1_T(L^1) \subset \widetilde{L^1_T}(\dot B^0_{1, \infty})$. Therefore, by Proposition \ref{p:heat},
    \begin{equation*}
        \left\| \int_0^t e^{(t - s)\Delta} \big( n^3 \big) \, {\rm d}s \right\|_{\widetilde{L^\infty_T}(\dot B^0_{1, \infty}) \cap \widetilde{L^2_T}(\dot B^1_{1, \infty})} \lesssim \| n \|_{L^3_T(L^3)}^3.
    \end{equation*}
    In addition, thanks to the inclusion $L^\infty([0, T]) \subset L^2([0, T])$, we have $\widetilde{L^\infty_T}(\dot B^0_{1, \infty}) \subset \widetilde{L^2_T}(\dot B^0_{1, \infty})$. This allows us to use the interpolation inequality of Proposition \ref{p:Interpolation}, which provides, for any $f$,
    \begin{equation*}
        \| f \|_{L^2_T(L^{6/5})} \lesssim \| f \|_{L^2_T(\dot B^{1/3}_{1, 1})} \lesssim \| f \|_{\widetilde{L^2_T}(\dot B^0_{1, \infty})}^{2/3} \| f \|_{\widetilde{L^2_T}(\dot B^1_{1, \infty})}^{1/3}.
    \end{equation*}
    This shows that we have in fact
    \begin{equation}\label{eq:VortEQ4}
        \left\| \int_0^t e^{(t - s)\Delta} \big( n^3 \big) \, {\rm d}s \right\|_{L^2_T(L^{6/5})} \lesssim \| n \|_{L^3_T(L^3)}^3.
    \end{equation}
    All that remains is to bound the initial datum term $e^{t \Delta}n_0$. In order to obtain a $L^{6/5}$ bound from Chemin-Lerner estimates, we invoke Proposition \ref{p:Interpolation}, which gives alongside Proposition \ref{p:heat},
    \begin{equation}\label{eq:VortEQ5}
        \begin{split}
            \| \nabla e^{t \Delta} n_0 \|_{L^1_T(L^{6/5})} & \leq \| \nabla e^{t \Delta} n_0 \|_{L^1_T(\dot B^0_{\frac{6}{5}, 1})} \\
            & \lesssim \| \nabla e^{t \Delta} n_0 \|_{\widetilde{L^1_T}(\dot B^{2/3}_{\frac{6}{5}, \infty})}^{1/3}  \| \nabla e^{t \Delta} n_0 \|_{\widetilde{L^1_T}( \dot B^{-1/3}_{\frac{6}{5}, \infty})}^{2/3} \\
            & \lesssim \| \nabla e^{t \Delta} n_0 \|_{\widetilde{L^1}(\dot B^1_{1, \infty})}^{1/3} \| \nabla e^{t \Delta} n_0 \|_{\widetilde{L^2}(\dot B^0_{1, \infty})}^{2/3} \\
            & \lesssim \| n_0 \|_{L^1}.
        \end{split}
    \end{equation}
    Putting estimates \eqref{eq:VortEQ1}, \eqref{eq:VortEQ2},  \eqref{eq:VortEQ3}, \eqref{eq:VortEQ4} and \eqref{eq:VortEQ5} together in Duhamel's formula gives the micro-swimmer density bound we were looking for, namely
    \begin{equation*}
        \begin{split}
            \| \nabla n \|_{L^1_T(L^{6/5})} & \lesssim \| n_0 \|_{L^1} + \| n \|_{L^3_T(L^3)} \| u \|_{L^\infty_T(L^2)} + \| n \|_{L^3_T(L^3)} \| c \|_{L^2(\dot H^1)} \\
            & \qquad \qquad \qquad \qquad  \qquad \qquad  + \| n \|_{L^\infty_T(L^1)}^{3/4} \| n \|_{L^3_T(L^3)}^{1/4} + \| n \|_{L^3_T(L^3)}^3.
        \end{split}
    \end{equation*}
    We may now use this bound in the vorticity equation. The bound reads
    \begin{equation*}
        \begin{split}
            \| \omega \|_{L^\infty_T(L^{6/5})} & \lesssim \| \omega_0 \|_{L^{6/5}} + \| n_0 \|_{L^1} \\
            & \qquad + \| n \|_{L^3_T(L^3)} \| u \|_{L^\infty_T(L^2)} + \| n \|_{L^3_T(L^3)} \| c \|_{L^2(\dot H^1)} \\
            &  \qquad + \| n \|_{L^\infty_T(L^1)}^{3/4} \| n \|_{L^3_T(L^3)}^{1/4} + \| n \|_{L^3_T(L^3)}^3.
        \end{split}
    \end{equation*}
    This is the slightly supercritical vorticity bound we were looking for.
\end{proof}

\subsection{Critical Concentration Estimates}

\begin{prop}\label{p:ConcLip}
    Let $T > 0$ be any time. Consider a smooth solution $(n, c, u)$ of \eqref{ieq:EC} associated to the initial data $(n_0, c_0, u_0) \in L^1 \times L^2 \times L^2$ such that $\omega_0 \in L^{6/5}$ and further assume that the initial concentration satisfies $c_0 \in \dot B^0_{\infty, 2}$. We also assume that the bounds of Propositions \ref{p:Basic}, \ref{p:CriticalVelEN} and \ref{p:VortSuper} hold, so that in particular:
    \begin{equation*}
        n \in L^3_T(L^3), \qquad c \in L^\infty(L^2) \cap L^2(\dot H^1) \qquad \text{and} \qquad \omega \in L^\infty_T(L^{6/5}).
    \end{equation*}
    Then the concentration satisfies the following bound:
    \begin{equation*}
        \begin{split}
            \| c \|_{L^\infty_T(\dot B^0_{\infty, 2})} +  \| \nabla c \|_{L^2_T(L^\infty)} & \lesssim \| c \|_{L^2_T(\dot B^1_{\infty, 1})} + \| c \|_{\widetilde{L^\infty_T}(\dot B^0_{\infty, 2})} \\
            & \lesssim \| c_0 \|_{\dot B^0_{\infty, 2}} + \| \omega \|_{L^\infty_T(L^{6/5})}^{3} \| c_0 \|_{L^2} + \| n \|_{L^3_T(L^3)}^{8/5} \| c_0 \|_{L^2}.
        \end{split}
    \end{equation*}
\end{prop}

This result is a direct consequence of the following proposition, which contains slightly more general estimates which will be useful for establishing the existence of weak solutions.
\begin{prop}
    Let $T > 0$ be a fixed time, $n \in L^3_T(L^3)$ be a non-negative function and $u \in L^\infty_T(L^2)$ be a divergence-free vector field with $\omega \in L^\infty_T(L^p)$ for some $p \in ]1; 2[$. Consider $c_0 \in L^2 \cap \dot B^0_{\infty, 2}$ a non-negative solution of the concentration equation
    \begin{equation*}
        \partial_t c - \Delta c + u \cdot \nabla c + nc = 0.
    \end{equation*}
    Then the solution satisfies the bounds
    \begin{equation*}
        \begin{split}
            \| \nabla c \|_{L^2_T(L^\infty)} & \lesssim \| c \|_{L^2_T(\dot B^1_{\infty, 1})} + \| c \|_{\widetilde{L^\infty_T}(\dot B^0_{\infty, 2})} \\
            & \lesssim \| c_0 \|_{\dot B^0_{\infty, 2}} + \| \omega \|_{L^\infty_T(L^p)}^{\frac{p'}{2}} \| c_0 \|_{L^2} + \| n \|_{L^3_T(L^3)}^{8/5} \| c_0 \|_{L^2}.
        \end{split}
    \end{equation*}

\end{prop}

\begin{proof}
    We write Duhamel's formula for the concentration equation, which we consider as a heat equation with $- u \cdot \nabla c - nc$ as a force term. This yields
    \begin{equation*}
        \begin{split}
            c & = e^{t \Delta}c_0 - \int_0^t e^{(t - s)\Delta}  u \cdot \nabla c  \, {\rm d s} - \int_0^t e^{(t - s)\Delta} nc \, {\rm d s}\\
            & = e^{t \Delta}c_0 - c_1 - c_2
        \end{split}
    \end{equation*}
    The term $c_1$ is related to convection, and the term $c_2$ is related to nutrient consumption.

    \medskip
    
    \textbf{STEP 1: Initial datum term.} On the one hand, Propositions \ref{p:AGHeat} and \ref{p:heat} shows that the initial datum term satisfies the estimate
    \begin{equation}\label{eq:ConcFinal1}
        \| e^{t \Delta}c_0 \|_{L^2(\dot B^1_{\infty, 1})} + \| e^{t \Delta}c_0 \|_{\widetilde{L^\infty}(\dot B^0_{\infty, 2})} \lesssim \| c_0 \|_{\dot B^0_{\infty, 2}}
    \end{equation}

    \medskip

    \textbf{STEP 2: Convective term.} We start with the $c_1$ term. Observe that it is scaling invariant with respect to the transformations \eqref{ieq:Scaling}, and this should be reflected in the estimate we find. On the one hand, using Minkowski's inequality (Proposition \ref{p:Min}) and the Besov heat estimates of Proposition \ref{p:heat}, we find that
    \begin{equation*}
        \begin{split}
            \left\| \int_0^t e^{(t - s)\Delta} \big( u \cdot \nabla c \big)  \, {\rm d s} \right\|_{L^2_T(\dot B^1_{\infty, 1})} & \leq \left\| \int_0^t e^{(t - s)\Delta} \big( u \cdot \nabla c \big) {\rm d s} \right\|_{\widetilde{L^2_T}(\dot B^1_{\infty, 1})} \\
            & \lesssim \| u \cdot \nabla c \|_{\widetilde{L^2_T}(\dot B^{-1}_{\infty, 1})}.
        \end{split}
    \end{equation*}
    We will estimate these last Besov norms by using the homogeneous Bony decomposition (see Subsection \ref{ass:Bony} in the appendix). Starting with the convection term, we use the divergence-free property of $u$, which yields $u \cdot \nabla c = \D (uc)$, to write
    \begin{equation}\label{eq:ConcBony1}
        u \cdot \nabla c = \D \dot{ \mc T}_u (c) + \dot{ \mc T}_{\nabla c}(u) + \D \dot{ \mc R}(u, c).
    \end{equation}
    As in the rest of this article, our main tool will be interpolation. More precisely, we use the two levels of regularity that are relevant for the concentration, the energy bound $c \in L^2_T(\dot H^1)$ and the sought bound $c \in L^2_T(\dot B^1_{\infty, 1})$. In order to guide our choice of interpolation exponents, note that the vorticity bound $\omega \in L^\infty_T(L^p)$ provides thanks the Besov embeddings of Propositions \ref{p:BesovFine} and the Calder\'on-Zygmund estimates of Proposition \ref{p:CZ},
    \begin{equation*}
        \| u \|_{L^\infty_T(\dot B^{1 - \frac{2}{p}}_{\infty, 2})} \lesssim \| u \|_{L^\infty_T(\dot B^1_{p, 2})} \lesssim \| \nabla u \|_{L^\infty_T(L^p)} \lesssim \| \omega \|_{L^\infty_T(L^p)}.
    \end{equation*}

    \medskip

    We now examine the terms in the Bony decomposition one by one, starting with the paraproduct $\D \dot {\mc T}_u(c)$, which carries the regularity of the concentration. Note that, thanks to the interpolation inequality in Proposition \ref{p:Interpolation} with $\theta = 2 \left( 1 - \frac{1}{p} \right) \in ]0, 1[$, we have a Chemin-Lerner bound on the concentration:
    \begin{equation*}
        \| c \|_{\widetilde{L^2_T}(\dot B^{- 1 + \frac{2}{p}}_{\infty, 1})} \lesssim \| c \|_{L^2_T(\dot B^0_{\infty, \infty})}^{2 \left( 1 - \frac{1}{p} \right)} \| c \|_{L^2_T(\dot B^1_{\infty, \infty})}^{\frac{2}{p} - 1}.
    \end{equation*}
    The embedding $\dot H^1 = \dot B^1_{2, 2} \subset \dot B^0_{\infty, \infty}$ then gives
    \begin{equation}\label{eq:ConcInterpolate}
        \| c \|_{\widetilde{L^2_T}(\dot B^{- 1 + \frac{2}{p}}_{\infty, 1})} \lesssim \| c \|_{L^2_T(\dot H^1)}^{2 \left( 1 - \frac{1}{p} \right)} \| c \|_{L^2_T(\dot B^1_{\infty, \infty})}^{\frac{2}{p} - 1}.
    \end{equation}
    Now, the paraproduct estimates of Proposition \ref{p:Bony} shows that the operator $\dot{\mc T}$ is bounded in the $\dot B^{1 - \frac{2}{p}}_{\infty, \infty} \times \dot B^{- 1 + \frac{2}{p}}_{\infty, 1} \longrightarrow \dot B^{0}_{\infty, 1}$ topology (recall that since $p < 2$, we have $1 - \frac{2}{p} < 0$). However, we cannot directly use this result, as we seek Chemin-Lerner estimates, which Proposition \ref{p:Bony} does not immediately provide. The appropriate computation is nevertheless straightforward and follows the same line: by using the definition of the paraproduct operator \eqref{aeq:Paraproduct},
    \begin{equation*}
        \begin{split}
            \| \dot {\mc T}_u (c) \|_{\widetilde{L^2_T}(\dot B^0_{\infty, 1})} & = \sum_{j \in \Z} \Big\| \dot \Delta_j \sum_{m \in \Z} S_{m-1} u \, \dot \Delta_m c \Big\|_{L^2_T(L^\infty)} \\
            & \lesssim \sum_{j \in \Z} \sum_{|j - m| \leq 4} \| S_{m-1} u \|_{L^\infty_T(L^\infty)} \| \dot \Delta_m c \|_{L^2_T(L^\infty)}.
        \end{split}
    \end{equation*}
    Above, we have used the fact that the functions $S_{m-1} u \dot \Delta_m c$ are spectrally supported in an annulus of size $\sim 2^m$. Now, since $1 - \frac{2}{p} < 0$, the following inequality holds:
    \begin{equation*}
        \| S_{m-1} u \|_{L^\infty} \leq \sum_{n \leq m-2} \| \dot \Delta_n u \|_{L^\infty_T(L^\infty)} \leq \sum_{n \leq m-2} 2^{-n \left( 1 - \frac{2}{p} \right)} \| u \|_{\dot B^{-1 + \frac{2}{p}}_{\infty, \infty}} \lesssim 2^{m \left( -1 + \frac{2}{p} \right)} \| u \|_{\dot B^{-1 + \frac{2}{p}}_{\infty, \infty}}.
    \end{equation*}
    Plugging this in the paraproduct bound gives:
    \begin{equation*}
        \| \dot {\mc T}_u (c) \|_{\widetilde{L^2_T}(\dot B^0_{\infty, 1})} \lesssim \sum_{j \in \Z} \sum_{|j - m| \leq 4} 2^{-m \left( -1 + \frac{2}{p} \right)} \| \omega \|_{L^\infty_T(L^p)} 2^{- m \left( -1 + \frac{2}{p} \right)} c_m(t) \| c \|_{\widetilde{L^2_T}(\dot B^{- 1 + \frac{2}{p}}_{\infty, 1})},
    \end{equation*}
    where the $c_m(t)$ are non-negative numbers such that $\sum_{m \in \Z} c_m(t) = 1$ for all $t \geq 0$. With \eqref{eq:ConcInterpolate}, this finally yields a bound for the first paraproduct:
    \begin{equation}\label{eq:ConCIN1}
        \begin{split}
            \| \dot {\mc T}_u (c) \|_{\widetilde{L^2_T}(\dot B^0_{\infty, 1})} & \lesssim \| \omega \|_{L^\infty_T(L^p)} \| c \|_{L^2_T(\dot B^0_{\infty, \infty})}^{2 \left( 1 - \frac{1}{p} \right)} \| c \|_{L^2_T(\dot B^1_{\infty, \infty})}^{\frac{2}{p} - 1} \\
            & \lesssim \| \omega \|_{L^\infty_T(L^p)} \| c \|_{L^2_T(\dot H^1)}^{2 \left( 1 - \frac{1}{p} \right)} \| c \|_{L^2_T(\dot B^1_{\infty, \infty})}^{-1 + \frac{2}{p}}
        \end{split}
    \end{equation}
    We will try to find similar inequalities to bound the second paraproduct and the remainder in \eqref{eq:ConcBony1}.
    
    \medskip

    Now, focus on the second paraproduct $\dot{\mc T}_{\nabla c}(u)$. There is a slight difficulty when compared to the previous term: we have shown above that the paraproduct operator $\dot{\mc T}_u$ sends Chemin-Lerner spaces to appropriate Chemin-Lerner spaces $\widetilde{L^2_T}(\dot B^{- 1 + \frac{2}{p}}_{\infty, 1}) \longrightarrow \widetilde{L^2_T}(\dot B^0_{\infty, 1})$, and we would expect an analogous result to hold for the operator $\dot{\mc T}_{\nabla c}$. Unfortunately, the best Besov velocity bound we have at our disposal $u \in L^\infty_T(\dot B^1_{p, 2})$ (recall that $L^p \subset \dot B^0_{p, 2}$ by Proposition \ref{p:BesovFine}) does not yield an estimate in a comparable Chemin-Lerner space, as Minkowski's inequality gives the reverse inclusion $\widetilde{L^\infty_T}(\dot B^1_{p, 2}) \subset L^\infty_T(\dot B^1_{p, 2})$. In order to still obtain a Chemin-Lerner estimate, we will use our interpolation result, Proposition \ref{p:Interpolation}, directly on the paraproduct $\dot{\mc T}_{\nabla c}(u)$. Proposition \ref{p:Bony} directly gives, on the one hand,
    \begin{equation*}
        \| \dot{\mc T}_{\nabla c}(u) \|_{L^2_T (\dot B^{1 - \frac{2}{p}}_{\infty, \infty})} \lesssim \| \nabla c \|_{L^2_T(L^\infty)} \| u \|_{L^\infty_T(\dot B^{1 - \frac{2}{p}}_{\infty, \infty} )} \lesssim \| c \|_{L^2_T(\dot B^1_{\infty, 1})} \| \omega \|_{L^\infty_T(L^p)},
    \end{equation*}
    and on the other hand,
    \begin{equation*}
        \| \dot{\mc T}_{\nabla c}(u) \|_{L^2_T (\dot B^{-2/p}_{\infty, \infty})} \lesssim \| \nabla c \|_{L^2_T(\dot B^{-1}_{\infty, \infty})} \| u \|_{L^\infty_T(\dot B^{1 - \frac{2}{p}}_{\infty, \infty})} \lesssim \| c \|_{L^2_T(\dot H^1)} \| \omega \|_{L^\infty_T(L^p)}.
    \end{equation*}
    Therefore, applying Proposition \ref{p:Interpolation} with regularity exponents $-\frac{2}{p} < -1 < 1 - \frac{2}{p}$ and interpolation parameter $\theta = 2 \left( 1 - \frac{1}{p} \right)$ gives the desired bound:
    \begin{equation}\label{eq:ConCIN2}
        \| \dot{\mc T}_{\nabla c}(u) \|_{\widetilde{L^2_T}(\dot B^{-1}_{\infty, 1})} \lesssim \| \omega \|_{L^\infty_T(L^p)} \| c \|_{L^2_T(\dot H^1)}^{2 \left( 1 - \frac{1}{p} \right)} \| c \|_{L^2_T(\dot B^1_{\infty, 1})}^{-1 + \frac{2}{p}}.
    \end{equation}
    Note that this is the same upper bound as \eqref{eq:ConCIN1}.

    \medskip

    Finally, we deal with the remainder term. Just as for the previous paraproduct, we will use interpolation directly on $\dot{\mc R}(u, c)$ to get a Chemin-Lerner estimate. To use Proposition \ref{p:Bony}, we need to make sure that the sum of regularities of $u$ and $c$ is positive. We have, firstly, thanks to the Besov embeddings \eqref{p:BesovEmbed} in Appendix \ref{A}, the embedding $\dot B^1_{p, 2} \subset \dot B^{2 \left( 1 - \frac{1}{p} \right)}_{2, 2}$, so
    \begin{equation*}
        \begin{split}
            \| \dot{\mc R}(u, c) \|_{L^2_T(\dot B^{1 - \frac{2}{p}}_{\infty, 1})} & \lesssim \| \dot{\mc R}(u, c) \|_{L^2_T(\dot B^{3 - \frac{2}{p}}_{1, 1})} \\ 
            & \lesssim \| u \|_{L^\infty_T \big(\dot B^{2 \left( 1 -  \frac{1}{p}\right)}_{2, 2} \big)} \| c \|_{L^2_T(\dot B^1_{2, 2})} \lesssim \| \omega \|_{L^\infty_T(L^p)} \| c \|_{L^2_T(\dot H^1)}.
        \end{split}
    \end{equation*}
    Secondly, we also have
    \begin{equation*}
        \| \dot{\mc R}(u, c) \|_{L^2_T \big( \dot B^{2 \left( 1 - \frac{1}{p} \right)}_{\infty, 1} \big)} \lesssim \| u \|_{L^\infty_T(\dot B^{1 - \frac{2}{p}}_{\infty, \infty})} \| c \|_{L^2_T(\dot B^1_{\infty, 1})} \lesssim \| \omega \|_{L^\infty_T(L^p)} \| c \|_{L^2_T(\dot B^1_{\infty, 1})}.
    \end{equation*}
    Both estimates, used in Proposition \ref{p:Interpolation}, with interpolation parameter $\theta = 2 \left( 1 - \frac{1}{p} \right)$ provide
    \begin{equation}\label{eq:ConcIN3}
        \| \dot{\mc R}(u, c) \|_{\widetilde{L^2_T}(\dot B^0_{\infty, 1})} \lesssim \| \omega \|_{L^\infty_T(L^p)} \| c \|_{L^2_T(\dot H^1)}^{2 \left( 1 - \frac{1}{p} \right)} \| c \|_{L^2_T(\dot B^1_{\infty, 1})}^{-1+\frac{2}{p}}.
    \end{equation}
    Putting together estimates \eqref{eq:ConCIN1}, \eqref{eq:ConCIN2} and \eqref{eq:ConcIN3} in the Bony decomposition \eqref{eq:ConcBony1} gives the bound we were seeking for the convective term
    \begin{equation}\label{eq:ConcFinal2}
        \| u \cdot \nabla c \|_{\widetilde{L^2_T}(\dot B^{-1}_{\infty, 1})} \lesssim \| \omega \|_{L^\infty_T(L^p)} \| c \|_{L^2_T(\dot H^1)}^{2 \left( 1 - \frac{1}{p} \right)} \| c \|_{L^2_T(\dot B^1_{\infty, 1})}^{-1+\frac{2}{p}}.
    \end{equation}
    The philosophy of this bound is that if $p \rightarrow 1^+$, the concentration energy norm $\| c_0 \|_{L^2}$ almost not involved. This is consistent with scaling of the convective term and the fact that the space $(u, c) \in L^\infty(L^1) \times L^2(\dot B^1_{\infty, 1})$ is scale-critical. However, taking $p = 1$ would have required using a super-critical norm for the concentration, which would have been undesirable as $\omega_0 \in L^\infty$ is required anyway for uniqueness.

    \medskip

    \textbf{STEP 3: consumption term.} The term $nc$ is a bit different, as it is subcritical. However, we only have subcritical bounds for the micro-swimmer density at our disposal $n \in L^\infty_T(L^1) \cap L^3_T(L^3) \cap L^1_T(\dot W^{1, \frac{6}{5}})$, the strongest of which is $n \in L^3_T(L^3)$. Our goal is the Chemin-Lerner estimate $nc \in \widetilde{L^2_T}(\dot B^{-1}_{\infty, 1})$, which will be achieved through our interpolation result, Proposition \ref{p:Interpolation}. On the one hand, the Gagliardo-Nirenberg inequality (Proposition \ref{p:GN}) provides
    \begin{equation*}
        \| c \|_{L^6(L^3)} \lesssim \| c \|_{L^\infty(L^2)}^{2/3} \| c \|_{L^2(\dot H^1)}^{1/3} \, ,
    \end{equation*}
    so that we get a supercritical estimate for the product
    \begin{equation}\label{eq:ConcConsumptionEQ1}
        \begin{split}
            \| nc \|_{L^2_T(\dot B^{-4/3}_{\infty, \infty})} & \lesssim \| nc \|_{L^2_T(L^{3/2})} \\ 
            & \lesssim \| n \|_{L^3_T(L^3)} \| c \|_{L^6(L^3)} \\
            & \lesssim \| n \|_{L^3_T(L^3)} \| c \|_{L^\infty(L^2)}^{2/3} \| c \|_{L^2(\dot H^1)}^{1/3}.
        \end{split}
    \end{equation}
    On the other hand, the usual characterization of Besov spaces as real interpolation spaces (Proposition \ref{p:BesovInterpolate}) applied twice gives
    \begin{equation*}
        \begin{split}
            \| c \|_{L^6_T(L^\infty)} & \lesssim \| c \|_{L^6_T(\dot B^0_{\infty, 1})} \\
            & \lesssim \| c \|_{L^\infty (\dot B^{-1}_{\infty, \infty})}^{1/4} \| c \|_{L^6_T(\dot B^{1/3}_{\infty, \infty})}^{3/4} \\
            & \lesssim \| c \|_{L^\infty(L^2)}^{1/4} \Big( \| c \|_{L^\infty_T(\dot B^0_{\infty, \infty})}^{2/3} \| c \|_{L^2_T(\dot B^1_{\infty, \infty})}^{1/3} \Big)^{3/4} \\
            & \lesssim \| c \|_{L^\infty(L^2)}^{1/4} \| c \|_{L^\infty_T(\dot B^0_{\infty, \infty})}^{1/2} \| \nabla c \|_{L^2_T(L^\infty)}^{1/4}.
        \end{split}
    \end{equation*}
    Therefore, we have a subcritical bound on the product, namely
    \begin{equation}\label{eq:ConcConsumptionEQ2}
        \begin{split}
            \| nc \|_{L^2_T(\dot B^{-2/3}_{\infty, \infty})} & \lesssim \| nc \|_{L^2_T(L^3)} \lesssim \| n \|_{L^3_T(L^3)} \| c \|_{L^6_T(L^\infty)} \\
            & \lesssim \| n \|_{L^3_T(L^3)} \| c \|_{L^\infty(L^2)}^{1/4} \| c \|_{L^\infty_T(\dot B^0_{\infty, \infty})}^{1/2} \| \nabla c \|_{L^2_T(L^\infty)}^{1/4}.
        \end{split}
    \end{equation}
    We may use both estimates \eqref{eq:ConcConsumptionEQ1} and \eqref{eq:ConcConsumptionEQ2} in the interpolation inequality of Proposition \ref{p:Interpolation}, and obtain the desired bound
    \begin{equation}\label{eq:ConcFinal3}
        \begin{split}
            \| nc \|_{\widetilde{L^2_T}(\dot B^{-1}_{\infty, 1})} & \lesssim \| nc \|_{L^2_T(\dot B^{-4/3}_{\infty, \infty})}^{1/2} \| nc \|_{L^2_T(\dot B^{-2/3}_{\infty, \infty})}^{1/2} \\
            & \lesssim \| n \|_{L^3_T(L^3)} \| c \|_{L^\infty(L^2)}^{11/24} \| c \|_{L^2(\dot H^1)}^{1/6} \| c \|_{L^\infty_T(\dot B^0_{\infty, \infty})}^{1/4} \| \nabla c \|_{L^2_T(L^\infty)}^{1/8} \\
            & \lesssim \| n \|_{L^3_T(L^3)} \| c_0 \|_{L^2}^{5/8} \| c \|_{L^\infty_T(\dot B^0_{\infty, \infty})}^{1/4} \| \nabla c \|_{L^2_T(L^\infty)}^{1/8}.
        \end{split}
    \end{equation}

    \medskip

    \textbf{STEP 4: final estimate.} We gather \eqref{eq:ConcFinal1}, \eqref{eq:ConcFinal2} and \eqref{eq:ConcFinal3}. The heat estimates of Proposition \ref{p:heat} provide the following bound for the concentration. Note that $L^\infty_T(\dot B^0_{\infty, \infty}) = \widetilde{L^\infty_T}(\dot B^0_{\infty, \infty})$, and that, due to Minkowski's inequality, we have $\widetilde{L^2_T}(\dot B^1_{\infty, 1}) \subset L^2_T(\dot B^1_{\infty, 1})$. This provides
    \begin{equation*}
        \begin{split}
            \| c \|_{L^2_T(\dot B^1_{\infty, 1})} + \| c \|_{\widetilde{L^\infty_T}(\dot B^0_{\infty, 2})} & \lesssim \| c_0 \|_{\dot B^0_{\infty, 2}} \\
            & \qquad + \| \omega \|_{L^\infty_T(L^p)} \| c_0 \|_{L^2}^{2 \left( 1 - \frac{1}{p} \right)} \| c \|_{L^2_T(\dot B^1_{\infty, 1})}^{-1+\frac{2}{p}} \\
            & \qquad \qquad + \| n \|_{L^3_T(L^3)} \| c_0 \|_{L^2}^{5/8} \| c \|_{L^\infty_T(\dot B^0_{\infty, \infty})}^{1/4} \| c \|_{L^2_T(\dot B^1_{\infty, 1})}^{1/8}
        \end{split}
    \end{equation*}
    We apply Young's inequality to the products, in order to absorb the norm $\| c \|_{L^2_T(\dot B^1_{\infty, 1})}$ in the left-hand side. This gives our final concentration estimation:
    \begin{equation*}
        \| c \|_{L^2_T(\dot B^1_{\infty, 1})} + \| c \|_{\widetilde{L^\infty_T}(\dot B^0_{\infty, 2})} \lesssim \| c_0 \|_{\dot B^0_{\infty, 2}} + \| \omega \|_{L^\infty_T(L^p)}^{\frac{p'}{2}} \| c_0 \|_{L^2} + \| n \|_{L^3_T(L^3)}^{8/5} \| c_0 \|_{L^2}.
    \end{equation*}
\end{proof}

\subsection{Energy Estimates: Critical micro-swimmer Density}

Now that we have proved the critical concentration bound we need $\nabla c \in L^2_T(L^2)$, we are able to close the energy estimates for the micro-swimmer density equation. This will finally give us access to the critical norm $n \in L^\infty_T(L^2) \cap L^2_T(\dot H^1)$.

\begin{prop}\label{p:BactEN}
    Let $(n_0, c_0, u_0) \in L^1 \cap L^2 \times L^2 \cap \dot B^0_{\infty, 2} \times L^2 \cap \dot W^{1, \frac{6}{5}}$ be a set of initial data and $(n, c, u)$ a smooth solution of \eqref{ieq:EC}. Assume that thou bounds of Propositions \ref{p:Basic}, \ref{p:CriticalVelEN}, \ref{p:VortSuper} and \ref{p:ConcLip} hold, so that in particular $\nabla c \in L^2_T(L^\infty)$. We then have $n \in L^\infty_T(L^2) \cap L^2_T(\dot H^1) \cap L^4_T(L^4)$ with the inequality
    \begin{equation*}
        \frac{1}{2}\| n \|_{L^\infty_T(L^2)} + \frac{1}{2} \| n \|_{L^4_T(L^4)}^4 + \| n \|_{L^2_T(\dot H^1)}^2 \leq \frac{1}{2} \| n_0 \|_{L^2}^2 \exp \big( 2 T + \| \nabla c \|_{L^2_T(L^\infty)} \big).
    \end{equation*}
\end{prop}

\begin{proof}
    The proof is a simple energy inequality. Multiplying the micro-swimmer density equation by $n$ and integrating yields
    \begin{equation*}
        \begin{split}
            \frac{1}{2} \frac{\rm d}{\dt} \int n^2 \dx + \int |\nabla n|^2 \dx + \int n^4 \dx & = \int n^2 \dx + \int n \nabla n \cdot \nabla c \dx \\
            & \leq \| n \|_{L^2}^2 + \| n \|_{L^2} \| \nabla n \|_{L^2} \| \nabla c \|_{L^\infty} \\
            & \leq \| n \|_{L^2}^2 + \frac{1}{2} \| \nabla n \|_{L^2}^2 + \| n \|_{L^2}^2 \frac{1}{2} \| \nabla c \|_{L^\infty}^2.
        \end{split}
    \end{equation*}
    Grönwall's lemma then gives the result.
\end{proof}

\subsection{Final Bounds for the micro-swimmer Density and Vorticity}

We now have critical bounds for $n$. Using interpolation, we will deduce that $\nabla n \in L^1_T(L^\infty)$. This is step $(iii)$ in our proof strategy. Using one last time the subcritical gravity coupling term, this will show that $\omega \in L^\infty_T(L^\infty)$ (step $(iv)$ in \eqref{ieq:Diagram}).

\begin{prop}
    Let $(n_0, c_0, u_0) \in L^1 \cap L^2 \times L^2 \cap \dot B^0_{\infty, 2} \times L^2 \cap \dot W^{1, \frac{6}{5}}$ be a set of initial data and $(n, c, u)$ a smooth solution of \eqref{ieq:EC}. Assume that the bounds Propositions \ref{p:Basic}, \ref{p:CriticalVelEN}, \ref{p:VortSuper}, \ref{p:ConcLip} and \ref{p:BactEN} hold, so that in particular
    \begin{equation*}
        n \in L^\infty_T(L^2) \cap L^2_T(\dot H^1), \qquad \nabla c \in L^2_T(L^\infty) \qquad \text{and} \qquad \omega \in L^\infty_T(L^{6/5}).
    \end{equation*}
    Then $n \in \widetilde{L^\infty_T}(\dot B^{-1}_{\infty, 1}) \cap \widetilde{L^2_T}(\dot B^0_{\infty, 1}) \cap L^1_T(\dot B^1_{\infty, 1})$ with the following bound,
        \begin{equation*}
        \begin{split}
            \| n \|_{L^1_T(\dot B^1_{\infty, 1})} & \lesssim \| n_0 \|_{\dot B^{-1}_{\infty, 1}} \\
            & \qquad + \| n \|_{L^4_T(L^4)}^4 \| n \|_{L^\infty_T(L^2)} + \| n \|_{L^\infty_T(L^1)}^{1/2} \| n \|_{L^2_T(\dot H^1)}^{1/2} \\
            & \qquad \qquad + \| \omega \|_{L^\infty_T(L^{6/5})}^3 \| n \|_{L^2_T(\dot H^1)}  + \| c \|_{L^2_T(\dot B^1_{\infty, 1})}^2 \| n \|_{L^\infty_T(L^2)}.
        \end{split}
    \end{equation*}
    In particular, $\nabla n \in L^1_T(L^\infty)$. Under the additional assumption that $\omega_0 \in L^\infty$, we consequently have the additional estimate 
    \begin{equation*}
        \| \omega \|_{L^\infty_T(L^\infty)} \leq \| \omega_0 \|_{L^\infty} + \| \nabla n \|_{L^1_T(L^\infty)}.
    \end{equation*}
\end{prop}

\begin{proof}
    We apply maximal regularity estimates in Besov spaces (Proposition \ref{p:heat}) to the micro-swimmer density equation. Duhamel's formula reads
    \begin{equation*}
        n = e^{t \Delta} n_0 + \int_0^t e^{(t - s) \Delta} \Big( n - n^3 - u \cdot \nabla n + \D (n \nabla c)  \Big) \ds,
    \end{equation*}
    and we examine the terms one after the other. We will bound them in $L^1_T(\dot B^{-1}_{\infty, 1}) = \widetilde{L^1_T}(\dot B^{-1}_{\infty, 1})$. This equality will allow us to use ``standard'' interpolation (Proposition \ref{p:BesovInterpolate}) to get our bounds, as opposed as needing the Chemin-Lerner interpolation estimates provided by Proposition \ref{p:Interpolation}.

    \medskip

    \textbf{STEP 1: Reproduction term.} We start by looking at the reproduction term $n - n^3$. Notice that the term $n^3$, which we consider first, has critical scaling. We have, thanks to Proposition \ref{p:BesovInterpolate},
    \begin{equation*}
        \begin{split}
            \| n^3 \|_{\dot B^{-1}_{\infty, 1}} & \lesssim \| n^3 \|_{\dot B^{-3/2}_{\infty, \infty}}^{2/3} \| n^3 \|_{\dot B^0_{\infty, \infty}}^{1/3} \\
            & \lesssim \| n^3 \|_{L^{4/3}}^{2/3} \| n^3 \|_{L^\infty}^{1/3} \\
            & \lesssim \| n \|_{L^4}^2 \| n \|_{L^\infty}.
        \end{split}
    \end{equation*}
    On the other hand, using Proposition \ref{p:BesovInterpolate} one more time, we have
    \begin{equation*}
        \begin{split}
            \| n \|_{L^\infty} & \lesssim \| n \|_{\dot B^{-1}_{\infty, \infty}}^{1/2} \| n \|_{\dot B^1_{\infty, \infty}}^{1/2} \\
            & \lesssim \| n \|_{L^2}^{1/2} \| n \|_{\dot B^1_{\infty, 1}}^{1/2},
        \end{split}
    \end{equation*}
    so that one gets
    \begin{equation}\label{eq:BactEQ1}
        \| n^3 \|_{L^1_T(\dot B^{-1}_{\infty, 1})} \lesssim \| n \|_{L^4_T(L^4)}^2 \| n \|_{L^\infty_T(L^2)}^{1/2} \| n \|_{L^1_T(\dot B^1_{\infty, 1})}^{1/2}.
    \end{equation}
    For the linear part of the reproduction term, which is largely subcritical, we use the interpolation inequality of Proposition \ref{p:BesovInterpolate} again, which gives
    \begin{equation}\label{eq:BactEQ2}
        \begin{split}
            \| n \|_{L^1_T(\dot B^{-1}_{\infty, 1})} & \lesssim \| n \|_{L^4_T(\dot B^{-1}_{\infty, 1})} \\
            & \lesssim \| n \|_{L^\infty_T(\dot B^{-2}_{\infty, \infty})}^{1/2} \| n \|_{L^2_T(\dot B^{0}_{\infty, \infty})}^{1/2} \\
            & \lesssim \| n \|_{L^\infty_T(L^1)}^{1/2} \| n \|_{L^2_T(\dot H^1)}^{1/2}.
        \end{split}
    \end{equation}

    \medskip

    \textbf{STEP 2: Convection term.} Secondly, we bound the convective term $u \cdot \nabla n$. This term is critical with respect to scaling and contains a derivative of the micro-swimmer density. This might lead the reader to worry about the use of interpolation to bound it, as $\nabla n \in L^1_T(\dot B^0_{\infty, 1})$ is an ``endpoint'' estimate when compared to the bounds we have. However, our saving grace is that we have a supercritical estimate at our disposal for the velocity, namely $\omega \in L^\infty_T(L^{6/5})$, which well help as it did in the proof of the concentration estimate in Proposition \ref{p:ConcLip}. We write the Bony decomposition for the product
    \begin{equation*}
        u \cdot \nabla n = \D \dot{\mc T}_u(n) + \dot{\mc T}_{\nabla n} (u) + \D \dot{\mc R}(n, u).
    \end{equation*}
    First of all, notice that
    \begin{equation*}
        u \in W^{1, \frac{6}{5}} \subset \dot B^1_{\frac{6}{5}, \infty} \subset \dot B^{-2/3}_{\infty, \infty} \,,
    \end{equation*}
    so that, by Proposition \ref{p:Bony} and the interpolation property of Besov spaces (Proposition \ref{p:BesovInterpolate}), we bound the first paraproduct by
    \begin{equation*}
        \begin{split}
            \| \dot{\mc T}_u(\nabla n) \|_{\dot B^{-1}_{\infty, 1}} & \lesssim \| u \|_{\dot B^{-2/3}_{\infty, \infty}} \| \nabla n \|_{\dot B^{-1/3}_{\infty, 1}} \\
            & \lesssim \| \omega \|_{L^{6/5}} \| \nabla n \|_{\dot B^0_{\infty, 1}}^{2/3} \| \nabla n \|_{\dot B^{-1}_{\infty, \infty}}^{1/3}\\
            & \lesssim \| \omega \|_{L^{6/5}} \| n \|_{\dot B^1_{\infty, 1}}^{2/3} \| n \|_{\dot H^1}^{1/3}.
        \end{split}
    \end{equation*}
    This gives the first paraproduct bound
    \begin{equation}\label{eq:BactParaEQ1}
        \| \dot{\mc T}_u(\nabla n) \|_{L^1_T(\dot B^{-1}_{\infty, 1})} \lesssim \| \omega \|_{L^\infty_T(L^{6/5})} \| n \|_{L^1_T(\dot B^1_{\infty, 1})}^{2/3} \| n \|_{L^2_T(\dot H^1)}^{1/3}.
    \end{equation}
    For the second paraproduct, there is a slight difference. The function $u \in L^\infty_T(\dot B^{-2/3}_{\infty, \infty})$ has a third Besov index which is different from $1$ (even with sharper embeddings, we would only get $u \in L^\infty_T(\dot B^{-2/3}_{\infty, 2})$). In order to find the bound we seek, we therefore interpolate: Proposition \ref{p:BesovInterpolate} gives\footnote{We could also have used the critical bound $u \in L^\infty_T(L^2)$ and interpolation on the velocity to find a third Besov index equal to $1$. The method is of equivalent complexity.}
    \begin{equation*}
        \begin{split}
            \| \dot {\mc T}_{\nabla n}(u) \|_{\dot B^{-1}_{\infty, 1}} & \lesssim \| \dot {\mc T}_{\nabla n}(u) \|_{\dot B^{-5/3}_{\infty, \infty}}^{1/3}  \| \dot {\mc T}_{\nabla n}(u) \|_{\dot B^{-2/3}_{\infty, \infty}}^{2/3} \\
            & \lesssim \big( \| \nabla n\|_{\dot B^{-1}_{\infty, \infty}} \| u \|_{\dot B^{-2/3}_{\infty, \infty}} \big)^{1/3} \big( \| \nabla n\|_{\dot B^{0}_{\infty, 1}} \| u \|_{\dot B^{-2/3}_{\infty, \infty}} \big)^{2/3} \\
            & \lesssim \| \omega \|_{L^{6/5}} \| n \|_{\dot H^1}^{1/3} \| n \|_{\dot B^1_{\infty, 1}}^{2/3}.
        \end{split}
    \end{equation*}
    This gives the bound
    \begin{equation}\label{eq:BactParaEQ2}
        \| \dot {\mc T}_{\nabla n}(u) \|_{L^1_T(\dot B^{-1}_{\infty, 1})} \lesssim \| \omega \|_{L^\infty_T(L^{6/5})} \| n \|_{L^2_T(\dot H^1)}^{1/3} \| n \|_{L^1_T(\dot B^1_{\infty, 1})}^{2/3}.
    \end{equation}
    We also use Proposition \ref{p:BesovInterpolate} to interpolate between two estimations for the remainder term. Only, we must be careful to always have a positive sum of regularities to insure the remainder exists. We have on the one hand, with Proposition \ref{p:Bony},
    \begin{equation*}
        \| \dot{\mc R}(n, u) \|_{\dot B^{1/3}_{\infty, 1}} \lesssim \| u \|_{\dot B^{-2/3}_{\infty, \infty}} \| n \|_{\dot B^1_{\infty, 1}},
    \end{equation*}
    and on the other hand
    \begin{equation*}
        \| \dot{\mc R}(n, u) \|_{\dot B^{-2/3}_{\infty, \infty}} \lesssim \| \dot{\mc R}(n, u) \|_{\dot H^{1/3}} \lesssim \| u \|_{\dot B^{-2/3}_{\infty, \infty}} \| n \|_{\dot H^1}.
    \end{equation*}
    Therefore Proposition \ref{p:BesovInterpolate} gives the remainder estimate
    \begin{equation}\label{eq:BactParaEQ3}
        \| \dot {\mc R} (n, u) \|_{L^1_T(\dot B^0_{\infty, 1})} \lesssim \| \omega \|_{L^\infty_T(L^{6/5})} \| n \|_{L^2_T(\dot H^1)}^{1/3} \| n \|_{L^1_T(\dot B^1_{\infty, 1})}^{2/3}.
    \end{equation}
    Putting together estimates \eqref{eq:BactParaEQ1}, \eqref{eq:BactParaEQ2} and \eqref{eq:BactParaEQ3} gives the sought bound for the convective term, namely
    \begin{equation}\label{eq:BactEQ3}
        \| u \cdot \nabla n \|_{L^1_T(\dot B^{-1}_{\infty, 1})} \lesssim \| \omega \|_{L^\infty_T(L^{6/5})} \| n \|_{L^2_T(\dot H^1)}^{1/3} \| n \|_{L^1_T(\dot B^1_{\infty, 1})}^{2/3}.
    \end{equation}

    \medskip

    \medskip

    \textbf{STEP 3: Chemotactic term.} We examine the chemotactic term $\D (n \nabla c)$. Just as the the non-linear part of the reproduction term, it is critical. We use the Bony decomposition to write the product $n \nabla c$,
    \begin{equation*}
        n \nabla c = \dot{\mc T}_{\nabla c}(n) + \dot{\mc T}_n(\nabla c) + \dot{\mc R}(n, \nabla c).
    \end{equation*}
    The first paraproduct is bounded by (through usage of Propositions \ref{p:Bony} and Proposition \ref{p:BesovInterpolate})
    \begin{equation*}
        \begin{split}
            \| \dot{\mc T}_{\nabla c}(n) \|_{\dot B^0_{\infty, 1}} & \lesssim \| \nabla c \|_{L^\infty} \| n \|_{\dot B^0_{\infty, 1}} \\
            & \lesssim \| c \|_{\dot B^1_{\infty, 1}} \| n \|_{\dot B^{-1}_{\infty, \infty}}^{1/2} \| n \|_{\dot B^1_{\infty, \infty}} \\
            & \lesssim \| c \|_{\dot B^1_{\infty, 1}} \| n \|_{L^2}^{1/2} \| n \|_{\dot B^1_{\infty, 1}}^{1/2},
        \end{split}
    \end{equation*}
    so that 
    \begin{equation}\label{eq:BactParaEQ4}
        \| \dot{\mc T}_{\nabla c}(n) \|_{L^1_T(\dot B^0_{\infty, 1})} \lesssim \| c \|_{L^2_T(\dot B^1_{\infty, 1})} \| n \|_{L^\infty_T(L^2)}^{1/2} \| n \|_{L^1_T(\dot B^1_{\infty, 1})}^{1/2}.
    \end{equation}
    The second paraproduct can be bounded in a similar fashion:
    \begin{equation*}
        \begin{split}
            \| \dot{\mc T}_n(\nabla c) \|_{\dot B^0_{\infty, 1}} & \lesssim \| n \|_{L^\infty} \| \nabla c \|_{\dot B^0_{\infty, 1}} \\
            & \lesssim \| c \|_{\dot B^1_{\infty, 1}} \| n \|_{\dot B^0_{\infty, 1}} \\
            & \lesssim \| c \|_{\dot B^1_{\infty, 1}} \| n \|_{L^2}^{1/2} \| n \|_{\dot B^1_{\infty, 1}}^{1/2},
        \end{split}
    \end{equation*}
    and we have exactly the same bound
    \begin{equation}\label{eq:BactParaEQ5}
        \| \dot{\mc T}_{n}(\nabla c) \|_{L^1_T(\dot B^0_{\infty, 1})} \lesssim \| c \|_{L^2_T(\dot B^1_{\infty, 1})} \| n \|_{L^\infty_T(L^2)}^{1/2} \| n \|_{L^1_T(\dot B^1_{\infty, 1})}^{1/2}.
    \end{equation}
    We are left with the remainder term $\dot {\mc R}(n, \nabla c)$. By using both estimates $n \in L^\infty_T(L^2)$ and $n \in L^1_T(\dot B^1_{\infty, 1})$, we may write, with Propositions \ref{p:BesovInterpolate} and \ref{p:Bony},
    \begin{equation*}
        \begin{split}
            \| \dot {\mc R}(n, \nabla c) \|_{\dot B^0_{\infty, 1}} & \lesssim \| \dot {\mc R}(n, \nabla c) \|_{\dot B^{-1}_{\infty, \infty}}^{1/2} \| \dot {\mc R}(n, \nabla c) \|_{\dot B^1_{\infty, \infty}}^{1/2} \\
            & \lesssim \big( \| \nabla c \|_{\dot B^0_{\infty, 1}} \| n \|_{\dot B^0_{2, \infty}} \big)^{1/2} \big( \| \nabla c \|_{\dot B^0_{\infty, 1}} \| n \|_{\dot B^1_{\infty, 1}} \big)^{1/2} \\
            & \lesssim \| c \|_{\dot B^1_{\infty, 1}} \| n \|_{L^2}^{1/2} \| n \|_{\dot B^1_{\infty, 1}}^{1/2}.
        \end{split}
    \end{equation*}
    This gives an identical bound to \eqref{eq:BactParaEQ4} and \eqref{eq:BactParaEQ4}, namely
    \begin{equation}\label{eq:BactParaEQ6}
        \| \dot {\mc R}(n, \nabla c) \|_{L^1_T(\dot B^0_{\infty, 1})} \lesssim \| c \|_{L^2_T(\dot B^1_{\infty, 1})} \| n \|_{L^\infty_T(L^2)}^{1/2} \| n \|_{L^1_T(\dot B^1_{\infty, 1})}^{1/2}.
    \end{equation}
    Finally, putting together \eqref{eq:BactParaEQ4}, \eqref{eq:BactParaEQ5} and \eqref{eq:BactParaEQ6}, we find that
    \begin{equation}\label{eq:BactEQ4}
        \| \D (n \nabla c) \|_{L^1_T(\dot B^{-1}_{\infty, 1})} \lesssim \| c \|_{L^2_T(\dot B^1_{\infty, 1})} \| n \|_{L^\infty_T(L^2)}^{1/2} \| n \|_{L^1_T(\dot B^1_{\infty, 1})}^{1/2}.
    \end{equation}

    \medskip

    \textbf{STEP 4: Overall bound.} By putting together \eqref{eq:BactEQ1}, \eqref{eq:BactEQ2}, \eqref{eq:BactEQ3} and \eqref{eq:BactEQ4} in Duhamel's formula and resorting to Proposition \ref{p:heat}, we find that
    \begin{equation*}
        \begin{split}
            \| n \|_{L^1_T(\dot B^1_{\infty, 1})} & \lesssim \| n_0 \|_{\dot B^{-1}_{\infty, 1}} \\
            & \qquad + \| n \|_{L^4_T(L^4)}^2 \| n \|_{L^\infty_T(L^2)}^{1/2} \| n \|_{L^1_T(\dot B^1_{\infty, 1})}^{1/2} + \| n \|_{L^\infty_T(L^1)}^{1/2} \| n \|_{L^2_T(\dot H^1)}^{1/2} \\
            & \qquad \qquad + \| \omega \|_{L^\infty_T(L^{6/5})} \| n \|_{L^2_T(\dot H^1)}^{1/3} \| n \|_{L^1_T(\dot B^1_{\infty, 1})}^{2/3} \\
            & \qquad \qquad \qquad + \| c \|_{L^2_T(\dot B^1_{\infty, 1})} \| n \|_{L^\infty_T(L^2)}^{1/2} \| n \|_{L^1_T(\dot B^1_{\infty, 1})}^{1/2}.
        \end{split}
    \end{equation*}
    We may now use Young's inequality to absorb the fractional powers of the norm $\| n \|_{L^1_T(\dot B^1_{\infty, 1})}$ appearing in the right-hand side. This gives a final bound for $n$,
    \begin{equation*}
        \begin{split}
            \| n \|_{L^1_T(\dot B^1_{\infty, 1})} & \lesssim \| n_0 \|_{\dot B^{-1}_{\infty, 1}} \\
            & \qquad + \| n \|_{L^4_T(L^4)}^4 \| n \|_{L^\infty_T(L^2)} + \| n \|_{L^\infty_T(L^1)}^{1/2} \| n \|_{L^2_T(\dot H^1)}^{1/2} \\
            & \qquad \qquad + \| \omega \|_{L^\infty_T(L^{6/5})}^3 \| n \|_{L^2_T(\dot H^1)}  + \| c \|_{L^2_T(\dot B^1_{\infty, 1})}^2 \| n \|_{L^\infty_T(L^2)}.
        \end{split}
    \end{equation*}
    It only remains to see that plugging the weaker estimate $\nabla n \in L^1_T(L^\infty)$ in the vorticity equation also gives the Yudovich-type bound we were looking for,
    \begin{equation*}
        \| \omega \|_{L^\infty_T(L^\infty)} \leq \| \omega_0 \|_{L^\infty} + \int_0^T \| \nabla n \|_{L^\infty} \dt.
    \end{equation*}
    This concludes the proof.
\end{proof}

\section{Uniqueness of Solutions}\label{s:Unique}

In this Section, we prove stability estimates for the solution. This will be used both to show uniqueness and to prove the convergence of approximate solutions, and hence obtain existence. The spaces involved in Theorem \ref{t:Stability} are below critical regularity (but this will be enough). A critical result will be proved in Section \ref{s:CritUn}.

\begin{thm}\label{t:Stability}
    Consider two sets of initial data $(n_{j, 0}, c_{j, 0}, u_{j, 0})_{j \in \{ 1, 2\}}$ as in \eqref{ieq:OurData}, and two associated solutions $(n_j, c_j, u_j)$ of \eqref{ieq:EC}. We introduce the differences
    \begin{equation*}
        \begin{split}
            & (\delta n_0, \delta c_0, \delta u_0) = (n_{2, 0} - n_{1, 0}, c_{2, 0} - c_{1, 0}, u_{2, 0} - u_{1, 0}) ,\\
            & (\delta n, \delta c, \delta u) = (n_2 - n_1, c_2 - c_1, u_2 - u_1).
        \end{split}
    \end{equation*}
    Then, we have the following stability estimate
    \begin{multline*}
        \| \delta c \|_{L^\infty_T(L^2) \cap L^2_T(\dot H^1 \cap \dot B^1_{\frac{2}{s}, 2})} + \| \delta n \|_{L^\infty_T(\dot H^{-s}) \cap L^2_T(\dot H^{1 - s})} + \| \delta u \|_{L^\infty_T(L^2)} \\
        \lesssim \Big( \| \delta c_0 \|_{L^2} + \| \delta c_0 \|_{\dot B^0_{\frac{2}{s}, 2}} + \| \delta n_0 \|_{\dot H^{-s}} + \| \delta u_0 \|_{L^2} \Big)^{e^{-KT}},
    \end{multline*}
    where $K > 0$ and the implicit constants only depend on $T > 0$ and on the norms of the solutions in the spaces \eqref{ieq:Bounds1}.
\end{thm}

\begin{proof}
        The proof is based on energy estimates on the difference equations. By taking the difference between the systems solved by $(n_j, c_j, u_j)$, we obtain
    \begin{equation}\label{eq:UniqueDiff}
        \begin{cases}
            \partial_t \delta u + (u_1 \cdot \nabla) \delta u + (\delta u \cdot \nabla) u_2 +\nabla \delta \pi = g\,\delta n\\
            \partial_t \delta n - \Delta \delta n = - \delta u \cdot \nabla n_2 - u_1 \cdot \nabla \delta n + \delta n - (n_2^2 + n_1 n_2 + n_1^2) \delta n - \D (\delta n \nabla c_2 + n_1 \nabla \delta c) \\
            \partial_t \delta c - \Delta \delta c = - \delta u \cdot \nabla c_2 - u_1 \cdot \nabla \delta c - \delta n c_2 - \delta c n_2.
        \end{cases}
    \end{equation}
    In the above, $\delta \pi$ denotes the difference between the pressures associated to both solutions. 
    
    \medskip
    
    Before performing our estimates, we discuss the functional setting these will take place in. First of all, because the vorticity is bounded $\omega_1 \in L^\infty_T(L^\infty)$, we will implement Yudovich's argument for uniqueness in the 2D Euler equations, which is based on $L^\infty_T(L^2)$ stability estimates for the velocity. Notice that this is a critical space. For $\delta c$ and $\delta n$, we will start by establishing subcritical stability estimates, which are easier, and the next section will be devoted to stability estimates in critical spaces. 
    
    A natural choice would be an energy space $L^\infty_T(L^2) \cap L^2_T(\dot H^1)$, however this space is critical for the micro-swimmer density $n$, although it is subcritical for the chemical concentration. We will therefore estimate $\delta c$ in the energy space $L^\infty_T(L^2) \cap L^2_T(\dot H^1)$, while we bound $\delta n$ in a subcritical space $L^\infty_T(\dot H^{-s}) \cap L^2_T(\dot H^{1 - s})$, where the exponent fulfills $0 < s < 1$. The bound $s > 0$ serves to make this space subcritical, while the bound $s < 1$ is simply due to the fact that our class of initial micro-organism densities do not lie $\dot H^{-1}$, as
    \begin{equation*}
        n_0 \in L^1 \cap L^2 \cap \dot B^{-1}_{\infty, 1} \not \subset \dot H^{-1},
    \end{equation*}
    whereas for $L^1$ we only have the weaker embedding $L^1 \subset \dot B^{-1}_{2, \infty}$.

    As a last remark, we point out that owing to the fact that $\delta n_0 \in \dot H^{-s}$ has $s$ less derivatives than critical regularity $n_0 \in L^2$, we will also have to bound $\delta c_0$ in a space with corresponding scaling. Since $c_0 \in L^2 \cap \dot B^0_{\infty, 2}$, Proposition \ref{p:ComplexBesovInterpolate} shows that an appropriate space is $\delta c_0 \in \dot B^0_{\frac{2}{s}, 2} \subset \dot B^{-s}_{\infty, 2}$.

    \medskip

    \textbf{STEP 1: $\delta c$ Estimate.} We start by performing energy estimates in the chemical concentration difference equation. By multiplying the last equation in \eqref{eq:UniqueDiff} by $\delta c$ and integrating, we obtain
    \begin{equation}\label{eq:UnSobEQ1}
        \frac{1}{2} \frac{\rm d}{\dt} \int |\delta c|^2 + \int |\nabla \delta c|^2 + \int n_2 |\delta c|^2 \leq \int |\delta u||\delta c||\nabla c_2| + \left| \int c_2 \delta n \delta c \, \right|.
    \end{equation}
    By applying Young's inequality, we obtain an upper bound for the first integral, namely
    \begin{equation}\label{eq:UnSobEQ2}
        \begin{split}
            \int |\delta u||\delta c||\nabla c_2| & \leq  \| \delta u \|_{L^2} \| \delta c \|_{L^2} \| \nabla c_2 \|_{L^\infty} \\
            & \leq \| \delta u \|_{L^2}^2 + \| \delta c \|_{L^2}^2 \| \nabla c_2 \|_{L^\infty}^2.
        \end{split}
    \end{equation}
    On the other hand, we treat the integral $\int c_2 \delta c \delta n$ as a $\dot H^{s} \times \dot H^{-s}$ duality bracket. This gives
    \begin{equation*}
        \left| \int c_2 \delta c \delta n \, \right| \lesssim \| \delta n \|_{\dot H^{-s}} \| c_2 \delta c \|_{\dot H^{s}}.
    \end{equation*}
    To finish our first step, we have to bound the product $c \delta c$ in the factional Sobolev space $\dot H^s$. For this, we will use the Bony decomposition of the product,
    \begin{equation*}
        c \delta c = \dot{\mc T}_{\delta c} (c) + \dot{\mc T}_c (\delta c) + \dot{\mc R}(c, \delta c).
    \end{equation*}
    For the first paraproduct, we directly use Proposition \ref{p:Bony}, to get
    \begin{equation*}
        \| \dot{\mc T}_{\delta c} (c) \|_{L^2_T(\dot H^s)} \lesssim \| \delta c \|_{L^\infty_T(L^2)} \| c_2 \|_{L^2_T(\dot B^s_{\infty, 2})}.
    \end{equation*}
    The $\dot B_{\infty, 2}$ norm can be bounded through interpolation between $\dot H^1 \subset \dot B^0_{\infty, 2}$ and $\dot B^1_{\infty, \infty}$. This provides a first inequality
    \begin{equation*}
        \| \dot{\mc T}_{\delta c} (c) \|_{L^2_T(\dot H^s)} \lesssim \| \delta c \|_{L^\infty_T(L^2)} \| c_2 \|_{L^2_T(\dot H^1 \cap \dot B^1_{\infty, \infty})}.
    \end{equation*}
    For the second paraproduct $\dot{\mc T}_c(\delta c)$, we start with the interpolation inequality $\| \delta c \|_{L^{2/s}_T(\dot H^s)} \lesssim \| \delta c \|_{L^\infty_T(L^2)}^{1 - s} \| \delta c \|_{L^2_T(\dot H^1)}^s$. If we can bound the concentration $c_2$ in the space $L^{2/(1 - s)}_T( L^\infty)$, then Proposition \ref{p:Bony} will immediately give the desired estimate. If $s \leq \frac{1}{2}$ (\textit{i.e.} when $\frac{2}{1-s} \leq 4$), then it suffice to notice that, through straightforward interpolation
    \begin{equation}\label{eq:UnSobEQ18}
        \begin{split}
            \| c_2 \|_{L^{2/(1 - s)}_T(L^\infty)} \leq \| c_2 \|_{L^{4}_T(\dot B^0_{\infty, 1})} & \lesssim \| c_2 \|_{L^\infty_T(\dot B^{-1}_{\infty, 2})}^{1/2} \| c_2 \|_{L^2_T(\dot B^1_{\infty, \infty})}^{1/2} \\
            & \lesssim \| c_2 \|_{L^\infty_T(L^2)}^{1/2} \| c_2 \|_{L^2_T(\dot B^1_{\infty, \infty})}^{1/2}.
        \end{split}
    \end{equation}
    In the case where $s > \frac{1}{2}$, we instead use two steps of interpolation to obtain the desired spacetime regularity. Note that $1 - \frac{1}{s} \in ]-1, 0[$. We have,
    \begin{equation}\label{eq:UnSobEQ19}
        \begin{split}
            \| c_2 \|_{L^{2/(1 - s)}_T(L^\infty)} \leq \| c_2 \|_{L^{2/(1 - s)}_T(\dot B^0_{\infty, 1})} & \lesssim \| c_2 \|_{L^2_T(\dot B^1_{\infty, \infty})}^{1 - s} \| c_2 \|_{L^\infty_T(\dot B^{1 - 1/s}_{\infty, \infty})}^s \\
            & \lesssim \| c_2 \|_{L^2_T(\dot B^1_{\infty, \infty})}^{1 - s} \| c_2 \|_{L^\infty_T(\dot B^{-1}_{\infty, 2})}^{s \left( \frac{1}{s} - 1 \right)} \| c_2 \|_{L^\infty_T(\dot B^{0}_{\infty, \infty})}^{s \left(2 - \frac{1}{s} \right)} \\
            & \lesssim \| c_2 \|_{L^2_T(\dot B^1_{\infty, \infty})}^{1 - s} \| c_2 \|_{L^\infty_T(L^2)}^{1 - s} \| c_2 \|_{L^\infty_T(\dot B^{0}_{\infty, \infty})}^{2s - 1}.
        \end{split}
    \end{equation}
    This provides the paraproduct bound
    \begin{equation*}
        \| \dot{\mc T}_{c_2}(\delta c) \|_{L^2_T(\dot H^s)} \lesssim \| \delta c \|_{L^\infty_T(L^2)}^{1 - s} \| \delta c \|_{L^2_T(\dot H^1)}^s \| \nabla c_2 \|_{L^2_T(L^\infty)}^\theta \| c_2 \|_{L^\infty_T(L^2 \cap \dot B^0_{\infty, \infty})}^{1 - \theta},
    \end{equation*}
    where $\theta = \min\{ \frac{1}{2}, 1 - s \} \in ]0, 1[$ is chosen so as to cover both cases. Finally, we are left with the remainder term $\dot {\mc R} (c_2, \delta c)$. Proposition \ref{p:Bony} directly provides, together with the bounds we have just proven
    \begin{equation*}
        \begin{split}
            \| \dot {\mc R} (c_2, \delta c) \|_{L^2_T(\dot H^s)} & \lesssim \| c_2 \|_{L^{2/(1 - s)}_T(\dot B^0_{\infty, \infty})} \| \delta c \|_{L^{2/s}_T(\dot H^s)} \\
            & \lesssim \| \delta c \|_{L^\infty_T(L^2)}^{1 - s} \| \delta c \|_{L^2_T(\dot H^1)}^s \| \nabla c_2 \|_{L^2_T(L^\infty)}^\theta \| c_2 \|_{L^\infty_T(L^2 \cap \dot B^0_{\infty, \infty})}^{1 - \theta} \, ,
        \end{split}
    \end{equation*}
    where $\theta$ is the same exponent as above. Putting together all these three bounds gives the product estimate
    \begin{equation*}
        \begin{split}
            \|c_2 \delta c \|_{L^2_T(\dot H^s)} & \lesssim \| \delta c \|_{L^\infty_T(L^2)} \| \nabla c_2 \|_{L^2_T(L^2 \cap L^\infty)} \\
             & \qquad \qquad + \| \delta c \|_{L^\infty_T(L^2)}^{1 - s} \| \delta c \|_{L^2_T(\dot H^1)}^s \| \nabla c_2 \|_{L^2_T(L^\infty)}^\theta \| c_2 \|_{L^\infty_T(L^2 \cap \dot B^0_{\infty, \infty})}^{1 - \theta} \,.
        \end{split}
    \end{equation*}
    Now, we use this bound and \eqref{eq:UnSobEQ2} in the differential inequality \eqref{eq:UnSobEQ1}, which we integrate with respect to time. This gives the bound
    \begin{equation*}
        \begin{split}
            \sup_{[0, T]} \| \delta c \|_{L^2}^2 + \frac{1}{2} \| \delta c \|_{L^2_T(\dot H^1)}^2 & \lesssim \| \delta c_0 \|_{L^2}^2 + \int_0^T \| \delta u \|_{L^2}^2 \| \nabla c_2 \|_{L^\infty}^2 \dt + \| \delta n \|_{L^2_T(\dot H^{-s})} \| c_2 \delta c \|_{L^2_T(\dot H^s)} \\
            & \lesssim \| \delta c_0 \|_{L^2}^2 + \int_0^T \big( \| \delta u \|_{L^2}^2 + \| \delta c \|_{L^2}^2 \| \nabla c_2 \|_{L^\infty}^2 \big) \dt + \int_0^T \| \delta n \|_{\dot H^{-s}}^2 \dt \\
            & \qquad \qquad\qquad+ \int_0^T \| \delta c \|_{L^2}^{2(1 - s)} \| \delta c \|_{\dot H^1}^{2s} \| \nabla c_2 \|_{L^2 \cap L^\infty}^{2\theta} \| c_2 \|_{L^2 \cap \dot B^0_{\infty, \infty}}^{2(1 - \theta)} \dt.
        \end{split}
    \end{equation*}
    To deal with this last integral, we apply Young's inequality to absorb the norm $\| \delta c \|_{\dot H^1}$ in the left-hand side. This provides the following inequality, with exponents depending on $\theta$.
    \begin{equation}\label{eq:UnSobEQ3}
        \begin{split}
            \| \delta c \|_{L^\infty_T(L^2)}^2 + \frac{1}{4} \| \delta c \|_{L^2_T(\dot H^1)}^2 & \lesssim \| \delta c_0 \|_{L^2}^2 + \int_0^T \big( \| \delta u \|_{L^2}^2 + \| \delta c \|_{L^2}^2 \| \nabla c_2 \|_{L^\infty}^2 + \| \delta n \|_{\dot H^{-s}}^2 \big) \dt \\
            & \qquad \qquad \qquad+ \int_0^T \| \delta c \|_{L^2}^2 \| \nabla c_2 \|_{L^2 \cap L^\infty}^{\frac{2 \theta}{1 - s}} \| c_2 \|_{L^2 \cap \dot B^0_{\infty, \infty}}^{2\frac{(1 - \theta)}{1 - s}} \dt.
        \end{split}
    \end{equation}
    Observe that, since we have set $\theta = \min \{ \frac{1}{2}, 1 - s\}$, the fraction $\theta/(1 - s)$ is always smaller than one. As the concentration satisfies $c_2 \in L^2_T( \dot H^1 \cap \dot W^{1, \infty}) \cap L^\infty_T(L^2 \cap \dot B^0_{\infty, \infty})$, this means that the factor of $\| \delta c \|_{L^2}^2$ in the last integral is always a $L^1$ function of time, and so this integral is a Grönwall term. Estimate \eqref{eq:UnSobEQ3} will give a concentration difference estimate, provided we can find $L^\infty_T(\dot H^{-s})$ bounds on $\delta n$.

    \medskip

    \textbf{STEP 2: $\delta n$ fractional Sobolev estimates.} We seek to bound the micro-organism difference function $\delta n$ in the space $L^\infty_T(\dot H^{-s}) \cap L^2_T(\dot H^{1 - s})$. For this, we multiply the $\delta n$ equation in \eqref{eq:UniqueDiff} by  $(- \Delta)^{-s} \delta n$  and integrate. We have, after integrating by parts,
    \begin{equation}\label{eq:UnSobEQ4}
        \begin{split}
            \frac{1}{2} \frac{\rm d}{\dt} \| \delta n \|_{\dot H^{-s}}^2 + \| \delta n  \|_{\dot H^{1 - s}}^2 & = \| \delta n \|_{\dot H^{-s}}^2 - \int (- \Delta)^{-s} \delta n (n_1^2 + n_1 n_2 + n_2^2) \delta n \\
            & \qquad \qquad+ \int \nabla (- \Delta)^{-s} \delta n \cdot \Big( \delta u \, n_2 + u_1 \delta n + \delta n \nabla c_2 + n_1 \nabla \delta c \Big)
        \end{split}
    \end{equation}
    Before proceeding, we remark that Proposition \ref{p:BesovInterpolate} gives the following interpolation inequalities, which we will use several times. We have:
    \begin{equation}\label{eq:UnSobEQ5LL}
        \begin{split}
            & \| \delta n \|_{L^{2/s}_T(L^2)} \lesssim \| \delta n \|_{L^\infty_T(\dot H^{-s})}^{1 - s} \| \delta n \|_{L^2_T(\dot H^{1 - s})}^s \\
            & \| \nabla (- \Delta)^{-s} \delta n \|_{L^{2/(1 - s)}_T(L^2)} \lesssim \| \nabla (- \Delta)^{-s} \delta n \|_{L^\infty_T(\dot H^{s - 1})}^s \| \nabla (- \Delta)^{-s} \delta n \|_{L^2_T(\dot H^s)}^{1 - s} \\
            & \qquad \qquad \qquad \qquad \qquad \lesssim \|\delta n\|_{L^\infty_T(\dot H^{-s})}^s \| \delta n \|_{L^2_T(\dot H^{1 - s})}^{1 - s} \, .
        \end{split}
    \end{equation}
    We have five terms to bound. We start with the last term in the first line of \eqref{eq:UnSobEQ4} just above. Our starting point is the Sobolev embedding $(- \Delta)^{-s} \delta n \in L^\infty_T(\dot H^{s}) \subset L^\infty_T(L^{2/(1 - s)})$. We wish to use this in combination with the first inequality in \eqref{eq:UnSobEQ5LL}, which means that we should estimate the quadratic quantities $n_k n_\ell$ (where $k, \ell \in \{1, 2\}$) in $L^{2/(1 - s)}_T(L^{2/s})$. Hölder's inequality then gives
    \begin{equation*}
        \begin{split}
            \left| \int (- \Delta)^{-s} \delta n \, n_k n_\ell \delta n \right| & \leq \big\| (- \Delta)^{-s} \delta n \big\|_{L^{2/(1 - s)}} \| \delta n \|_{L^2} \big\| n_k n_\ell \big\|_{L^{2/s}} \\
                & \lesssim \| \delta n \|_{\dot H^{-s}}^{2 - s} \| \delta n \|_{\dot H^{1 - s}}^s \big\| n_k n_\ell \big\|_{L^1}^{s/2} \big\| n_k n_\ell \big\|_{L^\infty}^{\frac{2-s}{2}},
        \end{split}
    \end{equation*}
    and we can use Young's inequality to write, for any $\epsilon > 0$,
    \begin{equation}\label{eq:UnSobEQ13}
        \left| \int (- \Delta)^{-s} \delta n \, n_k n_\ell \delta n \right| \lesssim \epsilon \| \delta n \|_{\dot H^{1- s}}^2 + C(\epsilon) \| \delta n \|_{\dot H^{-s}}^2 \| n \|_{L^2}^{\frac{2s}{2 - s}} \| n \|_{L^\infty}^2,
    \end{equation}
    where we have simply denoted $n$ the generic element $n_k$ for $k \in \{1, 2\}$. Remark that, because the solution satisfies $n \in L^\infty_T(L^2) \cap L^2_T(L^\infty)$, the factor of $\| \delta n \|_{\dot H^{-s}}^2$ is integrable with respect to time. This integral will be a Grönwall term.

    \medskip

    We continue with the next term, involving the product $\delta u n_2$. Resorting again to \eqref{eq:UnSobEQ5LL}, we see that 
    \begin{equation*}
        \begin{split}
            \left| \int \nabla (- \Delta)^{-s} \delta n \cdot \delta u \, n_2 \right| & \leq \| \nabla (- \Delta)^{-s} \delta n \|_{L^2} \| \delta u \|_{L^2} \| n \|_{L^\infty} \\
            & \lesssim \| \delta n \|_{\dot H^{-s}}^s \| \delta n  \|_{\dot H^{1 - s}}^{1 - s} \| \delta u \|_{L^2} \| n_2 \|_{L^\infty}.
        \end{split}
    \end{equation*}
    Using Young's inequality, we may first isolate the $\| \delta u \|_{L^2}$ norm, and then absorb the high order norm $\| \delta n \|_{\dot H^{1 - s}}$ in the left-hand side. This gives, for any $\epsilon > 0$,
        \begin{equation}\label{eq:UnSobEQ14}
        \begin{split}
            \left| \int \nabla (- \Delta)^{-s} \delta n \cdot \delta u \, n_2 \right| & \lesssim \| \delta u \|_{L^2}^2 \| n_2 \|_{L^\infty}^2 + \| \delta n \|_{\dot H^{-s}}^{2s} \| \delta n \|_{\dot H^{1 - s}}^{2(1 - s)} \\
            & \lesssim \| \delta u \|_{L^2}^2 \| n_2 \|_{L^\infty}^2 + \epsilon \| \delta n \|_{\dot H^{1 - s}}^2 + C(\epsilon) \| \delta n \|_{\dot H^{-s}}^2.
        \end{split}
    \end{equation}

    The newt integral has the terms $\nabla (- \Delta)^{-s} \delta n \cdot u_1 \delta n$. Here, the inequalities in \eqref{eq:UnSobEQ5LL} directly provide a bound, since $u_1 \in L^\infty_T(L^2 \cap \dot B^1_{\infty, \infty}) \subset L^\infty_T(L^\infty)$. We have:
    \begin{equation}\label{eq:UnSobEQ17}
        \begin{split}
            \left| \int \nabla (- \Delta)^{-s} \delta n \cdot u_1 \delta n \right| & \leq \| \nabla (- \Delta)^{-s} \delta n \|_{L^2} \| u_1 \|_{L^\infty} \| \delta n \|_{L^2} \\
            & \leq \| \delta n \|_{\dot H^{-s}} \| \delta n \|_{\dot H^{1 - s}} \| u_1 \|_{L^\infty} \\
            & \leq \epsilon \| \delta n \|_{\dot H^{1 - s}}^2 + \| \delta n \|_{\dot H^{-s}}^2 \| u_1 \|_{L^\infty}^2.
        \end{split}
    \end{equation}

    \medskip

    We examine the integral involving $\nabla (- \Delta)^{-s} \delta n \cdot \delta n \nabla c$. Resorting again to \eqref{eq:UnSobEQ5LL}, we proceed similarly to get
    \begin{equation}\label{eq:UnSobEQ15}
        \begin{split}
            \left| \int \nabla (- \Delta)^{-s} \delta n \cdot \delta n \nabla c_2 \right| & \leq \| \nabla (- \Delta)^{-s} \delta n \|_{L^2} \| \delta n \|_{L^2} \| \nabla c_2 \|_{L^\infty} \\
            & \lesssim \| \delta n \|_{\dot H^{-s}} \| \delta n \|_{\dot H^{1 - s}} \| \nabla c_2 \|_{L^\infty} \\
            & \lesssim \epsilon \| \delta n \|_{\dot H^{1 - s}}^2 + C(\epsilon) \| \delta n \|_{\dot H^{-s}}^2 \| \nabla c_2 \|_{L^\infty}^2.
        \end{split}
    \end{equation}

    \medskip

    Finally, we are left with the integral involving $\nabla \delta c$, which is the most difficult of the lot. This is because the bounds \eqref{eq:UnSobEQ3} we have for $\delta c$ are of lower scaling than the Sobolev norm $\delta n \in L^\infty_T(\dot H^{-s}) \cap L^2_T(\dot H^{1 - s})$, while the chemotaxis term is critical. To conclude, we will have to resort to a higher order bound, which will be proved in the third step, namely
    \begin{equation}\label{eq:UnSobEQ6}
        \delta c \in \widetilde{L^\infty_T}(\dot B^0_{\frac{2}{s}, 2}) \cap L^2_T(\dot B^1_{\frac{2}{s}, 2}).
    \end{equation}
    There are two reasons for wanting to prove such a Besov bound. The first one is that it is compatible with the initial datum space: Proposition \ref{p:ComplexBesovInterpolate} provides the embedding $\delta c_0 \in \dot B^0_{2/s, 2} \subset L^2 \cap \dot B^0_{\infty, 2}$. On the other hand, the above implies that $\nabla \delta c \in L^2_T(L^{2/s})$ by Proposition \ref{p:BesovFine}, which is very convenient to bound the last integral. Indeed, using \eqref{eq:UnSobEQ5LL}, we have
    \begin{equation}\label{eq:UnSobEQ16}
        \begin{split}
            \left| \int \nabla (- \Delta)^{-s} \delta n \cdot n_1 \nabla \delta c \right| & \leq \| \nabla (- \Delta)^{-s} \delta n \|_{L^2} \| n_1 \|_{L^{2/(1 - s)}} \| \nabla \delta c \|_{L^{2/s}} \\
            & \lesssim \| \delta n \|_{\dot H^{-s}}^s \| \delta n \|_{\dot H^{1 - s}}^{1 - s} \| n_1 \|_{L^2}^{1 - s} \| n_1 \|_{L^\infty}^{s} \| \delta c \|_{\dot B^1_{\frac{2}{s}, 2}} \\
            & \lesssim \epsilon \| \delta n \|_{\dot H^{1 - s}}^2 + C(\epsilon) \| \delta n \|_{\dot H^{-s}}^2 \| n_1 \|_{L^2}^{2 \left( \frac{1}{s} - 1 \right)} \| n_1 \|_{L^\infty}^2 + \| \nabla \delta c \|_{L^{2/s}}^2
        \end{split}
    \end{equation}
    Putting all the bounds \eqref{eq:UnSobEQ13}, \eqref{eq:UnSobEQ14}, \eqref{eq:UnSobEQ17}, \eqref{eq:UnSobEQ15} and \eqref{eq:UnSobEQ16} together in \eqref{eq:UnSobEQ4}, we have, after integration in time, an integral estimate for the fractional Sobolev norms, namely
    \begin{equation}\label{eq:UnSObEQ22}
        \begin{split}
            & \| \delta n \|_{L^\infty_T(\dot H^{-s})}^2 + \| \delta n \|_{L^2_T(\dot H^{1 - s})} \lesssim \| \delta n_0 \|_{\dot H^{-s}}^2 \\ 
            & \qquad + C(\epsilon) \int_0^T \| \delta n \|_{\dot H^{-s}}^2 \Big( \| n \|_{L^2}^{\frac{2s}{2-s}} \| n \|_{L^\infty}^2 + 1 + \| u_1 \|_{L^\infty}^2 + \| \nabla c_2 \|_{L^\infty}^2 + \| n_1 \|_{L^2}^{2 \left( \frac{1}{s} - 1 \right)} \| n_1\|_{L^\infty}^2 \Big) \dt \\
            & \qquad \qquad + \int_0^T \| \delta u \|_{L^2}^2 \| n_2 \|_{L^\infty}^2 \dt + \| \nabla \delta c \|_{L^2_T(L^{2/s})}
        \end{split}
    \end{equation}
    
    We now have to bound $\delta c$ in an appropriate space, which must be a subspace of $L^2_T(\dot W^{1, \frac{2}{s}})$. 

    \medskip

    \textbf{STEP 3: $\delta c$ Besov estimates.} For convenience of notation, we set $p = \frac{2}{s}$. To prove estimates of the form \eqref{eq:UnSobEQ6}, we will use Proposition \ref{p:heat}. We rewrite the $\delta c$ equation using Duhamel's formula:
    \begin{equation}\label{eq:UnSobEQ11}
        \delta c = e^{t \Delta} \delta c_0 - \int_0^t e^{(t - s) \Delta} \Big\{ \delta u \cdot \nabla c_2 + u_1 \cdot \nabla \delta c + c_1 \delta n + n_2 \delta c \Big\} \ds.
    \end{equation}
    From Proposition \ref{p:heat}, we see that finding $L^2_T(\dot B^{-1}_{p, 2})$ bounds for the terms in the brackets is enough to obtain \eqref{eq:UnSobEQ6}. Once again, we have to examine the different terms one after the other. We start by $\delta u \cdot \nabla c_2$. We will obtain our bound by interpolation: remark that we have the inequalities
    \begin{equation*}
        \begin{split}
            & \| \delta u \cdot \nabla c_2 \|_{L^2_T(\dot B^{-1 + \frac{2}{p}}_{p, 2})} \lesssim \| \delta u \cdot \nabla c_2 \|_{L^2_T(L^2)} \lesssim \| \delta u \|_{L^\infty_T(L^2)} \| \nabla c_2 \|_{L^2_T(L^\infty)} \\
            &  \| \delta u \cdot \nabla c_2 \|_{L^2_T(\dot B^{-2 + \frac{2}{p}}_{p, 2})} \lesssim \| \delta u \cdot \nabla c_2 \|_{L^2_T(L^1)} \lesssim \| \delta u \|_{L^\infty_T(L^2)} \| \nabla c_2 \|_{L^2_T(L^2)}.
        \end{split}
    \end{equation*}
    Therefore, interpolating with Proposition \ref{p:BesovInterpolate} provides 
    \begin{equation}\label{eq:UnSobEQ7}
        \| \delta u \cdot \nabla c_2 \|_{L^2_T(\dot B^{-1}_{p, 2})}^2 = \int_0^T \| \delta u \cdot \nabla c_2 \|_{\dot B^{-1}_{p, 2}}^2 \dt \lesssim \int_0^T \| \delta u \|_{L^2}^2 \| c_2 \|_{\dot H^1}^{4/p} \| c_2 \|_{\dot W^{1, \infty}}^{2 \left( 1 - \frac{2}{p} \right)} \dt  
    \end{equation}
    Secondly, for the term $u_1 \cdot \nabla \delta c$, we use the fact that $u_1$ is divergence free to write $u_1 \cdot \nabla \delta c = \D (u_1 \delta c)$, so that
    \begin{equation*}
        \begin{split}
            & \| u_1 \cdot \nabla \delta c \|_{L^2_T(\dot H^{-1})} = \| u_1 \delta c \|_{L^2_T(L^2)} \lesssim \| u_1 \|_{L^\infty_T(L^2)} \| \nabla \delta c \|_{L^2_T(\dot H^1)}\\
            & \| u_1 \cdot \nabla \delta c \|_{L^2_T(L^2)} \lesssim \| u_1 \|_{L^\infty_T(L^\infty)} \| \nabla \delta c \|_{L^2_T(\dot H^1)},
        \end{split}
    \end{equation*}
    Interpolating between $L^2$ and $\dot H^{-1}$ gives a $\dot H^{-s}$ bound, and the embedding $\dot H^{-s} \subset \dot B^{-1}_{\frac{2}{s}, 2} = \dot B^{-1}_{p, 2}$ in turn provides the bound we were looking for, namely
    \begin{equation}\label{eq:UnSobEQ8}
        \begin{split}
            \| u_1 \cdot \nabla \delta c \|_{L^2_T(\dot B^{-1}_{p, 2})}^2 & = \int_0^T \| u_1 \cdot \nabla \delta c \|_{\dot B^{-1}_{p, 2}}^2 \dt \lesssim \int_0^T \| u_1 \|_{L^\infty}^2 \| \delta c \|_{L^2}^{2s} \| \delta c \|_{\dot H^1}^{2(1 - s)} \dt \\
             & \lesssim \epsilon \| \delta c \|_{L^2_T(\dot H^1)}^2 + C(\epsilon) \int_0^T \| u_1 \|_{L^\infty}^2 \| \delta c \|_{L^2}^2 \dt.
        \end{split}
    \end{equation}
    For the term $n_2 \delta c$, we proceed as in \eqref{eq:UnSobEQ7}, interpolating between $L^1$ and $L^2$. On the one hand, we have the inequalities
    \begin{equation*}
        \begin{split}
            & \| n_2 \delta c \|_{L^2_T(L^2)} \lesssim \| n_2 \|_{L^4_T(L^4)} \| \delta c \|_{L^4_T(L^4)} \lesssim \| n_2 \|_{L^\infty_T(L^2)}^{1/2} \| n_2 \|_{L^2_T(\dot H^1)}^{1/2} \| \delta c \|_{L^\infty_T(L^2)}^{1/2} \| \delta c \|_{L^2_T(\dot H^1)}^{1/2} \\
            & \| n_2 \delta c \|_{L^\infty_T(L^1)} \leq \| n_2 \|_{L^\infty_T(L^2)} \| \delta c \|_{L^\infty_T(L^2)}.
        \end{split}
    \end{equation*}
    By applying the same interpolation inequalities as above, we have a $L^{2p/(p-2)}_T(\dot B^{-1}_{p, 2})$ estimate for $n_2 \delta c$. Now, since $p > 2$, we have $\frac{2p}{p-2} > 2$, and this bound is stronger than our target estimate:
    \begin{equation}\label{eq:UnSobEQ10}
        \begin{split}
            \int_0^T \| n_2 \delta c \|_{\dot B^{-1}_{p, 2}}^2 \dt  &\lesssim \int_0^T \| n_2 \|_{L^2}^{2 \left( \frac{1}{2} + \frac{1}{p} \right)} \| n_2 \|_{\dot H^1}^{2 \left( \frac{1}{2} - \frac{1}{p} \right)} \| \delta c \|_{L^2}^{2 \left( \frac{1}{2} + \frac{1}{p} \right)} \| \delta c \|_{\dot H^1}^{2 \left( \frac{1}{2} - \frac{1}{p} \right)} \dt \\
            & \lesssim \epsilon \| \delta c \|_{L^2_T(\dot H^1)}^2 + C(\epsilon) \int_0^T \| \delta c \|_{L^2}^2 \| n_2 \|_{L^2}^2 \| n_2 \|_{\dot H^1}^{2 \frac{p-2}{p+2}}\dt 
        \end{split}
    \end{equation}
    Finally, we have to deal with the term $c_1 \delta n$. Even though $\delta n$ has been estimated in fractional Sobolev spaces, the first inequality in \eqref{eq:UnSobEQ5LL} will allow us to proceed exactly as above. Therefore, we note that
    \begin{equation*}
        \begin{split}
            & \| c_1 \delta n \|_{L^1} \leq \| c_1 \|_{L^2} \| \delta n \|_{L^2}\\
            & \| c_1 \delta n \|_{L^2} \leq \| c_1 \|_{L^\infty} \| \delta n \|_{L^2}.
        \end{split}
    \end{equation*}
    This provides, with the help of \eqref{eq:UnSobEQ5LL}, the inequality
    \begin{equation*}
        \begin{split}
            \| c_1 \delta n \|_{\dot B^{-1}_{p, 2}} & \lesssim \| c_1 \|_{L^2}^{2/p} \| c_1 \|_{L^\infty}^{1 - 2/p} \| \delta n \|_{L^2} \\
            & \lesssim \| c_1 \|_{L^2}^{2/p} \| c_1 \|_{L^\infty}^{1 - 2/p} \| \delta n \|_{\dot H^{-s}}^{1 - s} \| \delta n \|_{\dot H^{1 - s}}^s.
        \end{split}
    \end{equation*}
    Therefore, integrating with respect to time and remembering that $\frac{2}{p} = s$, we write
    \begin{equation}\label{eq:UnSobEQ9}
        \begin{split}
            \int_0^T \| c_1 \delta n \|_{\dot B^{-1}_{p, 2}}^2 \dt & \lesssim \int_0^T \| c_1 \|_{L^2}^{2s} \| c_1 \|_{L^\infty}^{2(1 - s)} \| \delta n \|_{\dot H^{-s}}^{2(1 - s)} \| \delta n \|_{\dot H^{1 - s}}^{2s} \dt \\
            & \lesssim\epsilon \| \delta n \|_{\dot H^{1 -s}}^2 + C(\epsilon) \int_0^T \| \delta n \|_{\dot H^{-s}}^2 \| c_1 \|_{L^2}^2 \| c_1 \|_{L^\infty}^{\frac{2s}{1-s}} \dt.
        \end{split}
    \end{equation}
    Before concluding the Besov estimates, we will simply check that this last integral is a Grönwall term. This means that we have to make sure that $c_1 \in L^{2s/(1 - s)}_T(L^\infty)$ for any $0 < s < 1$. Recalling the estimates we have made in \eqref{eq:UnSobEQ18} and \eqref{eq:UnSobEQ19}, we see that 
    \begin{equation*}
        \| c_1 \|_{L^{2s/(1 - s)}_T(L^\infty)} \lesssim \| c_1 \|_{L^{2s/(1-s)}_T(L^\infty)} \lesssim \| c_1 \|_{L^2_T(\dot B^1_{\infty, 1})}^\theta \| c_1 \|_{L^\infty_T(L^2 \cap \dot B^0_{\infty, \infty})}^{1 - \theta},
    \end{equation*}
    with $\theta = \min \{ \frac{1}{2}, 1 - s \}$. This insures that the factor is indeed integrable, and that Grönwall's lemma will be usable.

    \medskip

    We now collect all estimates \eqref{eq:UnSobEQ7}, \eqref{eq:UnSobEQ8}, \eqref{eq:UnSobEQ10} and \eqref{eq:UnSobEQ9} in Duhamel's formula \eqref{eq:UnSobEQ11}. Proposition \ref{p:heat} then gives the estimate we were looking for:
    \begin{equation}\label{eq:UnSobEQ12}
        \begin{split}
            \| \delta c \|_{\widetilde{L^\infty_T}(\dot B^0_{p, 2})}^2 + \| \delta c \|_{L^2_T(\dot B^1_{p, 2})}^2 & \lesssim \| \delta c \|_{\dot B^0_{p, 2}}^2 + \epsilon \| \delta c \|_{L^2_T(\dot H^1)}^2 + \epsilon \| \delta n \|_{L^2_T(\dot H^{1 - s})}^2 \\
            & \qquad+ \int_0^T \| \delta u \|_{L^2} \| c_2 \|_{L^2}^{2s} \| c_2 \|_{\dot W^{1, \infty}}^{2(1 - s)} \dt \\
            & \qquad \qquad  + C(\epsilon) \int_0^T \| u_1 \|_{L^\infty}^2 \| \delta c \|_{L^2}^2 \dt \\
            & \qquad \qquad \qquad + C(\epsilon) \int_0^T \| \delta c \|_{L^2}^2 \| n_2 \|_{L^2}^2 \| n_2 \|_{\dot H^1}^{2 \frac{1-s}{1+s}} \dt \\
            & \qquad \qquad \qquad \qquad + C(\epsilon) \int_0^T \| \delta n \|_{\dot H^{-s}}^2 \| c_1 \|_{L^2}^2 \| c_1 \|_{L^\infty}^{2s/(1 - s)} \dt.
        \end{split}
    \end{equation}

    \medskip

    \textbf{STEP 4: Yudovich Argument.} Following the famous papers of Yudovich \cite{Yudovich, YudovichEN}, this part of the argument is based on an energy estimate in the Euler equation. My multiplying the $\delta u$ equation in \eqref{eq:UniqueDiff} by $\delta u$ and integrating, we find that
    \begin{equation}\label{eq:UnSobEQ20}
        \frac{1}{2} \frac{\rm d}{\dt} \| \delta u \|_{L^2}^2 \leq \int |\delta u|^2 |\nabla u_2| + \int |\delta u| |\delta n|,
    \end{equation}
    where the pressure and convective terms have vanished thanks to integration by parts. On the one hand, the last integral above is bounded by using \eqref{eq:UnSobEQ5LL}. We have
    \begin{equation*}
        \begin{split}
            \int |\delta u||\delta  n| & \leq \| \delta u \|_{L^2} \| \delta n \|_{L^2} \\
            & \lesssim \| \delta u \|_{L^2}\| \delta n \|_{\dot H^{-s}}^{1-s} \| \delta n \|_{\dot H^{1-s}}^s \\
            & \lesssim \| \delta u \|_{L^2} + \epsilon \| \delta n \|_{\dot H^{1-s}}^2 + C(\epsilon) \| \delta n \|_{\dot H^{-s}}^2.
        \end{split}
    \end{equation*}
    We look at the first integral in the right-hand side of \eqref{eq:UnSobEQ20}. In principle, to estimate the velocity difference function in $\delta u \in L^\infty_T(L^2)$, we would write a bound of the form
    \begin{equation*}
        \int |\delta u|^2 |\nabla u_2| \leq \| \delta u \|_{L^2}^2 \| \nabla u_2 \|_{L^\infty}.
    \end{equation*}
    However, it is not known that $\nabla u_2 \in L^1_T(L^\infty)$. The best bound we have at our disposal is $\omega_2 \in L^\infty_T(L^\infty)$, which leads to the inferior estimate $\nabla u_2 \in L^\infty_T({\rm BMO})$. The idea of Yudovich is to use the fact that ${\rm BMO} \subset L^p_{\rm loc}$ for every $p < + \infty$, and say that while it is not true that $\nabla u_2$ is $L^\infty$, it is \textit{almost} true. Quantitatively, we have, by the Calder\'on-Zygmund estimates of Proposition \ref{p:CZ},
    \begin{equation}\label{eq:UniqueCZ}
        \forall p \in \Big] \frac{6}{5}, + \infty \Big[, \qquad \| \nabla u_2 \|_{L^p} \leq C \, p \, \| \omega_2 \|_{L^p} \leq C \, p \, \| \omega_2 \|_{L^{6/5} \cap L^\infty}.
    \end{equation}
    We will perform estimates while using this bound with $p$ taken very large. Before diving in the computations, we agree that whenever $q \in [1, +\infty]$, $q'$ is the conjugate exponent, such that $\frac{1}{q} + \frac{1}{q'} = 1$. We have, to start with, by Hölder's inequality
    \begin{equation*}
        \int |\delta u|^2 |\nabla u_2| \leq \| \nabla u_2 \|_{L^p} \big\| |\delta u|^2 \big\|_{L^{p'}} \leq \| \nabla u_2 \|_{L^p} \| \delta u \|_{L^{2p'}}^2.
    \end{equation*}
    We may bound the $L^{2p'}$ norm of $\delta u$ by interpolation between $L^2$ (the norm we want to bound), and $L^\infty$. This provides
    \begin{equation*}
        \begin{split}
            \| \delta u \|_{L^{2p'}} & \leq \| \delta u \|_{L^2}^{1/p'} \| \delta u \|_{L^\infty}^{1/p} \\
            & \leq \| \delta u \|_{L^2}^{1 - \frac{1}{p}} \big( \| u_1 \|_{L^\infty} + \| u_2 \|_{L^\infty} \big)^{1/p}
        \end{split}
    \end{equation*}
    Further using Proposition \ref{p:BesovInterpolate} gives a bound of the $L^\infty$ norm by the velocity: for $u = u_1$ or $u = u_2$, we have
    \begin{equation*}
        \| u \|_{L^\infty} \leq \| u \|_{\dot B^0_{\infty, 1}} \lesssim \| u \|_{\dot B^{-2/3}_{\infty, \infty}}^{3/5} \| u \|_{\dot B^1_{\infty, \infty}}^{2/5} \lesssim \| \omega \|_{L^{6/5}}^{3/5} \| \omega \|_{L^\infty}^{2/5}.
    \end{equation*}
    Consequently, by plugging this and the Calder\'on-Zygmund estimate of \eqref{eq:UniqueCZ}, we have
    \begin{equation}\label{eq:UniqueEQ4}
        \begin{split}
            \int |\delta u|^2 |\nabla u_2| & \lesssim p \| \delta u \|_{L^2}^{2 \left( 1 - \frac{1}{p} \right)} \Big( \| \omega_1 \|_{L^{6/5}} + \| \omega_2 \|_{L^\infty} \Big)^{1 + \frac{2}{p}} \\
            & \lesssim p \| \delta u \|_{L^2}^{2\left( 1 - \frac{1}{p} \right)} K(T)^{1 + \frac{2}{p}}.
        \end{split}
    \end{equation}
    In order to obtain explicit estimates, we will optimize the choice of $p$. We fix
    \begin{equation*}
        p := \log \left( \frac{8 K(T)^2}{\| \delta u \|_{L^2}^2} \right).
    \end{equation*}
    Note that since $\| \delta u \|_{L^2} \leq 2 K(T)$, by definition of $K(T)$, we have $p \geq \log(4) > \frac{6}{5}$, so our choice of $p$ is indeed admissible. By plugging this in \eqref{eq:UniqueEQ4}, we obtain a log-Lipschitz estimate for the the integral
    \begin{equation*}
        \int |\delta u|^2 |\nabla u_2| \lesssim \| \delta u \|_{L^2}^2 K(T) \log \left( \frac{8 K(T)^2}{\| \delta u \|_{L^2}^2} \right). 
    \end{equation*}
    This is reassuring, as $r \mapsto r \big(1 - \log(r)\big)$ is an Osgood modulus of continuity, so we expect no problem to come from this term in the stability estimates. Overall, concerning the velocity difference function, we have the bound
    \begin{equation}\label{eq:UnSobEQ21}
        \| \delta u \|_{L^2}^2 \lesssim \| \delta u_0 \|_{L^2}^2 + \epsilon \| \delta n \|_{\dot H^{1-s}}^2 + \int_0^T \| \delta n \|_{\dot H^{-s}}^2 \dt + \int_0^T \| \delta u \|_{L^2}^2 K(T) \log \left( \frac{8 K(T)^2}{\| \delta u \|_{L^2}^2} \right) \dt .
    \end{equation}
    
    \medskip

    \textbf{STEP 5: Summary.} It is now time to put all of our estimates together. Given that we have numerous equations, we will adopt the following notation for better readability: in the sequel, $M_1(t)$ denotes a generic $L^1([0, T])$ function depending only on $T$ and the solutions $(n_j, c_j, u_j)$ with $j \in \{ 1, 2 \}$. Gathering the bounds \eqref{eq:UnSobEQ3}, \eqref{eq:UnSObEQ22}, \eqref{eq:UnSobEQ12} and \eqref{eq:UnSobEQ21}, we have the final estimate
    \begin{equation}\label{eq:UnSobEQ23}
        \begin{split}
            & \| \delta c \|_{L^\infty_T(L^2)}^2 + \| \delta c \|_{L^2_T(\dot H^1)}^2 + \| \delta c \|_{L^2_T(\dot B^1_{\frac{2}{s}, 2})}^2 + \| \delta n \|_{L^\infty_T(\dot H^{-s})}^2 + \| \delta n \|_{L^2_T(\dot H^{1-s})}^2 + \| \delta u \|_{L^\infty_T(L^2)}^2 \\
            & \qquad \lesssim \| \delta c \|_{L^2}^2 + \| \delta c \|_{\dot B^0_{\frac{2}{s}, 2}}^2 + \| \delta n \|_{\dot H^{-s}}^2 + \| \delta u \|_{L^2}^2 \\
            & \qquad \qquad +\int_0^T M_1(t) \Big( \| \delta u \|_{L^2}^2 + \| \delta n \|_{\dot H^{1-s}}^2 + \| \delta c \|_{L^2}^2 \Big) \dt + \int_0^T \| \delta u \|_{L^2}^2 K(T) \log \left( \frac{8 K(T)^2}{\| \delta u \|_{L^2}^2} \right) \dt
        \end{split}
    \end{equation}
    This is encouraging, as $r \mapsto (1-r) \log(r)$ is an Osgood modulus of continuity, allowing us to conclude thanks to a Grönwall-type argument. We will replicate the computations in the proof of Lemma 3.4 in \cite{BCD} to extract our final estimate from the above. We define the quantity
    \begin{equation*}
        E(t) := \sup_{[0, t]} \Big[ \| \delta n \|_{\dot H^{-s}} + \| \delta c \|_{L^2} + \| \delta u \|_{L^2} \Big]
    \end{equation*}
    as well as the functions
    \begin{equation*}
        \begin{split}
            & \mu(r) := \frac{r}{K(T)} \left(1 - \log\left(\frac{K(T)}{r} \right) \right) \qquad \text{for } r \leq K(T),\\
            & M(r) = \int_r^{K(T)} \frac{1}{\mu} = \frac{1}{K} \log \left[ 1 + \log \left( \frac{K(T)}{r} \right) \right] + {\rm Cst}.
        \end{split}
    \end{equation*}
    Now, remark that, up to adding non-negative terms to our estimate \eqref{eq:UnSobEQ23}, we may rewrite it as
    \begin{equation*}
        E(T) \leq E(0) + \int_0^T \mu \big( E(t) \big) \dt := Q(T).
    \end{equation*}
    The right-hand side defines a function $Q(t)$ whose derivative is is simply $\mu(E(t))$ and which satisfies, by construction, $E(t) \leq Q(t)$. In particular, since $\mu$ is a non-decreasing function, we have $Q'(t) \leq \mu(E(t)) \leq \mu(Q(t))$, and therefore
    \begin{equation*}
        - \frac{\rm d}{{\rm d}T} M \big( Q(T) \big) = \frac{Q'(T)}{\mu (Q(T))} \leq 1.
    \end{equation*}
    By integrating with respect to time, we have
    \begin{equation*}
        - M\big( Q(T) \big) + M \big( Q(0) \big) \leq T,
    \end{equation*}
    which gives, after exponentiating twice
    \begin{equation*}
        Q(T) \lesssim Q(0)^{e^{-KT}},
    \end{equation*}
    and, remembering that $E(T) \leq Q(T)$, we have our final estimate
    \begin{multline*}
        \| \delta c \|_{L^\infty_T(L^2) \cap L^2_T(\dot H^1 \cap \dot B^1_{\frac{2}{s}, 2})} + \| \delta n \|_{L^\infty_T(\dot H^{-s}) \cap L^2_T(\dot H^{1 - s})} + \| \delta u \|_{L^\infty_T(L^2)} \\
            \lesssim \Big( \| \delta c_0 \|_{L^2} + \| \delta c_0 \|_{\dot B^0_{\frac{2}{s}, 2}} + \| \delta n_0 \|_{\dot H^{-s}} + \| \delta u_0 \|_{L^2} \Big)^{e^{-KT}},
    \end{multline*}
    where the implicit constants depend on $K(T)$. This concludes the proof of the stability estimates.
\end{proof}

\section{Existence and Uniqueness of Solutions}\label{s:Existence}

In this short paragraph, we explain how the existence of solutions may be obtained by means of an iterative scheme. Since the method is extremely standard, we only give a summary of the proof, only providing details for the difficult part, which is the convergence of approximate solutions to a strong solution. We refer for example to \cite{LZ2026} for such a construction on a nearly identical system.

\begin{prop}
    Under the assumptions \eqref{ieq:OurData} for the initial data $(n_0, c_0, u_0)$, there exists a unique solution (in the sense of distributions) $(n, c, u)$ of \eqref{ieq:EC} that satisfies all the bounds \eqref{ieq:Bounds1} and \eqref{ieq:Bounds2}.
\end{prop}

\begin{proof}
    The first step is the construction of a family of approximate solutions $(n_k, c_k, u_k)$, by means of an iterative scheme. We consider the following system
    \begin{equation}\label{eq:ExistenceEQ1}
        \begin{cases}
            \partial_t n_{k+1} - \Delta n_{k+1} = n_{k+1}(1 - n_k^2) - u_k \cdot \nabla n_{k+1} - \D (n_k \nabla c_k) \\
            \partial_t c_{k+1} - \Delta c_{k+1} = - u_k \cdot \nabla c_{k+1} - n_k c_k \\
            \partial_t u_{k+1} + (u_k \cdot \nabla) u_{k+1} + \nabla \pi_{k+1} = n_k g
        \end{cases}
    \end{equation}
    with regularized initial data $(S_k n_0, S_k c_0, S_k u_0)$. Solutions for this regularized system can in turn be constructed by using, for example, a spectral projection and a hyper-viscous term $\Delta^N n_{k+1}$. A straightforward adaptation of the \textit{a priori} above show that the sequence $(n_k, c_k, u_k)$ is bounded in the Banach space
    \begin{equation*}
        \begin{split}
            & (n_k) \in L^\infty_T(L^1 \cap L^2 \cap \dot B^{-1}_{\infty, 1}) \cap L^2_T(\dot H^1 \cap \dot B^0_{\infty, 1}) \cap L^1_T(L^1) \\
            & (c_k) \in L^\infty(L^2 \cap \dot B^0_{\infty, \infty}) \cap L^2_T(\dot H^1 \cap \dot B^1_{\infty, 1}) \\
            & (u_k) \in L^\infty_T(L^2 \cap W^{1, \frac{6}{5}}), \quad \omega \in L^\infty_T(L^\infty).
        \end{split}
    \end{equation*}
    Likewise, the stability estimates show that we have the following strong convergences: for any $0 < s < 1$, as $k \rightarrow + \infty$,
    \begin{equation}\label{eq:ExistenceEQ2}
        \begin{split}
            & n_k \longrightarrow n \qquad \text{in } L^\infty_T(\dot H^{-s}) \cap L^2_T(\dot H^{1 - s}) \\
            & c_k \longrightarrow c \qquad \,\, \text{in } L^\infty_T(L^2 \cap \dot B^0_{\frac{2}{s}, 2}) \cap L^2_T(\dot H^1 \cap \dot B^1_{\frac{2}{s}, 2}) \\
            & u_k \longrightarrow u \qquad \, \text{in } L^\infty_T(L^2).
        \end{split}
    \end{equation}
    Using the convergences above, we show that the different terms in \eqref{eq:ExistenceEQ1} converge, so that the limit $(n, c, u)$ solves the PDE problem. Only two terms require special attention, the first one being the cubic term $n_{k+1}n_k^2$. Using the fact that \eqref{eq:ExistenceEQ2} holds for any $0 < s < 1$, we take in particular $s = 1/3$ and use interpolation: the embeddings
    \begin{equation*}
        L^\infty_T(\dot H^{-1/3}) \cap L^2_T(\dot H^{2/3}) \subset L^3_T(\dot H^{1/3}) \subset L^3_T(L^3)
    \end{equation*}
    then prove that the convergence $n_{k+1}n_k^2 \longrightarrow n^3$ holds in the space $L^1_T(L^1)$. The second part which needs some care is the chemotaxis term $\D(n_k \nabla c_k)$, which is the non-linear term with the highest order derivatives. The first two lines in \eqref{eq:ExistenceEQ2}, along with \eqref{eq:UnSobEQ5LL}, show that
    \begin{equation*}
        \begin{split}
            & n_k \longrightarrow n \qquad \quad \,\,\text{in } L^{2/s}_T(L^2) \\
            & \nabla c_k \longrightarrow \nabla c \qquad \text{in } L^2_T(L^{2/s}).
        \end{split}
    \end{equation*}
    Therefore, we have strong convergence $n_k \nabla c_k \longrightarrow n \nabla c$ in the space $L^{2/(s+1)}_T(L^{2/(s+1)})$, and therefore distributional convergence of $\D(n_k \nabla c_k)$ to $\D(n\nabla c)$. For the remaining non-linear terms, we have
    \begin{equation*}
        \begin{split}
            & u_k n_k  \longrightarrow un \qquad \quad\qquad \,\,\, \text{in } L^{2/s}_T(L^1)\\
            & u_k c_{k+1} \longrightarrow uc \qquad \quad \qquad \text{in } L^\infty_T(L^1) \\
            & n_k c_k \longrightarrow nc \qquad \quad \qquad \;\;\,\, \text{in } L^{2/s}_T(L^1)\\
            & u_k \otimes u_{k+1} \longrightarrow u \otimes u \,\,\qquad \text{in } L^\infty_T(L^1).
        \end{split}
    \end{equation*}
    The last line also insures the convergence of the pressures, since they are given by Leray projection $\pi_{k+1} = (- \Delta)^{-1} \nabla^2 : (u_k \otimes u_{k+1})$ and by virtue of the non-homogeneous embedding $L^1 \subset H^{-1-\epsilon}$. We therefore have
    \begin{equation*}
        \nabla \pi_{k+1} \longrightarrow \nabla \pi \qquad \text{in } L^\infty_T(H^{-2-\epsilon}).
    \end{equation*}
    All these convergences put together in the iterative scheme are enough to show that system \eqref{eq:ExistenceEQ1} converges to the PDE we are studying \eqref{ieq:EC}. Uniqueness follows from the stability estimates.
\end{proof}

\section{Flow Map Continuity in Critical Spaces}\label{s:CritUn}

\begin{thm}\label{t:t:StabilityCritical}
    Under the assumptions and notation of Theorem \ref{t:Stability}, we have the following stability estimate:
    \begin{equation*}
        \begin{split}
            \| \delta n \|_{L^\infty_T(L^2)}^2 + \| \delta n \|_{L^2_T(\dot H^1)}^2 & \lesssim \| \delta n_0 \|_{L^2}^2 + \| \delta c_0 \|_{\dot B^0_{\infty, 2}} \\
             & \qquad + F \left( \left[ \| \delta c_0 \|_{L^2} + \| \delta c_0 \|_{\dot B^0_{\frac{2}{s}, 2}} + \| \delta n_0 \|_{\dot H^{-s}} + \| \delta u_0 \|_{L^2} \right]^{e^{-KT}} \right).
        \end{split}
    \end{equation*}
    In the above, $F(r) = r \log(1 + \frac{K}{r})$.
    \end{thm}

\begin{proof}
    To prove critical estimates, we simply have to find an energy estimate $\delta n \in L^\infty_T(L^2) \cap L^2_T(\dot H^1)$ and a Besov bound $\delta c \in L^2_T(\dot B^1_{\infty, 1})$ (recall that the space $\delta u \in L^\infty_T(L^2)$ is already scaling critical).

    \medskip

    \textbf{STEP 1: Micro-swimmer density bound.} We multiply the $\delta n$ equation in \eqref{eq:UniqueDiff} by $\delta n$ and integrate. This provides
    \begin{equation}\label{eq:UnCrEQ6}
        \begin{split}
            & \frac{1}{2} \frac{\rm d}{\dt} \int |\delta n|^2 + \int |\nabla \delta n|^2 + \int (n_1^2 + n_1 n_2 + n_2^2)|\delta n|^2 \\
            & \qquad\qquad \leq \int |\delta n|^2 + \int |\delta u| |\delta n| |\nabla n_2| + \int |\delta n| |\nabla \delta n| |\nabla c_2| + \int n_1 |\nabla \delta n| |\nabla \delta c|.
        \end{split}
    \end{equation}
    The first three integrals in the right-hand side are fairly easy to bound. We may write
    \begin{equation}\label{eq:UnCrEQ7}
        \begin{split}
            & \int |\delta u| |\delta n| |\nabla n_2| \leq \| \delta u \|_{L^2} \| \delta n \|_{L^2} \| \nabla n_2 \|_{L^\infty} \leq \big( \| \delta u \|_{L^2}^2 + \| \delta n \|_{L^2}^2 \big) \| \nabla n \|_{L^\infty}\\
            & \int |\delta n| |\nabla \delta n| |\nabla c_2| \leq \| \delta n\|_{L^2} \| \nabla \delta n \|_{L^2} \| \nabla c_2 \|_{L^\infty} \leq \frac{1}{4} \| \nabla \delta n \|_{L^2}^2 + \| \delta n \|_{L^2}^2 \| \nabla c \|_{L^\infty}^2.\\
        \end{split}
    \end{equation}
    However, to bound the last integral requires a $L^2_T(\dot W^{1, \infty})$ estimate on the concentration difference function $\delta c$. This requires examining the $\delta c$ equation separately.

    \medskip

    \textbf{STEP 2: critical $\delta c$ bound.} We use Duhamel's formula for the $\delta c$ equation from \eqref{eq:UnSobEQ11}. In order to find $L^2_T(\dot B^1_{\infty, 1})$ estimates on $\delta c$, we seek to bound the terms in the brackets in $\widetilde{L^2_T}(\dot B^{-1}_{\infty, 1})$. Since we work in comparable spaces as those of the previous stability estimates, we might have expected many of the bounds to be the same. However, the fact that $\widetilde{L^2_T}(\dot B^{-1}_{\infty, 1})$ is an endpoint space means that additional work is required. As an illustration, the embedding $\dot B^0_{p, 2} \subset L^p$ does not hold for $p = +\infty$ (we instead have $\dot B^0_{\infty, 2} \subset {\rm BMO}$).

    \medskip

    First of all, thanks to Proposition \ref{p:AGHeat}, we have
    \begin{equation}\label{eq:UnCrEQ1}
        \| e^{t \Delta} \delta c_0 \|_{L^2_T(\dot B^{1}_{\infty, 1})} + \| e^{t \Delta} \delta c_0 \|_{\widetilde{L^\infty_T}(\dot B^{0}_{\infty, 2})} \lesssim \| \delta c_0 \|_{\dot B^0_{\infty, 2}}.
    \end{equation}
    In \eqref{eq:UnSobEQ7} we leave the term $\delta u \cdot \nabla c_2$ for later, as it involves more technical computations. We examine the convective term $u_1 \cdot \nabla \delta c$ first. To bound this, we can simply resort to inequality  \eqref{eq:ConcFinal2} obtained in STEP 2, where the concentration $c$ is replaced by the difference function $\delta c$, namely
    \begin{equation}\label{eq:UnCrEQ3}
        \| u_1 \cdot \nabla \delta c \|_{\widetilde{L^2_T}(\dot B^{-1}_{\infty, 1})} \lesssim \| \omega_1 \|_{L^\infty_T(L^{6/5})} \| \delta c \|_{L^2_T(\dot H^1)}^{1/3} \| \delta c \|_{L^2_T(\dot B^1_{\infty, 1})}^{2/3}
    \end{equation}
    The next terms we have to estimate are $n_2 \delta c$ and $c_1 \delta n$. For the first of these two, we notice that
    \begin{equation*}
        \begin{split}
            & \| n_2 \delta c\|_{L^\infty_T(\dot B^{-2}_{\infty, \infty})} \lesssim \| n_2 \delta c\|_{L^\infty_T(L^1)} \leq \| n_2 \|_{L^\infty_T(L^2)} \| \delta c\|_{L^\infty_T(L^2)} \\
            & \| n_2 \delta c\|_{L^{4}_T(L^\infty)} \lesssim \| n_2 \|_{L^\infty_T(L^2)} \| \delta c\|_{L^4_T(L^\infty)} \leq \| n_2 \|_{L^\infty_T(L^2)} \| \delta c\|_{L^\infty_T(L^2)}^{1/2} \| \delta c\|_{L^2_T(\dot B^1_{\infty, 1})}^{1/2}.
        \end{split}
    \end{equation*}
    Interpolation (with Proposition \ref{p:Interpolation}) therefore gives an appropriate estimate
    \begin{equation}\label{eq:UnCrEQ2}
        \| n_2 \delta c \|_{ \widetilde{L^{2}_T}(\dot B^{-1}_{\infty, 1})} \lesssim \| n_2 \delta c \|_{ \widetilde{L^{8}_T}(\dot B^{-1}_{\infty, 1})} \lesssim \| n_2 \|_{L^\infty_T(L^2)} \| \delta c \|_{L^\infty_T(L^2)}^{3/4} \| \delta c \|_{L^2_T(\dot B^1_{\infty, 1})}^{1/4}.
    \end{equation}
    Next, for the term $c_1 \delta n$, we recall the following embedding inequalities from the proof of Theorem \ref{t:Stability}, as well as write a new one with Proposition \ref{p:ComplexBesovInterpolate}. We have, for all $0 < s < 1$,
    \begin{equation*}
        \begin{split}
            & \| \delta n \|_{L^{2/s}_T(L^2)} \lesssim \| \delta n \|_{L^\infty_T(\dot H^{-s})}^{1-s} \| \delta n \|_{L^2_T(\dot H^{1-s})}^s \\
            & \| \delta n \|_{L^2_T(L^{2/s})} \lesssim \| \delta n \|_{L^2_T(\dot H^{1-s})} \\
            & \| c_1 \|_{L^{2/(1 - s)}_T(L^{2/s})} \lesssim \| c_1 \|_{L^{2/(1 - s)}_T(\dot H^{1-s})} \lesssim \| c_1 \|_{L^\infty_T(L^2)}^s \| c_1 \|_{L^2_T(\dot H^1)}^{1-s} \\
            & \| c_1 \|_{L^\infty_T(L^{4/(1-s)})} \lesssim \| c_1 \|_{L^\infty_T(L^2)}^{\frac{1-s}{2}} \| c_1 \|_{L^\infty_T(\dot B^0_{\infty, 2})}^{\frac{1+s}{2}}.
        \end{split}
    \end{equation*}
    This provides the following bounds for the product $c_1 \delta n$,
    \begin{equation*}
        \begin{split}
            & \| c_1 \delta n \|_{L^2_T(L^{2/(s+1)})} \lesssim \| c_1 \|_{L^{2/(1-s)}_T(L^{2/s})} \| \delta n \|_{L^{2/s}_T(L^2)} \lesssim \| c_1 \|_{L^{2/(1-s)}_T(L^{2/s})} \| \delta n \|_{L^\infty_T(\dot H^{-s})}^{1-s} \| \delta n \|_{L^2_T(\dot H^{1-s})}^s  \\
            & \| c_1 \delta n \|_{L^2_T(L^{4/(1+s)})} \lesssim \| c_1 \|_{L^\infty_T(L^{4/(1-s)})} \| \delta n \|_{L^2_T(L^{2/s})} \lesssim \| c_1 \|_{L^\infty_T(L^{4/(1-s)})} \| \delta n \|_{L^2_T(\dot H^{1-s})}.
        \end{split}
    \end{equation*}
    Therefore, by interpolating with the help of Proposition \ref{p:Interpolation}, we have a bound in the right space
    \begin{equation}\label{eq:UnCrEQ4}
        \begin{split}
            \| c_1 \delta n \|_{\widetilde{L^2_T}(\dot B^{-1}_{\infty, 1})} & \lesssim \| c_1 \delta n \|_{L^2_T(\dot B^{-(s+1)}_{\infty, \infty})}^{\frac{1-s}{1+s}} \| c_1 \delta n \|_{L^2_T(\dot B^{- (s+1)/2}_{\infty, \infty})}^\frac{2s}{1+s} \\
            & \lesssim \| c_1 \|_{L^\infty_T(L^2 \cap \dot B^0_{\infty, 2})} \| \delta n \|_{L^\infty_T(\dot H^{-s})}^{\frac{(1-s)^2}{1+s}} \| \delta n \|_{L^2_T(\dot H^{1-s})}^{\frac{3s - s^2}{1+s}} \\
            & \lesssim \| c_1 \|_{L^\infty_T(L^2 \cap \dot B^0_{\infty, 2})} \| \delta n \|_{L^\infty_T(\dot H^{-s}) \cap L^2_T(\dot H^{1-s})}.
        \end{split}
    \end{equation}
    Finally, we have to bound the term involving the velocity difference function $\delta u \cdot \nabla c_2$, which creates two kind of difficulties. The first one is related to criticality: it is clear from Yudovich's argument that it is the $L^\infty_T(L^2)$ norm of $\delta u$ that we want to bound. However, by using this bound we cannot bound the product $\delta u \cdot \nabla c_2$ in a space that is better than $L^2_T(L^2) \not \subset \widetilde{L^2_T}(\dot B^{-1}_{\infty, 1})$, and this is insufficient for using Besov maximal regularity (Proposition \ref{p:heat}). The second difficulty is related ot the nature of the spaces we handle. Even if we had a $\delta u$ bound in the critical space $L^\infty_T(\dot B^0_{2, 1})$ (in addition to the bound $c_2 \in L^2_T(\dot B^0_{\infty, 1})$ which is already known), we could not estimate the product $\delta u \cdot \nabla c_2$ in a Chemin-Lerner space, as we need. Both problems will be dealt with by resorting to logarithmic interpolation. We remark that we have the following basic inequalities
    \begin{equation}\label{eq:UniqueEQ6}
        \begin{split}
            & \| \delta u \cdot \nabla c_2 \|_{L^2_T(\dot B^{-2}_{\infty, \infty})} \lesssim \| \delta u \cdot \nabla c_2 \|_{L^2_T(L^2)} \lesssim \| \delta u \|_{L^\infty_T(L^2)} \| \nabla c_2 \|_{L^2_T(L^2)}  \\
            & \| \delta u \cdot \nabla c_2 \|_{L^2_T(\dot B^{-2/3}_{\infty, \infty})} \lesssim \| \delta u \|_{L^\infty_T(L^3)} \| \nabla c_2 \|_{L^2_T(L^\infty)} \\
            & \qquad \qquad \qquad \qquad \qquad \lesssim \big( \| \omega_1 \|_{L^\infty_T(L^{6/5})} + \| \omega_2 \|_{L^\infty_T(L^{6/5})} \big) \| \nabla c_2 \|_{L^2_T(L^\infty)} \\
            & \| \delta u \cdot \nabla c_2 \|_{L^2_T(\dot B^{-1}_{\infty, \infty})} \lesssim \| \delta u \cdot \nabla c_2 \|_{L^2_T(L^2)} \lesssim \| \delta u \|_{L^\infty_T(L^2)} \| \nabla c_2  \|_{L^2_T(L^\infty)}.
        \end{split}
    \end{equation}
    The last inequality is critical, but insufficient, whereas the first two inequalities give better bounds for the low and high frequencies respectively. We examine high frequencies first. Recall the high-frequency localizer $\dot S_0$ from the appendix. We have, for any non-negative $N \in \N$, by invoking the last two inequalities in \eqref{eq:UniqueEQ6},
    \begin{equation*}
        \begin{split}
            \big\| \dot S_0 \big( \delta u \cdot \nabla c_2 \big) \big\|_{\widetilde{L^2_T}(\dot B^{-1}_{\infty, 1})} & = \sum_{j = 0}^N 2^{-j} \big\| \dot \Delta_j \big( \delta u \cdot \nabla c_2 \big) \big\|_{L^2_T(L^\infty)} + \sum_{j = N+1}^\infty 2^{-j} \big\| \dot \Delta_j \big( \delta u \cdot \nabla c_2 \big) \big\|_{L^2_T(L^\infty)} \\
            & \lesssim N \,  \| \delta u \cdot \nabla c_2 \|_{L^2_T(\dot B^{-1}_{\infty, \infty})} + \sum_{j = N+1}^\infty 2^{-j/3} \| \delta u \cdot \nabla c_2 \|_{L^2_T(\dot B^{-2/3}_{\infty, \infty})} \\
            & \lesssim N \, \| \delta u \|_{L^\infty_T(L^2)} \| \nabla c_2 \|_{L^2_T(L^\infty)} + 2^{-N/3} K(T) \| \nabla c_2 \|_{L^2_T(L^\infty)}.
        \end{split}
    \end{equation*}
    We optimize the choice of $N \in \N$ by setting
    \begin{equation*}
        N \approx \log \left( 1 + \frac{K(T)}{\| \delta u \|_{L^\infty_T(L^2)}} \right).
    \end{equation*}
    This gives the high-frequency logarithmic interpolation inequality
    \begin{equation}\label{eq:UniqueEQ7}
        \big\| \dot S_0 ( \delta u \cdot \nabla c_2 ) \big\|_{\widetilde{L^2_T}(\dot B^{-1}_{\infty, 1})} \lesssim \| \delta u \|_{L^\infty_T(L^2)} \| \nabla c_2 \|_{L^2_T(L^\infty)} \log \left( 1 +  \frac{K(T)}{\| \delta u \|_{L^\infty_T(L^2)}} \right).
    \end{equation}
    We may proceed exactly in the same way for low frequencies, with the first and third inequalities in \eqref{eq:UniqueEQ6}. This gives, for $N \in \N$,
    \begin{equation*}
        \begin{split}
            \big\| ({\rm Id}  - \dot S_0) (\delta u \cdot \nabla c_2) \big\|_{\widetilde{L^2_T}(\dot B^{-1}_{\infty, 1})} & \lesssim N \| \delta u \cdot \nabla c_2 \|_{L^2_T(\dot B^{-1}_{\infty, \infty})} + 2^{-N} \| \delta u \cdot \nabla c_2 \|_{L^2_T(\dot B^{-2}_{\infty, \infty})} \\
            & \lesssim \| \delta u \|_{L^\infty_T(L^2)} \Big( N \| \nabla c_2 \|_{L^2_T(L^\infty)} + 2^{-N} \| \nabla c_2 \|_{L^2_T(L^2)} \Big).
        \end{split}
    \end{equation*}
    By fixing again $N \in \N$ such that
    \begin{equation*}
        N \approx \log \left( 1 + \frac{\| c_2 \|_{L^2_T(\dot H^1)}}{\| c_2 \|_{L^2_T(\dot B^1_{\infty, 1})}} \right),
    \end{equation*}
    we have a low-frequency logarithmic inequality
    \begin{equation}\label{eq:UniqueEQ8}
        \big\| ({\rm Id}  - \dot S_0) (\delta u \cdot \nabla c_2) \big\|_{\widetilde{L^2_T}(\dot B^{-1}_{\infty, 1})} \lesssim \| \delta u \|_{L^\infty_T(L^2)} \| \nabla c_2 \|_{L^2_T(L^\infty)} \log \left( 1 + \frac{\| c_2 \|_{L^2_T(\dot H^1)}}{\| c_2 \|_{L^2_T(\dot B^1_{\infty, 1})}} \right).
    \end{equation}
    Finally, by putting \eqref{eq:UniqueEQ7} and \eqref{eq:UniqueEQ8} together, we have an overall Chemin-Lerner bound for $\delta u \cdot \nabla c_2$, which involves an expression of the velocity difference function through an Osgood modulus of continuity, that is
    \begin{equation}\label{eq:UniqueEQ10}
        \| \delta u \cdot \nabla c_2 \|_{\widetilde{L^2_T}(\dot B^{-1}_{\infty, 1})} \lesssim \| \delta u \|_{L^\infty_T(L^2)} \| \nabla c_2 \|_{L^2_T(L^\infty)} \log \left( 1 + \frac{ K(T) \| c_2 \|_{L^2_T(\dot H^1)} }{ \| \delta u \|_{L^\infty_T(L^2)} \| c_2 \|_{L^2_T(\dot B^1_{\infty, 1})}} \right).
    \end{equation}
    
    Putting all the $\delta c$ estimates \eqref{eq:UnCrEQ1}, \eqref{eq:UnCrEQ3}, \eqref{eq:UnCrEQ2}, \eqref{eq:UnCrEQ4} along with the logarithmic interpolation inequality \eqref{eq:UniqueEQ10} immediately above, we have 
    \begin{equation*}
        \begin{split}
            \| \delta c \|_{L^2_T(\dot B^1_{\infty, 1})} & \lesssim \| \delta c_0 \|_{\dot B^0_{\infty, 2}} \\
            & \qquad+ \| \omega_1 \|_{L^\infty_T(L^{6/5})} \| \delta c \|_{L^2_T(\dot H^1)}^{1/3} \| \delta c \|_{L^2_T(\dot B^1_{\infty, 1})}^{2/3} \\
             &\qquad \qquad  + \| n_2 \|_{L^\infty_T(L^2)} \| \delta c \|_{L^\infty_T(L^2)}^{3/4} \| \delta c \|_{L^2_T(\dot B^1_{\infty, 1})}^{1/4} \\
              & \qquad \qquad \qquad + \| c_1 \|_{L^\infty_T(L^2 \cap \dot B^0_{\infty, 2})} \| \delta n \|_{L^\infty_T(\dot H^{-s}) \cap L^2_T(\dot H^{1-s})} \\
             & \qquad \qquad \qquad \qquad + \| \delta u \|_{L^\infty_T(L^2)} \| \nabla c_2 \|_{L^2_T(L^\infty)} \log \left( 1 + \frac{ K(T) \| c_2 \|_{L^2_T(\dot H^1)} }{ \| \delta u \|_{L^\infty_T(L^2)} \| c_2 \|_{L^2_T(\dot B^1_{\infty, 1})}} \right)
        \end{split}
    \end{equation*}
    Then, using Young's inequality to absorb the norm $\| \delta c \|_{L^2_T(\dot B^1_{\infty, 1})}$ in the left-hand side, we have the critical bound we were looking for:
    \begin{equation}\label{eq:UnCrEQ5}
        \begin{split}
            \| \delta c \|_{L^2_T(\dot B^1_{\infty, 1})} & \lesssim \| \delta c_0  \|_{\dot B^0_{\infty, 2}} + \| \omega_1 \|_{L^\infty_T(L^{6/5})}^3 \| \delta c \|_{L^2_T(\dot H^1)} \\
            & \qquad+ \| n_2 \|_{L^\infty_T(L^2)}^{4/3} \| \delta c \|_{L^\infty_T(L^2)} + \| c_1 \|_{L^\infty_T(L^2 \cap \dot B^0_{\infty, 2})} \| \delta n \|_{L^\infty_T(\dot H^{-s}) \cap L^2_T(\dot H^{1-s})} \\ 
            & \qquad \qquad + \| \delta u \|_{L^\infty_T(L^2)} \| \nabla c_2 \|_{L^2_T(L^\infty)} \log \left( 1 + \frac{ K(T) \| c_2 \|_{L^2_T(\dot H^1)} }{ \| \delta u \|_{L^\infty_T(L^2)} \| c_2 \|_{L^2_T(\dot B^1_{\infty, 1})}} \right)
        \end{split}
    \end{equation}
    Remark that the terms in the left-hand side are already bounded by the computations in Theorem \ref{t:Stability}.

    \medskip

    \textbf{STEP 3: Conclusion.} We are ready to conclude. We plug inequality \eqref{eq:UnCrEQ7} into the energy estimate \eqref{eq:UnCrEQ6} for $\delta n$. This provides, after integration over the time interval $[0, T]$,
    \begin{equation*}
        \begin{split}
            \| \delta n \|_{L^\infty_T(L^2)}^2 + \| \delta n \|_{L^2_T(\dot H^1)}^2 & \lesssim \| \delta n_0 \|_{L^2}^2 + \int_0^T \| \delta n \|_{L^2}^2 \Big( \| \nabla c_2 \|_{L^\infty}^2 + \| \nabla n_2 \|_{L^\infty} \Big) \dt \\
            & \qquad + \| \nabla n_2 \|_{L^1_T(L^\infty)} \| \delta u \|_{L^\infty_T(L^2)} + \| n_1 \|_{L^2_T(L^\infty)} \| \nabla \delta c \|_{L^2_T(L^\infty)}.
        \end{split}
    \end{equation*}
    From this point on, as in the proof of Theorem \ref{t:Stability}, we will call $M_1(t)$ any generic $L^1([0, T])$ function which depends only on $T$ and the solutions $(n_j, c_j, u_j)$, and treat any norms of the solutions as constants. Plugging \eqref{eq:UnCrEQ5} in the above gives
    \begin{equation*}
        \begin{split}
            \| \delta n \|_{L^\infty_T(L^2)}^2 + \| \delta n \|_{L^2_T(\dot H^1)}^2 & \lesssim \| \delta n_0 \|_{L^2}^2 + \| \delta c_0 \|_{\dot B^0_{\infty, 2}} \\
            & \qquad + \int_0^T M_1(t) \| \delta n \|_{L^2}^2 \dt \\
            & \qquad \qquad + \| \delta c \|_{L^\infty_T(L^2) \cap L^2_T(\dot H^1)} + \| \delta n \|_{L^\infty_T(\dot H^{-s}) \cap L^2_T(\dot H^{1-s})} \\
            & \qquad \qquad \qquad + \| \delta u \|_{L^\infty_T(L^2)} \log \left( 1 + \frac{K(T)}{\| \delta u \|_{L^\infty_T(L^2)}} \right).
        \end{split}
    \end{equation*}
    Applying Grönwall's lemma to the above gives a bound of the desired type. Set $F(r) = r \log(1 + \frac{K}{r})$ for convenience. We have
    \begin{equation*}
        \begin{split}
            \| \delta n \|_{L^\infty_T(L^2)}^2 + \| \delta n \|_{L^2_T(\dot H^1)}^2 & \lesssim \| \delta n_0 \|_{L^2}^2 + \| \delta c_0 \|_{\dot B^0_{\infty, 2}} \\
             & \qquad + F \left( \left[ \| \delta c_0 \|_{L^2} + \| \delta c_0 \|_{\dot B^0_{\frac{2}{s}, 2}} + \| \delta n_0 \|_{\dot H^{-s}} + \| \delta u_0 \|_{L^2} \right]^{e^{-KT}} \right),
        \end{split}
    \end{equation*}
    which, together with \eqref{eq:UnCrEQ5}, is the critical stability bound we were looking for.
\end{proof}

\appendix

\section{Appendix: Some Useful Results in Harmonic Analysis}\label{A}

In this section, we provide, for reference and ease of reading, some of the harmonic analysis results that are used throughout the text. Most of them are extremely standard, although some are less so: refined embeddings of Besov spaces, sharp maximal regularity properties and interpolation.

\subsection{Littlewood-Paley Analysis and Besov Spaces}

We recall the (homogeneous) Littlewood-Paley decomposition of functions. Consider a truncation function $\chi \in C^\infty_c(\R^d)$ which is radially symmetric, bounded by $0 \leq \chi \leq 1$, is such that $r \mapsto \chi(re)$ is non-increasing for $e \in \R^d$. We define a dyadic partition of the frequency space $1 = \sum_{j \in \Z} \varphi(2^{-j} \xi)$ with the function $\varphi$ constructed from the truncation function by $\varphi(\xi) = \chi(\xi) - \chi(2 \xi)$, so that $\varphi$ is supported in an annulus. We can write the associated \textit{formal} Fourier decomposition of a function
\begin{equation}\label{aeq:LP}
    \forall f, \qquad f = \sum_{j \in \Z} \dot \Delta_j f,
\end{equation}
where $\dot \Delta_j = \varphi(2^{-j}D)$ is the Fourier multiplier defined by $\what{\dot \Delta_j f}(\xi) = \varphi(2^{-k} \xi) \what{f}(\xi)$. We attract the reader's attention to the fact that the decomposition \eqref{aeq:LP} is formal. In fact, since the frequency partition $1 = \sum_{j \in \Z} \varphi(2^{-j}\xi)$ fails for $\xi = 0$, identity \eqref{aeq:LP} cannot be true for all $f \in \mc S'$. It fails for all polynomials $P \in \R[X]$, as well as a broad class of functions \cite{Cobb2024}. However, we may eliminate this low-frequency issue by working in the quotient $\mc S' / \, \R[X]$ of distributions modulo polynomials. In this class, the decomposition \eqref{aeq:LP} holds, and the series therein converges in the $\mc S' / \, \R[X]$ topology \cite{Sawano}.

\medskip

Note that when $p \in [1, + \infty[$, the spaces $L^p$ and $L^p / \,\R[X]$ are naturally isomorphic because $L^p$ contains no non-trivial polynomial. We will therefore allow a slight abuse of notation and omit the quotient when there is no ambiguity. 


\medskip

Also note that since $\dot \Delta_j$ are convolution operators whose kernels have a $L^1$ norm independent of $j \in \Z$, they are bounded in the $L^p \longrightarrow L^p$ topologies independently of $j \in \Z$ and $1 \leq p \leq + \infty$. The same applies for the operators Fourier truncation operators
\begin{equation}\label{aeq:FTruncation}
    \dot S_j := \sum_{m \leq j-1} \dot \Delta_m \qquad \text{and} \qquad S_j = \chi(2^{-j}D)
\end{equation}
although in the case $p = +\infty$, one has to quotient out constants for the homogeneous operator $\dot S_j : L^\infty / \, \R \longrightarrow L^\infty / \, \R$. Remark that, contrary to the non-homogeneous operator in \cite{BCD}, the Fourier truncation operator $S_j$ is non-zero for all $j \in \Z$.

\medskip

One of the major interests of the Littlewood-Paley decomposition is the behavior of the frequency localizations $\dot \Delta_j f$ under differentiation, and more generally, under the operation  of Fourier multipliers. Since the Fourier transform of $\dot \Delta_j f$ is supported in an annulus of radius $2^j$, differentiation $\nabla$ essentially behaves as a multiplication by $2^j$.

\begin{prop}[Bernstein inequalities -- Lemma 2.1 p. 52 in \cite{BCD}]\label{p:Bern}
    Consider $1 \leq p \leq q \leq + \infty$ and $f \in \mc S'$. Then, for any $j \in \Z$ the following assertions are true;
    \begin{enumerate}[(i)]
        \item For any $k \in \N$, we have
        \begin{equation*}
            \| \nabla^k \dot \Delta_j f \|_{L^p} \approx 2^{kj} \| \dot \Delta_j f \|_{L^p},
        \end{equation*}
        where the implicit constants depend only on $\varphi$ and $k$.

        \item We also have
        \begin{equation*}
            \| \dot \Delta_j f \|_{L^q} \lesssim 2^{jd \left( \frac{1}{p} - \frac{1}{q} \right)} \| \dot \Delta_j f \|_{L^p}.
        \end{equation*}
    \end{enumerate}
\end{prop}

\begin{prop}\label{p:FourierMultipliers}
    Let $\sigma(\xi) \in C^{2 \lfloor 1 + d/2 \rfloor}(\R^d \setminus \{ 0 \})$ be function such that there is a $m \in \R$, with, for all $k \leq 2 \lfloor 1 + d/2 \rfloor$, 
    \begin{equation*}
        |\nabla^k \sigma (\xi)| \leq C(k) |\xi|^{m - \alpha}.
    \end{equation*}
    Then, for all $j \in \Z$, $p \in [1, + \infty]$ and $f \in \mc S'$,
    \begin{equation*}
        \| \dot \Delta_j \sigma(D) f \|_{L^p} \lesssim 2^{jm} \| \dot \Delta_j f \|_{L^p}.
    \end{equation*}
\end{prop}

The Littlewood-Paley decomposition also allows to naturally define an appropriate family of function spaces, (homogeneous) Besov spaces.

\begin{defi}
    Consider $s \in \R$ and $p, r \in [1, +\infty]$. We define the homogeneous Besov space $\dot B^s_{p, r}$ as the set of (classes of) functions $f \in \mc S' / \, \R[X]$ such that
    \begin{equation*}
        \| f \|_{\dot B^s_{p, r}} := \Big\| \big( 2^{js} \| \dot \Delta_j f \|_{L^p} \big)_{j \in \Z} \Big\|_{\ell^r(\Z)} < + \infty \,.
    \end{equation*}
    This space is Banach for all values of $(s, p, r)$ (see \cite{Sawano}).
\end{defi}

Given the Littlewood-Paley decomposition, one can write embeddings of Besov spaces. For example, for $1 \leq p < + \infty$,
\begin{equation*}
    \dot B^0_{p, 1} \hookrightarrow L^p \hookrightarrow \dot B^0_{p, \infty}.
\end{equation*}
The same holds for $p = +\infty$, taking care to quotient constant functions. We also have, thanks to the first Bernstein inequality, for $s \in \R$ and $1 \leq p \leq q \leq + \infty$,
\begin{equation}\label{p:BesovEmbed}
    \dot B^s_{p, r} \hookrightarrow \dot B^{s - d \left( \frac{1}{p} - \frac{1}{q} \right)}_{q, r}.
\end{equation}
However, much more precise embeddings are available, and are needed in our analysis.

\begin{prop}[Theorems 2.40 and 2.41 in \cite{BCD}]\label{p:BesovFine}
    The following assertions are true:
    \begin{enumerate}[(i)]
        \item For $p \in [2, + \infty[$ we have $\dot B^0_{p, 2} \subset L^p$, and for $p \in [1, 2]$ we have $\dot B^0_{p, p} \subset L^p$.
        \item For $p \in [2, + \infty]$ we have $L^p \subset \dot B^0_{p, p}$, and for $p \in ]1, 2]$ we have $L^p \subset \dot B^0_{p, 2}$.
    \end{enumerate}
\end{prop}

\begin{rmk}
    It is important to note that the endpoint embeddings above do not hold. Notably $\dot B^0_{\infty, 2}$ does \textit{not} embed in $L^\infty$ (or even in $L^\infty / \, \R[X]$). We instead have $\dot B^0_{\infty, 2} \hookrightarrow {\rm BMO}$. Likewise, $L^1$ does not embed in $\dot B^0_{1, 2}$. Instead, we have a comparable embedding for the Hardy space $\dot F^0_{1, 2} \approx \mc H^1 \hookrightarrow \dot B^0_{1, 2}$. For $L^1$, the best Besov embedding available is the trivial one $L^1 \hookrightarrow \dot B^0_{1, \infty}$. We refer to 2.2.2 p. 105 in \cite{GrafakosM}.
\end{rmk}

\subsection{Paradifferential Calculus}\label{ass:Bony}

This paragraph centers around the Bony decomposition. These are now standard tools, which are presented in Chapter 2 of \cite{BCD}. Here, we give a slight generalization of Theorems 2.47 and 2.52 from \cite{BCD}, as well as a presentation of the paraproducts and remainder operators accounting for our definition of homogeneous spaces (as sets of distributions modulo polynomials).

\medskip

Recall the Bony decomposition of a product $fg = \mc {\dot T}_f(g) + \mc {\dot T}_g(f) + \mc {\dot R}(f, g)$ after \cite{Bony}, where the paraproduct operator is (formally) defined by
\begin{equation}\label{aeq:Paraproduct}
    \mc {\dot T}_f(g) := \sum_{m \in \Z} S_{m-1} f \, \dot \Delta_m g,
\end{equation}
and the remainder operator is defined by
\begin{equation*}
    \mc {\dot R}(f, g) := \sum_{|m - n| \leq 1} \dot \Delta_m f \, \dot \Delta_n g \, .
\end{equation*}
We attract the reader's attention on the fact that the paraproduct is defined with the help of the \textit{non-homogeneous} operator $S_{m-1}$ from \eqref{aeq:FTruncation}.

\begin{prop}\label{p:Bony}
    Let $1 \leq p, a, b, r \leq + \infty$ with $\frac{1}{a} + \frac{1}{b} = \frac{1}{p}$. The following assertions hold.
    \begin{enumerate}[(i)]
        \item Assume that $t > 0$ and $s \in \R$. For all functions $(f, g) \in \dot B^{-t}_{a, \infty} \times \dot B^s_{b, r}$ such that $S_j f \in L^a$ for some $j \in \Z$, we have the following paraproduct bound:
        \begin{equation*}
            \big\| \mc {\dot T}_f(g) \big\|_{\dot B^{s-t}_{p, r}} \lesssim \frac{1}{t} \| f \|_{\dot B^{-t}_{a, \infty}} \| g \|_{\dot B^s_{b, r}}.
        \end{equation*}

        \item Assume that $(f, g) \in L^a \times B^s_{b, r}$. Then we have the following paraproduct bound
        \begin{equation*}
            \big\| \mc {\dot T}_f(g) \big\|_{\dot B^{s}_{p, r}} \lesssim \| f \|_{L^a} \| g \|_{\dot B^s_{b, r}}.
        \end{equation*}

        \item Assume that $(f, g) \in \dot B^s_{a, r_1} \times \dot B^t_{b, r_2}$ with $\frac{1}{r} \leq \frac{1}{r_1} + \frac{1}{r_2}$ and $s + t > 0$. Then
        \begin{equation*}
            \big\| \mc {\dot R}(f, g) \big\|_{\dot B^{s+t}_{p, r}} \lesssim \| f \|_{\dot B^s_{a, r_1}} \| g \|_{\dot B^t_{b, r_2}}.
        \end{equation*}
        In the case where $s+t = 0$ and $1 \leq \frac{1}{r_1} + \frac{1}{r_2}$, we have
        \begin{equation*}
            \big\| \mc {\dot R}(f, g) \big\|_{\dot B^0_{p, \infty}} \lesssim \| f \|_{\dot B^s_{a, r_1}} \| g \|_{\dot B^t_{b, r_2}}.
        \end{equation*}
    \end{enumerate}
\end{prop}

\begin{rmk}
    In the previous proposition, we see that the paraproduct $\mc {\dot T}_f(g)$ is defined as an element of a homogeneous space (and hence an element of $\mc S' / \, \R[X]$) provided that the subscript $f$ (or some regularized function $S_j f$) lies in some Lebesgue space. In particular, $g$ may be defined only up to polynomials, and we do not require the Littlewood-Paley decomposition of $f$ to converge to $f$ (that is, we do not need $f \in \mc S'_h$).
\end{rmk}

\begin{rmk}
    On the other hand, the remainder $\mc {\dot R}(f, g)$ is defined as an element of a homogeneous Besov space as long as the regularity requirements of $(iii)$ are fulfilled. These do not specify that either $f$ or $g$ need to be actual distributions, they can be distributions modulo polynomials.
\end{rmk}

\subsection{Maximal Regularity Results for the Heat Equation}

This paragraph is centered around estimates for the solution of the heat equation. For some given $f$ and $u_0$ (to be specified), we consider the solution $u$ of the following initial value problem:
\begin{equation}\label{aeq:heat}
    \begin{cases}
        \big( \partial_t - \Delta \big) u = f\\
        u(0) = u_0.
    \end{cases}
\end{equation}
Formally, the solution of this problem is given by Duhamel's formula
\begin{equation*}
    u = e^{t \Delta} u_0 + \int_0^t e^{(t-s) \Delta} f(s) \ds,
\end{equation*}
and we wish to state optimal regularity results for $u$ for given regularities of $f$ and $u_0$. One of the best known theorems is related to the fact that the symbol of the second summand in the equation above reads
\begin{equation*}
    \int_0^t e^{-(t-s)|\xi|^2} \ds = \frac{1}{|\xi|^2} e^{-t|\xi|^2},
\end{equation*}
which is the product of a smooth factor with a homogeneous function of degree $-2$. As we have seen in Proposition \ref{p:FourierMultipliers}, homogeneous symbols play an important role in harmonic analysis, and have been intensively studied. For example, we have the famous Calder\'on-Zygmund theorem. We refer to \cite{SteinRed} for a proof and a much more general statement.

\begin{prop}\label{p:CZ}
    Consider $m \in C^\infty(\R^d \setminus \{ 0 \})$ to be a bounded homogeneous function of degree zero. Then, for $1 < p < + \infty$, the operator $m(D) : L^p \longrightarrow L^p$ is bounded and the inequality
    \begin{equation*}
        \| m(D) f \|_{L^p} \lesssim \frac{p^2}{p-1} \| f \|_{L^p}
    \end{equation*}
    holds for all $f \in L^p$.
\end{prop}

This result indicates that the operator $\nabla^2 \int_0^t e^{(t-s)\Delta} \ds$ may have the same properties (while the time aspect creates some new difficulties). This is in essence what the following very well-known statement says. We refer to the pedagogical article \cite{M2024} for a summary and examples of applications.

\begin{prop}\label{p:MaxRegularity}
    Consider $1 < p, q < + \infty$ and $f \in L^p(L^q)$ and assume $u_0 = 0$. Then, the solution of problem \eqref{aeq:heat} satisfies the bounds
    \begin{equation*}
        \| \partial_t u \|_{L^p(L^q)} + \| \nabla^2 u \|_{L^p(L^q)} \lesssim \| f \|_{L^p(L^q)}.
    \end{equation*}
\end{prop}

The obvious drawback of this estimate is it's failure for $p, q \in \{ 1, + \infty\}$. As for Calder\'on-Zygmund operators, one rather has ${\rm BMO}$ estimates instead of $L^\infty$ ones \cite{OS2021, DWY}. Another solution \cite{Chemin1999} is to replace Lebesgue spaces by Besov spaces. The procedure involves applying the $\dot \Delta_j$ operators to the heat equation, and taking the $L^p(L^q)$ norm. We then have (see Corollary 2.5 in \cite{BCD}),
\begin{equation*}
    \| \dot \Delta_j u \|_{L^p(L^q)} \lesssim 2^{-2j/p} \| \dot \Delta_j u_0 \|_{L^q} + 2^{-2j} \| \dot \Delta_j f \|_{L^p(L^q)}.
\end{equation*}
The Bernstein inequalities (Proposition \ref{p:Bern}) show that the solution $u$ does indeed gain two derivatives when compared to $f$. The issue is that, if we were to sum with respect to $j \in \Z$ and try to obtain a Besov norm, we would not get spaces of the type $L^p(\dot B^s_{q, r})$, as the norm $L^p(L^q)$ would stay inside the $\ell^r$ norm. This motivates the introduction of a different family of spaces, Chemin-Lerner spaces after \cite{CL1995}.

\begin{defi}
    Consider $1 \leq p, q, r \leq + \infty$ and $s \in \R$. We define the norm 
    \begin{equation*}
        \| f \|_{\widetilde{L^p}(\dot B^s_{q, r})} := \Big\| \big( 2^{js} \| \dot \Delta_j f \|_{L^p(L^q)} \big)_{j \in \Z} \Big\|_{\ell^r(\Z)}.
    \end{equation*}
\end{defi}

These norms relate to usual Besov spaces \textit{via} the Minkowski inequality.

\begin{prop}[see Subsection 2.6 in \cite{BCD}]\label{p:Min}
    Consider $s \in \R$ and $1 \leq p, q, r \leq + \infty$. Then,
    \begin{equation*}
        \| f \|_{\widetilde{L^p}(\dot B^s_{p, r})} \leq \| f \|_{L^p(\dot B^s_{p, r})} \text{ if } r \geq q \qquad \text{and} \qquad \| f \|_{L^p(\dot B^s_{p, r})} \leq \| f \|_{\widetilde{L^p}(\dot B^s_{p, r})} \text{ if } r \leq q.
    \end{equation*}
\end{prop}

All in all, we have the following estimate.

\begin{prop}[Subsubsection 3.4.1. in \cite{BCD}]\label{p:heat}
    Consider a solution $u$ of \eqref{aeq:heat}. Then, for any $p, r \in [1, \infty]$ and $\rho \leq q$, the following estimate holds: for all times $T > 0$,
    \begin{equation*}
        \| u \|_{\widetilde{L^q_T}(\dot{B}^{s + 2/q}_{p, r})} \lesssim \| u_0 \|_{\dot{B}^s_{p, r}} + \| f \|_{\widetilde{L^\rho_T}(\dot{B}^{s - 2 + 2/\rho}_{p, r})}.
    \end{equation*}
\end{prop}

\subsection{Sharp Results on the Heat Semi-Group}

In this paragraph, we gather some sharp results about the heat equation. We start by the heat flow characterization of some Besov spaces. The following proposition is a direct application of Theorem 2.34 in \cite{BCD}.

\begin{prop}\label{p:BesovChar}
    We have the following equivalences of norms: for any $1 \leq p \leq + \infty$,
    \begin{equation*}
        \| f \|_{\dot B^0_{p, 2}}^2 \approx \int_0^{+ \infty} \| \nabla e^{t \Delta} f \|_{L^p}^2 \dt
    \end{equation*}
    and
    \begin{equation*}
        \| f \|_{\dot B^{-1}_{p, 1}} \approx \int_0^{+ \infty} \| \nabla e^{t \Delta} \|_{L^p} \dt.
    \end{equation*}
\end{prop}

\begin{proof}
    Theorem 2.34 from \cite{BCD} asserts that for all $s > 0$ and $1 \leq r \leq + \infty$, we have the norm equivalence
    \begin{equation*}
        \| f \|_{\dot B^{-2s}_{p, r}} \approx \bigg( \int_0^{+ \infty} \big\| e^{t \Delta} f\big\|_{L^p}^r \, t^{sr-1} \dt  \bigg)^{1/r}.
    \end{equation*}
    In particular, taking $(s, r) = (\frac{1}{2}, 2)$ and $(s, r) = (1, 1)$, and applying the inequalities above to $\| \nabla f \|_{\dot B^{- 2s}_{p, r}} \approx \| f \|_{\dot B^{1-2s}_{p, r}}$ gives the result.
\end{proof}

The proposition above shows that the Besov norm $f \in \dot B^0_{\infty, 2}$ is characterized by the heat condition $e^{t \Delta}f \in L^2(\dot W^{1, \infty})$. In fact, it is possible to be more precise, as is shown by the following result of Arsénio and Houamed \cite{AH2025}.

\begin{prop}[Proposition 2.2 in \cite{AH2025}]\label{p:AGHeat}
    Consider $q < + \infty$, $s \in \R$ and $p \in [1, + \infty]$. Then, for all $f \in \dot{B}^s_{p, q}$. Then we have the inequality $\| e^{t \Delta} f \|_{L^q(\dot B^{s + 2/q}_{p, 1})} \lesssim \| f \|_{B^s_{p, q}}$. In particular we have
    \begin{equation*}
        \| \nabla e^{t \Delta} f \|_{L^2(\dot B^0_{\infty, 1})} \lesssim \| f \|_{\dot B^0_{\infty, 2}}.
    \end{equation*}
\end{prop}


\subsection{Interpolation Inequalities}

To conclude this appendix, we provide some interpolation inequalities. We start with reminding the classical Gagliardo-Nirenberg estimate (Corollary 1.2 of \cite{CDGG}).

\begin{prop}[Theorem 2.44 in \cite{BCD}]\label{p:GN}
    Consider $2 \leq p < + \infty$ such that $\frac{1}{d} + \frac{1}{p} > \frac{1}{2}$. Then for all $f \in H^1$, one has
    \begin{equation*}
        \| f \|_{L^p} \leq C(p) \| f \|_{L^2}^\theta \| \nabla f \|_{L^2}^{1 - \theta} \qquad \text{where } \theta = \frac{d(p-2)}{2p}.
    \end{equation*}
\end{prop}

The next statement essentially the characterization of Besov spaces as real interpolation spaces. We refer to the book \cite{Tartar} for an introduction to interpolation theory, and in particular Chapter 32, which explains why 
\begin{equation*}
    \dot B^s_{p, r} = \big[ L^p, \dot W^{1, p} \big]_{s, r} \qquad \text{for } s \in ]0, 1[.
\end{equation*}
However, we will never need to do abstract interpolation. The following result is enough for our purposes. It's proof is essentially identical to that of Proposition \ref{p:Interpolation}, ignoring time dependency, so we will skip it.

\begin{prop}\label{p:BesovInterpolate}
    Consider $s_1 < s_2$ and $1 \leq p \leq + \infty$. Then, the following (real) interpolation inequality holds
    \begin{equation*}
        \| f \|_{\dot B^s_{p, 1}} \lesssim \| f \|_{\dot B^{s_1}_{p, \infty}}^\theta \| f \|_{\dot B^{s_2}_{p, \infty}}^{1 - \theta},
    \end{equation*}
    where the regularity exponent is given by $s = \theta s_1 + (1 - \theta)s_2$.
\end{prop}

We will also need the following inequality, which is related to complex interpolation (and simply follows from Hölder's inequality).

\begin{prop}\label{p:ComplexBesovInterpolate}
    Consider $p \in [2, + \infty[$. Then we have the embeddings $L^p \subset \dot B^0_{p, 2} \subset L^2 \cap \dot B^0_{\infty, 2}$. In particular, the following inequality holds: 
    \begin{equation*}
        \| f \|_{L^p} \lesssim \| f \|_{\dot B^0_{p, 2}} \leq \| f \|_{L^2}^{2/p} \| f \|_{\dot B^0_{\infty, 2}}^{1 - \frac{2}{p}}
    \end{equation*}
\end{prop}

\begin{proof}
    First of all, Hölder's inequality gives
    \begin{equation*}
        \sum_{j \in \Z} \| \dot \Delta_j f \|_{L^p}^2 \leq \sum_{j \in \Z} \| \dot \Delta_j f \|_{L^2}^{4/p} \| \dot \Delta_j f \|_{L^\infty}^{2 \left( 1 - 2/p \right)}.
    \end{equation*}
    We then apply Hölder's inequality in the sequence spaces $\ell^{2/p}(\Z)$ and $\ell^{p/(p-2)} (\Z)$, and this gives
    \begin{equation*}
        \sum_{j \in \Z} \| \dot \Delta_j f \|_{L^p}^2 \leq \Big( \sum_{j \in \Z} \| \dot \Delta_j f \|_{L^2}^2 \Big)^{2/p} \Big( \sum_{j \in \Z} \| \dot \Delta_j f \|_{L^2}^2 \Big)^{1 - 2/p}
        \end{equation*}
        which is the square of the second inequality. The embedding $\dot B^0_{p, 2}$ simply comes from Proposition \ref{p:BesovFine}.
\end{proof}

\begin{rmk}
    In Proposition \ref{p:ComplexBesovInterpolate}, if one only aims at $L^p$ bounds, then the spaces $L^2$ and $\dot B^0_{\infty, 2}$ can be replaced by $L^{2, \infty}$ and ${\rm BMO}$ respectively, although the proof is more involved \cite{MRR}. This is consistent with the embedding $\dot B^0_{\infty, 2} \subset {\rm BMO}$. However, as we have explained in the introduction, the space $\dot B^0_{\infty, 2}$ is both natural and optimal regarding our energy estimates, so we are content with the simple statement above.
\end{rmk}

We finish this appendix by addressing a common problem of Chemin-Lerner spaces: they are rather difficult to use in Grönwall-type arguments. A possible remedy to this has been given by Arsénio and Gallagher \cite{AG2020}, and consists in a maximal regularity property for the spaces $L^q(\dot B^s_{p, r})$ which holds for all $s \in \R$, all $r \leq q$ and $p \in [1, \infty]$, but is limited to $1 < q < + \infty$.

In the present article, we adopt a different approach (which is ultimately equivalent, up to some abstract interpolation). The Chemin-Lerner can in fact be seen as real interpolation spaces
\begin{equation}\label{aeq:CLinterpolation}
    \widetilde{L^q}(\dot B^{\theta s_1 + (1 - \theta s_2)}_{p, 1}) = \left[ \widetilde{L^q}(\dot B^{s_1}_{p, \infty}), \widetilde{L^q}(\dot B^{s_2}_{p, \infty}) \right]_{\theta, 1} \, .
\end{equation}
Note that $L^q(\dot B^\sigma_{p, \infty}) \subset \widetilde{L^q}(\dot B^\sigma_{p, \infty})$ by Proposition \ref{p:Min}. We owe these helpful remarks to Anatole Gaudin. In this article, we will not use abstract interpolation \eqref{aeq:CLinterpolation}, and instead prove the simpler real interpolation inequality below.

\begin{prop}\label{p:Interpolation}
    Consider $s_1 < s_2$ and $p, q \in [1, + \infty]$. For all $0 < \theta < 1$ set $s = \theta s_1 + (1 - \theta)s_2$. The following interpolation inequality holds: for all $f \in \widetilde{L^q_T}(\dot B^{s_1}_{p, \infty}) \cap \widetilde{L^q_T}(\dot B^{s_2}_{p, \infty})$,
    \begin{equation*}
        \| f \|_{\widetilde{L^q_T}(\dot B^s_{p, 1})} \lesssim \| f \|_{\widetilde{L^p_T}({\dot B^{s_1}_{p, \infty}})}^\theta \| f \|_{\widetilde{L^p_T}({\dot B^{s_2}_{p, \infty}})}^{1 - \theta} \leq \| f \|_{L^p_T{(\dot B^{s_1}_{p, \infty}})}^\theta \| f \|_{L^p_T({\dot B^{s_2}_{p, \infty}})}^{1 - \theta}
    \end{equation*}
\end{prop}

\begin{proof}
    The proof is a standard application of real interpolation techniques. Consider $N \in \Z$ to be fixed later on. We decompose the Chemin-Lerner norm $\widetilde{L^q_T}(\dot B^s_{p, 1})$ as
    \begin{equation*}
        \begin{split}
            \| f \|_{\widetilde{L^q_T}(\dot B^s_{p, 1})} & = \sum_{j \in \Z} 2^{js} \| \dot \Delta_j f \|_{L^q_T(L^p)} \\
            & \sum_{- \infty}^N 2^{js} \| \dot \Delta_j f \|_{L^q_T(L^p)} + \sum_{N+1}^\infty 2^{js} \| \dot \Delta_j f \|_{L^q_T(L^p)}.
        \end{split}
    \end{equation*}
    We now estimate low and high frequencies with the $\widetilde{L^q_T}(\dot B^{s_1}_{p, \infty})$ and $\widetilde{L^q_T}(\dot B^{s_2}_{p, \infty})$ norms respectively. This leads to 
    \begin{equation*}
        \begin{split}
            \| f \|_{\widetilde{L^q_T}(\dot B^s_{p, 1})} & \leq \| f \|_{\widetilde{L^q_T}(\dot B^{s_1}_{p, \infty})} \sum_{- \infty}^N 2^{j(s - s_1)} + \| f \|_{\widetilde{L^q_T}(\dot B^{s_2}_{p, \infty})} \sum_{- \infty}^N 2^{j(s - s_2)} \\
            & \lesssim 2^{N(s - s_1)} \| f \|_{\widetilde{L^q_T}(\dot B^{s_1}_{p, \infty})} + 2^{N(s - s_2)} \| f \|_{\widetilde{L^q_T}(\dot B^{s_2}_{p, \infty})}.
        \end{split}
    \end{equation*}
    Note that the two sums in the inequality above are convergent, due to the fact that $s_1 < s < s_2$, thanks to our assumptions. We now fix the value of $N \in \Z$ in order to make the two terms in the upper bound equal. in other words, we choose $N$ so that
    \begin{equation*}
        2^N \approx \left( \frac{\| f \|_{\widetilde{L^q_T}(\dot B^{s_2}_{p, \infty})}}{\| f \|_{\widetilde{L^q_T}(\dot B^{s_1}_{p, \infty})}} \right)^{\frac{1}{s_2 - s_1}},
    \end{equation*}
    which leads to 
    \begin{equation*}
        \| f \|_{\widetilde{L^q_T}(\dot B^{s}_{p, 1})} \lesssim \| f \|_{\widetilde{L^q_T}(\dot B^{s_1}_{p, \infty})}^{1 - \frac{s - s_1}{s_2 - s_1}} \| f \|_{\widetilde{L^q_T}(\dot B^{s_1}_{p, \infty})}^{\frac{s - s_1}{s_2 - s_1}}.
    \end{equation*}
    Noticing that $\theta = 1 - \frac{s - s_1}{s_2 - s_1}$ gives the first inequality. The second one is the Minkowski's inequality from Proposition \ref{p:Min}.
\end{proof}

\begin{rmk}
    We believe the computations in this article could have been done with the maximal regularity result of Arsénio and Gallagher \cite{AG2020}. However, Proposition \ref{p:AGHeat} uses simple bounds, and no maximal functions or abstract interpolation, although it certainly increases the volume of our computations.
\end{rmk}




\end{document}